# A Dynamical Key to the Riemann Hypothesis

Chris King Emeritus, University of Auckland 11 May – 16 Oct 2011

v1.48 ( http://arxiv.org/abs/1105.2103 )

**Abstract**: *We investigate a dynamical basis for the Riemann hypothesis[1] (RH) that the non-trivial zeros of the Riemann zeta function lie on the critical line $x = ½$. In the process we graphically explore, in as rich a way as possible, the diversity of zeta and L-functions, to look for examples at the boundary between those with zeros on the critical line and otherwise. The approach provides a dynamical basis for why the various forms of zeta and L-function have their non-trivial zeros on the critical line. It suggests RH is an additional unprovable postulate of the number system, similar to the axiom of choice, arising from the asymptotic behavior of the primes as $n \to \infty$.*

The images in the figures are generated using a Mac software research application developed by the author, which is available at: http://dhushara.com/DarkHeart/RZV/. It includes open source files for XCode compilation for flexible research use and scripts for the open source math packages PARI-GP and SAGE to generate *L*-functions of elliptic curves and modular forms.

The Riemann zeta function $\zeta(z) = \sum_{n=1}^{\infty} n^{-z} = \prod_{p \text{ prime}} \left(1 - p^{-z}\right)^{-1}$ Re($z$)>1 is defined as either a sum of complex exponentials over integers, or as a product over primes, due to Euler's prime sieving. The zeta function is a unitary example of a Dirichlet series $\sum_{n=1}^{\infty} a_n n^{-z}$, which are similar to power series except that the terms are complex exponentials of integers, rather than being integer powers of a complex variable as with power series. We shall examine a variety of Dirichlet series to discover which, like zeta, have their non-real zeros on the critical line $x = ½$ and which don't.

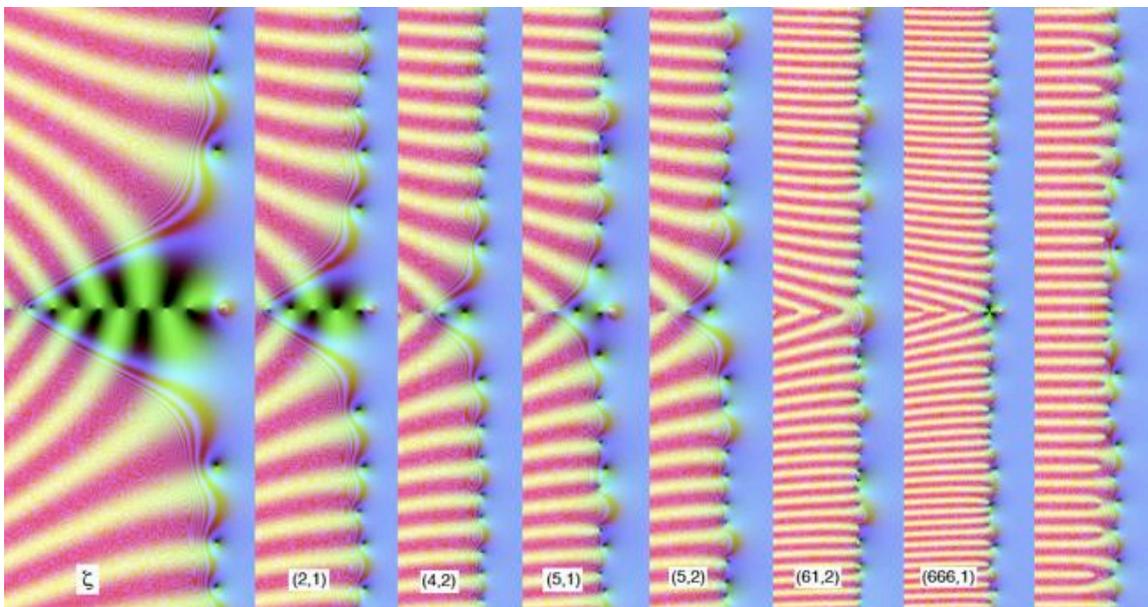

Fig 1: Riemann zeta and a selection of Dirichlet *L*-functions with a non-*L* function for comparison: $L(2,1)$ and $L(5,1)$ have regular zeros on $x = 0$ as well as non-trivial zeros on $x = ½$, due to their being equal to zeta with additional prime product terms. While $L(4,2)$ is symmetric with real coefficients, $L(5,2)$ and $L(61,2)$ have asymmetric non-trivial zeros on $x = ½$, having conjugate *L*-functions. $L(666,1)$ is similar to $L(2,1)$ and $L(5,1)$, but has a central third-order zero due to 666 being the product of three distinct primes $666=2.3^2.37$. Far right the period 10 non-*L*-function with $\chi = \{0,1,0,-1,0,0,0,1,0,-1\}$ (portrayed naked of any functional equation for 100 terms) has zeros in the critical strip $0<x<1$ manifestly varying from the critical line. Images generated using the author's application RZViewer for Mac (http://www.dhushara.com/DarkHeart/RZV/RZViewer.htm ).

---

[1] Links in blue which do not have an htm reference refer to Wikipedia search items, or in the case of Collatz conjecture, Mellin transform, prime counting and grand unitary ensemble to the live htm link in King 2009.

In historical terms, there is a unique class of such series, which do appear to have their unreal zeros on the critical line - the Dirichlet *L*-series, or when extended to the complex plane, L-functions:

$$L(z,\chi) = \sum_{n=1}^{\infty} \chi(n) n^{-z} = \prod_{p \text{ prime}} \left(1 - \chi(p) p^{-z}\right)^{-1} \text{ where } \chi(n), n = 0, \cdots, k-1 \text{ is a Dirichlet character.}$$

It was originally proven by Dirichlet that $L(1, \chi) \neq 0$ for all Dirichlet characters $\chi$, allowing him to establish his theorem on primes in arithmetic progressions. While $\zeta(1)$ is singular, $L(1,\chi)$ for non-trivial characters is known to be transcendental (Gun et. al.). For example

$$1 - \frac{1}{3} + \frac{1}{5} - \frac{1}{7} \cdots = L(1, \chi(4,2)) = \frac{\pi}{4},$$ leading to the study of special values of *L*-functions.

A Dirichlet character is any function $\chi$ from the integers to the complex numbers, such that:
1) **Periodic:** There exists a positive integer k such that $\chi(n) = \chi(n+k)$ for all *n*.
2) **Relative primality:** If $\gcd(n,k) > 1$ then $\chi(n) = 0$; if $\gcd(n,k) = 1$ then $\chi(n) \neq 0$.
3) **Completely multiplicative:** $\chi(mn) = \chi(m) \chi(n)$ for all integers m and n.

Consequently $\chi(1)=1$ and since only numbers relatively prime to *k* have non-zero characters, there are $\phi(k)$ of these where $\phi$ is the totient function, consisting of the number integers less than *n* coprime to n, and each non-zero character is a $\phi$-th complex root of unity. These conditions lead to the possible characters being determined by the finite commutative groups of units in the quotient ring *Z/kZ*, the residue class of an integer *n* being the set of all integers congruent to *n* modulo *k*.

As a consequence of the particular definition of each $\chi$, $L(z, \chi)$ is also expressible as a product over a set of primes $p_i$ with terms depending on the Dirichlet characters of $p_i$. As well as admitting an Euler product, oth Riemann zeta and the Dirichlet *L*-functions (*DL*-functions) also have a generic functional equation enabling them to be extended to the entire complex plane minus a simple infinity at *z* = 1 for the principal characters, whose non-zero terms are 1, as is the case of zeta.

Extending RH to the *L*-functions gives rise to the generalized Riemann hypothesis - that for all such functions, all zeros on the critical strip $0 < x < 1$ lie on $x = ½$. Examining where the functional boundaries lie, beyond which the unreal zeros depart from the critical line, has become one major avenue of attempting to prove or disprove RH, as noted in Brian Conrey's (2003) review. Some of these involve considering wider classes of functions such as the *L*-functions associated with cubic curves, echoing Andre Weil's (1948) proving of RH for zeta-functions of (quadratic) function fields. Here, partly responding to Brian Conrey's claim of a conspiracy among abstract *L*-functions, we will restrict ourselves to the generalized RH in the standard complex function setting, to elucidate dynamic principles using Dirichlet series inside and outside the *L*-function framework.

**The Impossible Coincidence**
To ensure convergence, zeta is expressed in terms of Dirichlet's eta function on the critical strip:

$$\zeta(z) = \left(1 - 2^{1-z}\right)^{-1} \sum_{n=1}^{\infty} (-1)^{n+1} n^{-z} = \left(1 - 2^{1-z}\right)^{-1} \eta(z)$$ and then in terms of the functional equation

$$\zeta(z) = 2^z \pi^{-1+z} \sin\left(\frac{\pi z}{2}\right) \Gamma(1-z) \zeta(1-z) \text{ where } \Gamma(z) = \int_0^{\infty} t^{z-1} e^{-t} dt \text{ in the half-plane real}(z) \leq 0.$$

In terms of investigating the convergence of the series to its zeros, eta is better placed than zeta because the convergence is more uniform, as shown in fig 2.

RH is so appealing, as an object of possible proof, because of the obvious symmetry in all the zeroes lying on the same straight line, reinforced by Riemann's reflectivity relation:

$$\Gamma\left(\frac{z}{2}\right)\pi^{-\frac{z}{2}}\zeta(z) = \Gamma\left(1-\frac{z}{2}\right)\pi^{-\frac{1-z}{2}}\zeta(1-z)$$

making it possible to express the zeros in terms of the function xi: $\xi(z) = \frac{1}{2}z(z-1)\Gamma\left(\frac{z}{2}\right)\pi^{-\frac{z}{2}}\zeta(z)$ for which $\xi(z) = \xi(1-z)$, so it is symmetric about $x = \frac{1}{2}$, leading to any off-critical zeros of zeta being in symmetrical pairs. The function $\Xi(z) = \Gamma\left(\frac{z}{2}\right)\pi^{-\frac{z}{2}}\zeta(z)$ also has this symmetry. It also applies to the L-functions, for which: $\xi(z,\chi) = \Gamma\left(\frac{z+a}{2}\right)(\pi/q)^{-\frac{z+a}{2}}L(z,\chi)$, $a = \{\chi(-1) \equiv -1\}$, with $i^{a}k^{1/2}\xi(z,\chi) = \sum_{n=1}^{k}\chi(n)e^{2\pi ni/k}\xi(1-z,\overline{\chi})$.

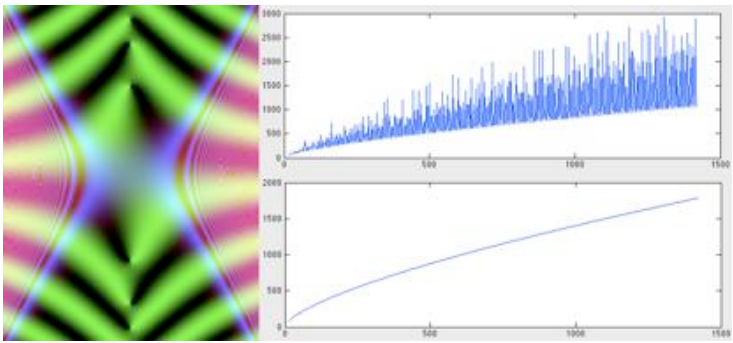

Fig 2 (Left) Symmetry in xi means any off-critical zeros have to be in symmetric pairs. (Top) The number of iteration steps in the eta-derived zeta series required to get 5 steps with 0.005 of 0 varies erratically from one zero to the next, but this is a disguised effect of the presence of the $1/(1-2^{1-z})$ term so becomes a smooth curve for eta (below).

However, when we come to examine the convergence in detail, this symmetry seems to be lost in the actual convergence process. Each term in the series for zeta is $n^{-x+iy} = n^{-x}(\cos(y\ln n) + i\sin(y\ln n))$, forming a series of superimposed logarithmic waves of wavelength $\lambda = \frac{2\pi}{\ln n}$, with the amplitude varying with $n^{-1/2}$ for points on the critical line. Unlike power series, which generally have coefficients tending to zero, Dirichlet L-functions have coefficients all of absolute value 1, which means all the wave functions are contributing in equal amplitude in the sum except for the fact that the real part forms an index determining the absolute convergence. So RH is equivalent to all the zeros being at the same real (absolute) address.

Fig 3: Top left sequence of iterates of eta for the 20,000[th] zero, showing winding into and out of a succession of spirals linking the real and imaginary parts of the iterates. Top right: The wave functions are logarithmic, leading to powers, but not multiples, having harmonic relationships. Below is shown the real and imaginary parts of the iterates (blue and red) overlaid on the phase angle of individual terms (yellow). The zero is arrived at only after a long series of windings interrupted by short phases of mode-locking in the phases of successive terms.

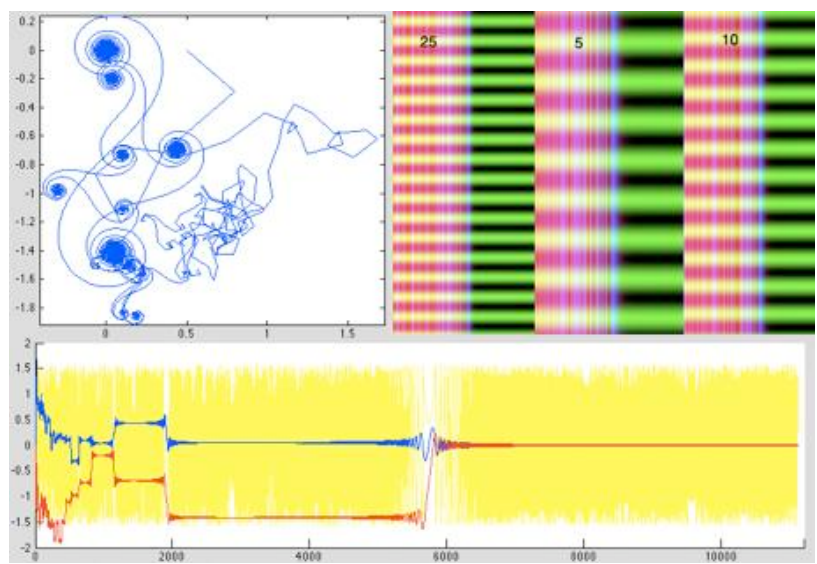

The logarithmic variation means that the wave functions are harmonic only in powers, e.g. 5, 25, 125 etc. and not in multiples. There is no manifest relationship between $\ln n$ and $n^{1/2}$ that explains why the zeros should be on $x = \frac{1}{2}$ and indeed we will find examples where they are not, so there is another factor involved - the primes. Powers of primes or their negation are reflected in both Riemann's primality proofs and other functions, such as the [Möbius function](#):

$$\frac{1}{\zeta(z)} = \sum_{n=1}^{\infty} \frac{\mu(n)}{n^z}, \quad \mu(n) = \begin{cases} (-1)^k, & n \text{ has } k \text{ distinct prime factors of multiplicity 1} \\ 0 & \text{otherwise} \end{cases}.$$

The iterative dynamics give an immediate clue to the potential uncomputability of this problem. If we take a given zero of eta, say the 20,000th, and plot the iterates, we find successive *n*-term approximations wind into and out of a series of spirals associated with non-phase locked epochs, where the angle of successive terms is rotating steadily, interrupted by briefer periods of phase locking, where the angles remain transiently static and hence the complex values of the iteration make a systematic translation. Eventual convergence to zero or another final value occurs only after the last of these mode-locking episodes (see appendix 1), whose iteration numbers can be calculated directly, by finding where the waves match phase:

$$y\ln(n+1) = y\ln(n) + k\pi \Leftrightarrow \ln\left(\frac{n+1}{n}\right) = \frac{k\pi}{y} \Leftrightarrow 1 + \frac{1}{n} = e^{k\pi/y} \Leftrightarrow n = \frac{1}{e^{k\pi/y} - 1}.$$

This corresponds also to the mode shifts in the phase-locking of the orbits in yellow in fig 3. Between the phase locked translations, the iterative value winds towards and then away from an equilibrium value because the angular rotation tends to periodically cancel the effects of intervening terms. After the last phase-translation, further terms simply cause asymptotic convergence to the equilibrium. These effects are all caused because we are dealing with a discrete sum $\sum_{n=1}^{\infty} a(n)n^{-z}$, rather than the continuous integral, which in the case of zeta would simply be polynomial integral $\int_0^{\infty} t^{-z} dt$. It is the transient discrete effects of the phase-locked translations, which determine the eventual value of any Dirichlet series at a given point, so effectively we have a discrete computational problem for each potential zero over the integers, at least up to the last phase translation. This suggests that although the zeros of zeta and the *L*-functions lie hovering temptingly on the critical line, their location can be determined within ε only by explicit computation over the sequence of terms, suggesting RH is a potentially unprovable problem of non-inductive integer computation just as simpler unproven conjectures such as the Collatz conjecture are, although palpably true in each finite case (King 2009).

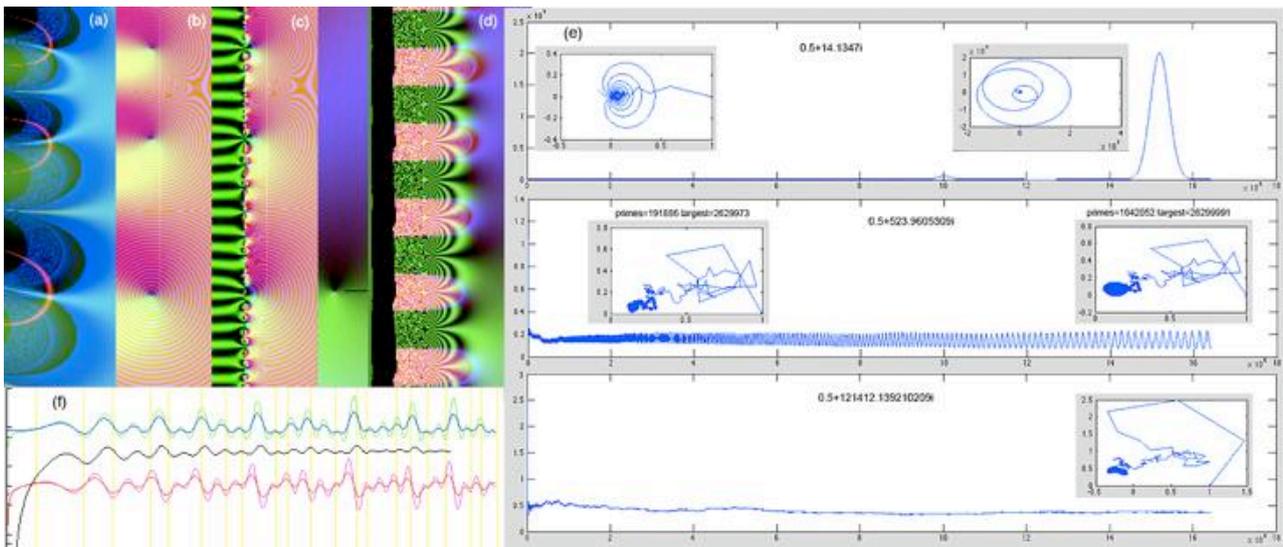

Fig 4 (a) X-ray view of zeta with curves of Re(ζ)=0 (red) and Im(ζ)=0 (cyan) show neither alone determines the location of the zeros. (b,c) log(abs(1/(1-ζ))) plot of the analytic and product forms of zeta show their divergence for *x*<1 and identity for *x*>1. (d) Large fluctuation at the first zeta zero for the product of 84270 primes due to the (red) tongue moving across the zero as the number of product terms increases. (e) Iterative dynamics of the product are radically unstable, leading eventually to exponentiating fluctuations even at the zeros, but these take an extreme number of primes to appear for higher zeros. (f) Fluctuations of real (blue,green) and imaginary (red,magenta) parts of zeta along *x*=1 approximate those of *x*=1/2, the zeros (yellow), and the Fourier sin transform (black) of an integer step function.

It is difficult to apply the Euler product directly to the zeros because it is radically non-convergent in the critical strip and equality with the Dirichlet series holds only for $x>1$ and although variations in values along the line $x=1$ where the sum and product formulations are equivalent do approximate the real and imaginary fluctuations along the critical line.

In fig 4 are shown some of the dynamic features of the Euler product of zeta in comparison with the analytic Dirichlet sum. The sum and product representations diverge in the half plane $x<1$ while being identical on $x>1$. In the critical strip, the iterated product has radical divergence with orbits at the zeros first erratically fractal before setting into exponentiating pulses of divergence, as tongues of large value move down the strip with escalating prime values. When we evaluate the cumulative product up to the 1,642,052$^{th}$ prime 26299991, we find the first zero $y\sim14$ (top) has grown to a peak of around 10 million, while the zero $y\sim523$ (middle) has only begun to enter its first oscillatory burst around the 200,000$^{th}$ prime of around 3 million and $y\sim121412$ is as yet showing no signs of having fully explored its fractal dynamics

However zeta values along $x=1$ do fluctuate in a way which approximates both the imaginary values of the zeros and a Fourier sin transform of an integer step function (the corresponding prime transform also reflects the zeta zeros - see Conrey), showing the distribution of the zeros is transform-based, as demonstrated in Riemann's original proof.

Generally the existence of an Euler product formulation for the sum is seen as a pre-condition for well-behaved *L*-functions and a way of generating new types of *L*-function through prime mediated generators such as elliptic curves which form Euler products determining sum coefficients through prime factorization, which also possess a functional equation representation in the left-half plane.

**Primes and Mediants - Equivalents of RH**
Riemann developed an explicit formula for the prime counting function $\pi(x)$ which is most easily expressed in terms of the related prime counting step function $\psi(x) = \sum_{n \leq x} \Lambda(x)$, the additive <u>von Mangoldt</u> function, where $\Lambda(x) = \log p$ if $x = p^k$ and 0 otherwise. Notice here the exclusive appearance of prime powers eliminated in the Möbius function. We then have the explicit formula
$$\psi(x) = x - \sum_{\substack{\zeta(\rho)=0 \\ 0<\text{Re}(\rho)<1}} \frac{x^\rho}{\rho} - \frac{1}{2}\log(1-x^{-2}) - \log(2\pi),$$ where $\rho = 1/2 + it$ are the zeros of $\zeta(z)$, and the summation is over zeros of increasing $|t|$.

From Ingham (1932 83), we have $\pi(x) = \text{li}(x) + O(x^\theta \ln x)$ where $\theta = \sup_{\rho:\zeta(\rho)=0}(real(\rho))$, $\text{li}\,x = \int_0^x \frac{dt}{\ln t}$.
Hence the asymptotic behavior of the primes is determined by the real sup of the zeros. This comes about because the explicit formula shows the magnitude of the oscillations of primes around their expected position is controlled by the real parts of the zeros of the zeta function, since
$$\frac{x^\rho}{\rho} = \frac{e^{(p+iq)\ln(x)}}{p+iq} = \frac{x^p(\cos(q\ln x) + i\sin(q\ln x))}{p+iq} = \frac{x^p}{p^2+q^2}(\cos(q\ln x) + i\sin(q\ln x))(p - iq)$$
so $\frac{x^\rho}{\rho} + \frac{x^{\bar\rho}}{\bar\rho} = 2\,\text{real}\left(\frac{x^\rho}{\rho}\right) = 2\frac{x^p}{p^2+q^2}(p\cos(q\ln x)) + q\sin(q\ln x)) \sim 2\frac{x^p}{p^2+q^2}q\sin(\ln(x)q)$.

Hence RH has been shown to be equivalent to the statement $|\pi(x) - li(x)| < x^{1/2}\log(x)/8\pi$.

A further equivalent of RH is that $M(x) = \sum_{n \leq x} \mu(n) = O(x^{1/2+\varepsilon})$, which would guarantee the Möbius function would converge for $x > \frac{1}{2}$, and show there were no infinite poles (and hence no zeta zeros). Likewise we have $\sum_{n \leq x} \lambda(n) = O(x^{1/2+\varepsilon})$, $\lambda(n) = (-1)^{\Omega(n)}$, $\Omega(n)$ = no prime factors with multiplicity , the [Liouville function]. Even more basic functional approximations have been found using the floor function (Cloitre). However Mertens conjecture that $M(n) = \sum_{k=1}^{n} \mu(k) < n^{1/2}$, which would have proved the Riemann hypothesis, was found false at a value of around $10^{30}$ by Odlyzko and Herman te Riele (1985), who also showed that $\pi(x) < \text{li}(x)$ fails for some unspecified $x < 6.69 \times 10^{370}$. Even more unachievable potential anomalies arise from considering the number of zeta zeros up to $T$: $N(T) = \frac{T}{2\pi} \log\left(\frac{T}{2\pi}\right) - \frac{T}{2\pi} + \frac{7}{8} + S(T) + O\left(\frac{1}{T}\right)$, $S(T) = \frac{1}{\pi} \arg \zeta\left(\frac{1}{2} + iT\right) = O(\log T)$.

If RH is true we have a much closer bound $S(T) = O\left(\frac{\log(T)}{\log(\log(T))}\right)$ (Ivic). Odlyzko (1992) showed that $S(T)/(\log(T))^{1/2}$ resembles a Gaussian random variable with mean 0 and variance $2\pi^2$, which means it is occasionally much larger than $(\log(\log(T)))^{1/2}$. These results suggest we may only see asymptotic behavior when $|S(T)|$ reaches beyond current limits of around 3.2 (Odlyzko 2002) to values such as 100, implying $T \sim 10^{10^{100}}$, beyond reach of current computational methods.

The Farey sequences appear in a third manifestation of RH (Franel and Landau 1924). These consist of all fractions with denominators up to $n$ ranked in order of magnitude - for example, $F_5 = \left\{\frac{0}{1}, \frac{1}{5}, \frac{1}{4}, \frac{1}{3}, \frac{2}{5}, \frac{1}{2}, \frac{3}{5}, \frac{2}{3}, \frac{3}{4}, \frac{4}{5}, \frac{1}{1}\right\}$. Each fraction is the mediant (see appendix 1) of its neighbours (i.e. $\frac{n_1 + n_2}{d_1 + d_2}$). For an adjacent pair $\frac{a}{b}, \frac{c}{d}$, $bc - ad = 1$. Because the sequence of fractions removes degenerate common factors from the numerator and denominator, they are relatively prime and hence $|F_n| = |F_{n-1}| + \varphi(n)$ since $F_n$ contains $F_{n-1}$ plus all fractions $\frac{p}{n}$ where $p$ is coprime to $n$.

Two [Farey sequence] equivalents of RH state:

(i) $\sum_{k=1}^{m_n} |d_{k,n}| = O(n^r)$, any $r > 1/2$ and (ii) $\sum_{k=1}^{m_n} d_{k,n}^2 = O(n^r)$, any $r > -1$

$d_{k,n} = a_{k,n} - \frac{k}{m_n}$, where $m_n$ is the length of the Farey sequence $\{a_{k,n}, k = 1, \cdots, m_n\}$

This is saying that the **Farey fractions are as evenly distributed as they can be** (to order $n^{1/2}$) **given that they are by definition not evenly distributed [1]**, but determined by fractions with all (prime) common factors removed.

The same consideration applies to the asymptotic **distribution of the primes - they are as evenly distributed as they can be** (to order $n^{1/2}$ from li($n$)) - **given that they are not evenly distributed [2]**, being those integers with no other factors.

This is reflected in other properties of the prime distribution, despite its manifest irregularity, in such processes as the quadratic Ulam spiral. For example, the Dirichlet prime number theorem, states that there are infinitely many primes which are congruent to *a* modulo *d* in the arithmetic progression *a*+*nd*. Stronger forms of Dirichlet's theorem state that different arithmetic progressions with the same modulus have approximately the same proportions of primes. Equivalently, the primes are evenly distributed (asymptotically) among each congruence class modulo *d*.

What RH - that **the non-trivial zeros of the zeta function are all on the critical line [3]** - shows us is the order to which these fluctuations approach an even distribution is inverse quadratic because all the zeros appear to lie on $x = ½$. However the lack of a proof of RH suggests that these three statements are encoded forms of one another and that the locations of the zeros are a **consequence** of the distribution of primes rather than proving their asymptotic distribution, or at best that the three statements are encoded versions of one another. Thus RH is either true but unprovable except in finite numerical approximations, or a type of additional axiom like the axiom of choice that arises from infinities in calculation, just as the Collatz, and other discrete infinity problems appear to be versions of the undecidable Turing halting problem. Turing himself tried to prove computationally that RH was false! (Booker 2006).

We now turn to examining how a dynamical interpretation of the zeta zeros can explain why zeta and the Dirichlet *L*-functions have their non-trivial zeros on the critical line as a result of the asymptotically even distribution of the primes avoiding mode-locking which could knock the zeros 'off-line', as is the case for related functions where mode-locking is more pronounced.

**A Mode-Locking View of Dirichlet *L*-functions and their Counterexamples**

When we look at the sum formula for zeta, it appears to be simply a sum of powers of integers without the primes we see in the product formula, however, immediately we turn to zeta variants such as $\dfrac{1}{\zeta(s)} = \sum_{n=1}^{\infty} \dfrac{\mu(n)}{n^s}$ and $\dfrac{\zeta'(s)}{\zeta(s)} = -\sum_{n=1}^{\infty} \dfrac{\Lambda(n)}{n^s}$, we see the primes reappearing the coefficients.

In the context of the natural numbers, the minimally mode-locked numbers are the primes, since the only common factor of a prime with any other number, apart from itself, is 1. If we turn to the *L*-functions, we see their characters are constructed to eliminate any form of mode locking in three distinct ways, while keeping all the non-zero contributions to the superimposed wave function of equal unit weight:
(1) All coefficients of the bases not relatively prime to the period *k* are set to zero,
    leaving $m = \phi(k)$ relatively prime coefficients.
(2) The remaining coefficients of the relatively-prime bases are distributed cyclically
    with equally weighted values of absolute value 1 in the *m*-th roots of unity, according to
    a power of a generator of the *m* units of *Z*/*Zk*.
(3) Since the group generators result in a sum that can also be represented as a product function
    over primes, the asymptotic distribution of primes places a final limit on any phase-locking.

The negation of the non-relatively prime bases is consistent with the removal of one or more series $\sum_{n=1}^{\infty} (\pm 1)^{n-1} (qn)^{-z}$, $q \setminus k$, for which RH applies, but the distribution around the relatively prime residues with rotating coefficients arises from the group generators and the product representation, which again shows the primes becoming evident in the sum formula. Thus although periodic solutions might appear to be mode-locked these periodic solutions are the least mode-locked coefficient series in terms of integrating in an equi-distributed way with the prime distribution.

These conditions have been abstractly generalized into the four axioms of the Selberg class, attempting to define the conditions causing a Dirichlet series *L*(*z*) to have zeros on the critical line:
(1) Functional equation and  (2) Euler product
(3) Coefficients of order 1. Ramanujan conjecture $a_1 = 1$, $a_n \ll n^\varepsilon$ $\forall \varepsilon > 0$.

(4) At most a single simple pole infinity at 1 i.e. $(z-1)^m L(z)$ analytic for some $m$.

Pivotally the existence of an Euler product is a signature of non-mode-locking because, in a product structure, each of the factors are acting independently with no feedback between them. We shall firstly look generally at Dirichlet series and then focus firstly on *L*-functions that do have both a product structure and a functional equation and then on other variants that arise from products. In abstract *L*-functions, the form of the functional equation varies discretely, with a finite number of gamma factors dependent on the underlying topology of the prime process generating the product. The Ramanujan conjecture separates functions with weight 1 from different weightings which have non-trivial zeros on a different critical line (see later).

To assess the status of RH, we thus consider a wider class of Dirichlet series functions, to explore the effects of mode-locking of the wave functions in the critical strip. As a starting point we look at series where the coefficients are all 0 or roots of unity, but do not satisfy *L*-function conditions. The only Dirichlet *L*-function solutions from the finite group theory are periodic, the period $kn$ consisting of characters in $k$ that are perfect periodic repeats of $k$ characters and not cyclic, or fractal permutations. Non-primitive characters are likewise generated from homologies of the residue groups $\chi_{kn}(p) = \chi_k(p \mod k)$, $\chi_{kn}(p) \neq 0$. Key here is the requirement for complete multiplicativity arising from the Euler product, each integer being a unique product of primes.

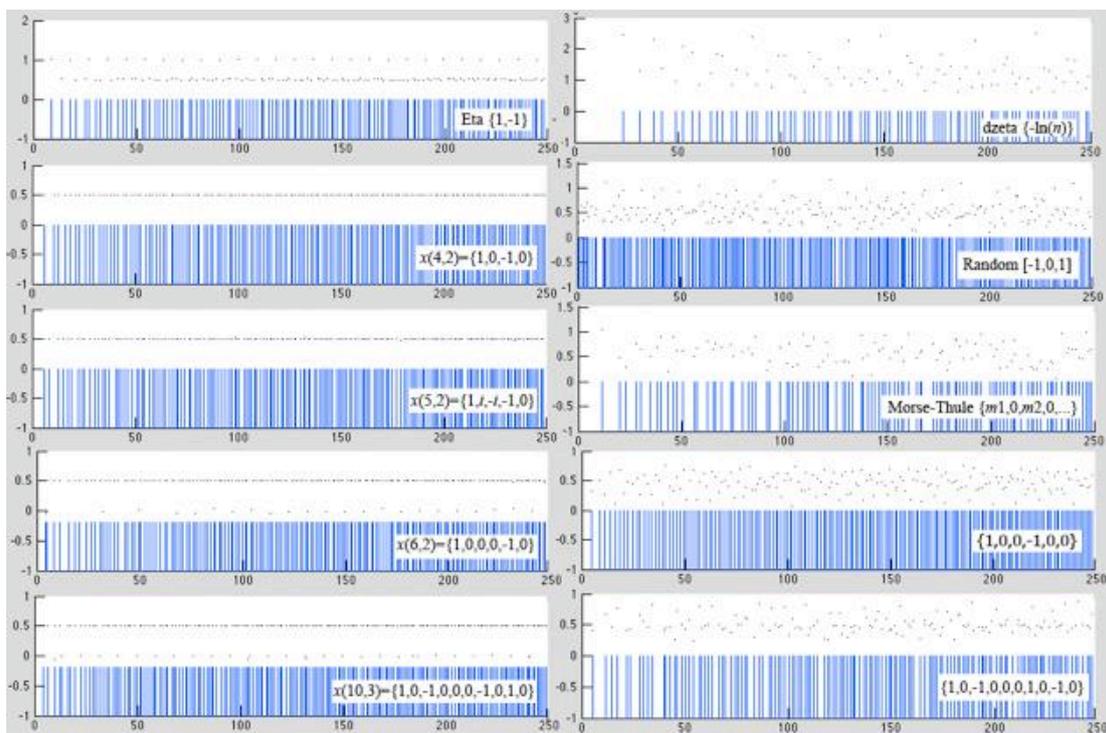

Fig 5: A series of *L*-functions and Eta (left) and RH-violating Dirichlet functions (right) whose critical strip zeros are illustrated by plotting their location (above) and their cumulative frequency, using a Matlab Newton's method scan.

In fig 5 on the left are shown the zeros of eta and a set of typical *L*-functions, confirming both the confinement of the zeros on $0 < x < 1$ to x = ½, and the $t.\ln(t)$ related cumulative frequency, discovered in Riemann's analysis of the zeta zeros. The method looks along a series of closely-spaced values running vertically for local absolute minima and then performs Newton's method using the approximate formal derivative for small $h$. On the right are shown a series of greater and lesser violations of the *L*-function / Selburg class conditions. Note that the derivative of zeta, despite not having an Euler product, does inherit a functional equation from zeta and its zeros are wide of the critical line, implying the functional equation is by no means sufficient, although it does define a symmetry about the critical line. Further examples are the Hurwitz and Davenport-Heilbron zeta functions (see later).

From the top down we have the derivative of zeta $\zeta'(z)$ by formal differentiation of the functional equation, which has terms effectively growing with $-\ln(n)$. Its zeros, corresponding to critical points of zeta, extend far out of the critical strip with an average real part of over 1. The next are Dirichlet series of random equi-distributed integers from $\{-1, 0$ and $1\}$. This shows zeros distributed with means close to $x = ½$. Morse-Thule is a fractal sequence with even coefficients zero and the vector of odd coefficients recursively generated by $v = [v, –v]$ with initial condition $v = 1$ viz $\{1,-1,-1,1,-1,1,1,-1 …\}$ Again this has a mean close to $x = ½$.

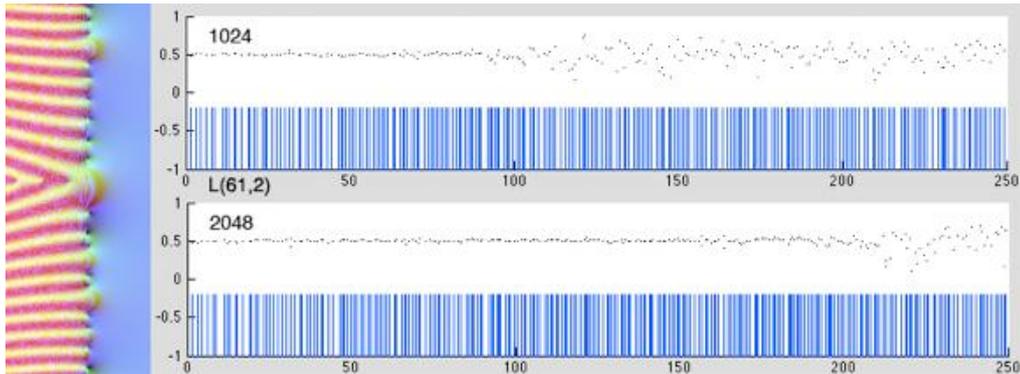

Fig 6: Even with a confirmed *L*-function such as Dirichlet $L(61,2)$ higher periods cause delayed convergence, requiring a disproportionate number of function terms to recognize zeros tending to the critical line.

The last two are variants of the *L*-functions on their left by minor substitution. The first is effectively an alternating arithmetic series in 3's similar to that in 2's of $\chi$ (4,2), namely $1^{-z} - 4^{-z} + 7^{-z} - \cdots$, showing arithmetic series of bases appear to have zeros on the critical line if and only if they correspond to *DL*-functions. In particular, these modified series are not necessarily completely multiplicative, as all *L*-functions are leading to them not having a straightforward expression as an Euler product of primes. This may itself be sufficient reason for the non-*L* functions to be off-critical.

| Function | Means over zeros in [0,1000] |
|---|---|
| Dzeta | 1.1174 |
| Random [-1,0,1] | 0.5306, 0.4891, 0.4905 |
| Morse-Thule ±{0,+/-1,0,-/+1} | 0.5161 |
| Golden Angle Rotation | 0.6290 |
| {1,0,0,-1,0,0} | 0.4761 |
| {1,0,-1,0,0,0,1,0,-1,0} | 0.4959 |

Table 1: Some average *x* coordinates in the critical strip

From table 1 we can also see that, although these variants may have neither a functional equation nor an exact symmetry around the line $x = ½$, the mean real value of their zeros still lie close to the critical line. This is also consistent with the average trends in zeta functions. For example, if we take the curve $f(x) = 1 - \text{geometric mean}(abs(\zeta(x+iy)-1))$, $y=20...120$ step $0.01$, we find it has a zero at ~0.5646, reflecting the innate symmetry of the xi function of fig 2.

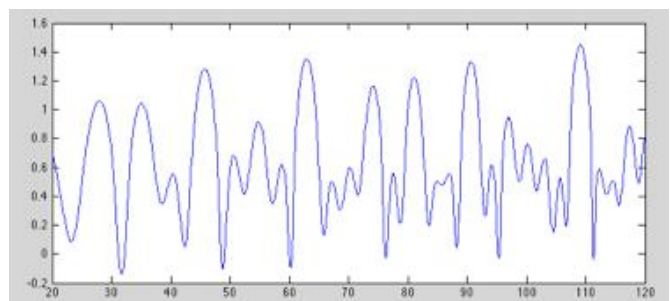

Fig 7: Function $p(y)$ showing the x-coordinate for each y, where the absolute value of zeta differs by 1 from 1.

Alternatively when we take the individual curve $p(y) = \{x : |\zeta(x+iy) - 1| = 1\}$ in the interval [20,120], as in fig 7, we find it has a geometric mean of 0.4965.

While these estimates are just very rough ad-hoc approximations because of the exponentiating irregularity of all these functions, they do indicate how zeros of Dirichlet functions can deviate significantly from the critical line while still having an averaged behavior closely spanning it. There is also no evidence for symmetric pairs of off-critical zeros, as would be required by the symmetry of the functional equations of zeta and the *L*-functions.

There are two additional ways we can compare ideas about the basis of the critical zeros. The first is the notion that the distribution of the zeta zeros reflects the statistics of random matrix theory. The zeros of zeta and their pair correlations have been shown to correspond to a GUE, or grand unitary ensemble. In fig 7b we thus compare these two statistics for the unreal zeros of *DL*(6,2) and the non-*L* function with quasi-character {0,1,0,0,-1,0} illustrated in fig 5 up to 2500*i*. Although it is true that *DL*(6,2), conforms a little more closely to the GUE statistic and there is more evidence for sustained phase-locking in the enhanced periodic fluctuations of the pair correlation, the idea that GUE is a defining indicator for criticality is less than convincing.

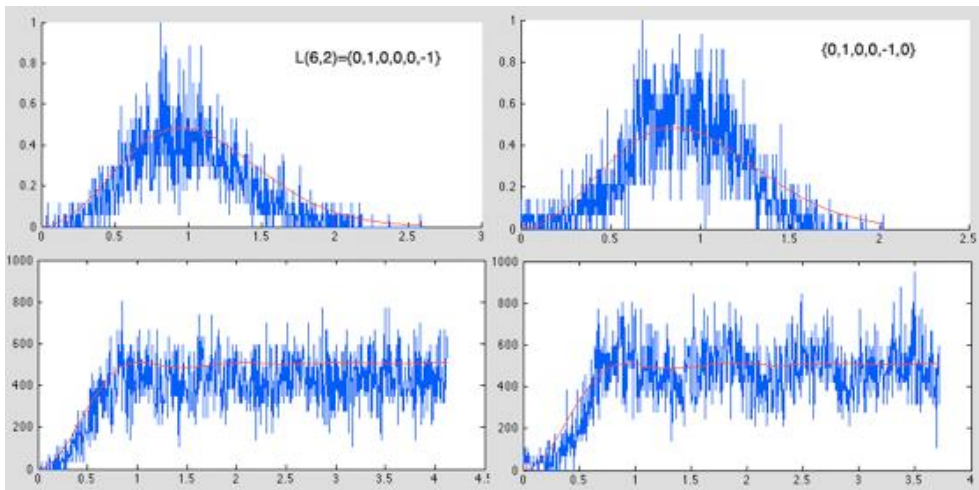

Fig 7b: (Above) distribution of the spacing of the zeros and (below) pair correlations for *DL*(6,2) with periodic zeros removed and the non *L* function with quasi-character {0,1,0,0,-1,0}. GUE distribution in red.

We can also examine the way in which convergent *DL* and non-*L* functions generate 'prime counting' functions using variants of the explicit formula above for zeta. We will use the simplified formula $\varphi(x) = \sum_{\substack{L(\rho)=0 \\ 0 \leq \text{Re}(\rho) < 1}} \frac{x^\rho}{\rho}$, $\rho = x + iy$ counting the zeros in the critical strip in order in both directions from $y = 0$. In fig 7c the results are illustrated. Notably, both (5,2) and (6,2) correctly shift at primes and prime powers relatively prime to the period, but (6,2) does this only when the periodic zeros on *x*=0 are also included. Even more intriguing, the non −*L* function (0,1,0,0,-1,0} also counts shifts unperturbed by its off-critical zeros and correctly deletes shifts for terms having more than one factor in the series – i.e. 28=4x7, 52=4x13, 70-7x10,76=4x19 and 91=7x13.

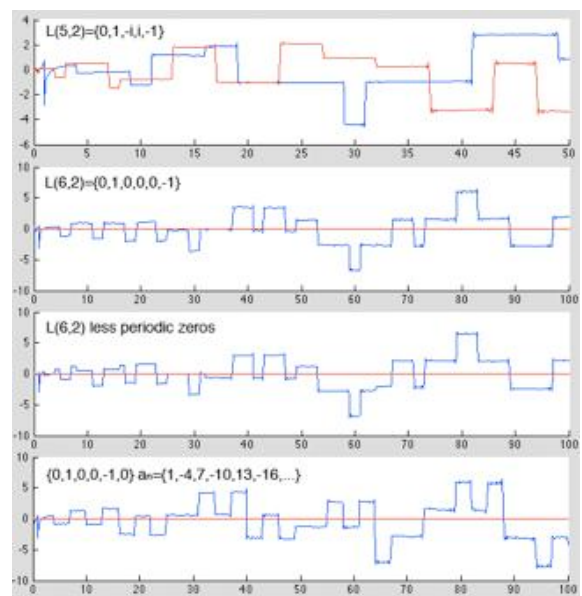

Fig 7c: *DL*(5,2) has a prime counting function with real and imaginary parts shifting precisely at primes $p \neq 5$, and at powers of these primes according to $s(p^n) = (\chi(p))^n$, $\chi = \{0,1,-i,i,-1\}$, reflecting the von Mangoldt definition above. There is no shift at integers with more than one prime factor. *DL*(6,2) has the same profile if the periodic zeros on *x*=0 are included, but if they are removed, spurious shifts occur at powers of 2. The non-*L* function (0,1,0,0,-1,0} forming an arithmetic progression $a_n=\{1,-4,7,-10,13,-16,\ldots\}$ has shifts at each of the $a_n$ except those which have more than one factor from the existing series.

We still lack a broad spectrum of examples lying outside zeta and the Dirichlet *L*-functions where the zeros are on the critical line or its displaced equivalent. Classically all the examples found comprise more general types of zeta and L-functions where the coefficients are determined by more arcane primal relationships, essentially guaranteeing the zeros are on-line through more veiled forms of primal non-phase-locking. In the following section we thus give a portrayal of the key types of abstract *L*-function, with a discussion of how their primal relationships arise.

**Widening the Horizon to other types of Zeta and *L*-Function**

To get a view of how *L*-functions can be extended beyond the context of Riemann and Dirichlet, a first stepping point is given by Dedekind zeta and Hecke *L*-functions of field extensions of the rationals ***Q*** (Garrett 2011). Here we look for the non-zero ideals of the ring of integers in a field extension. These also share features of analytic continuation using functional equations and Euler products. Some such as $Q[\sqrt{-5}]$ do not have unique prime factorizations and require consideration of the so-called class number, in this case 2, as $6 = 2.3 = (1+\sqrt{-5})(1-\sqrt{-5})$.

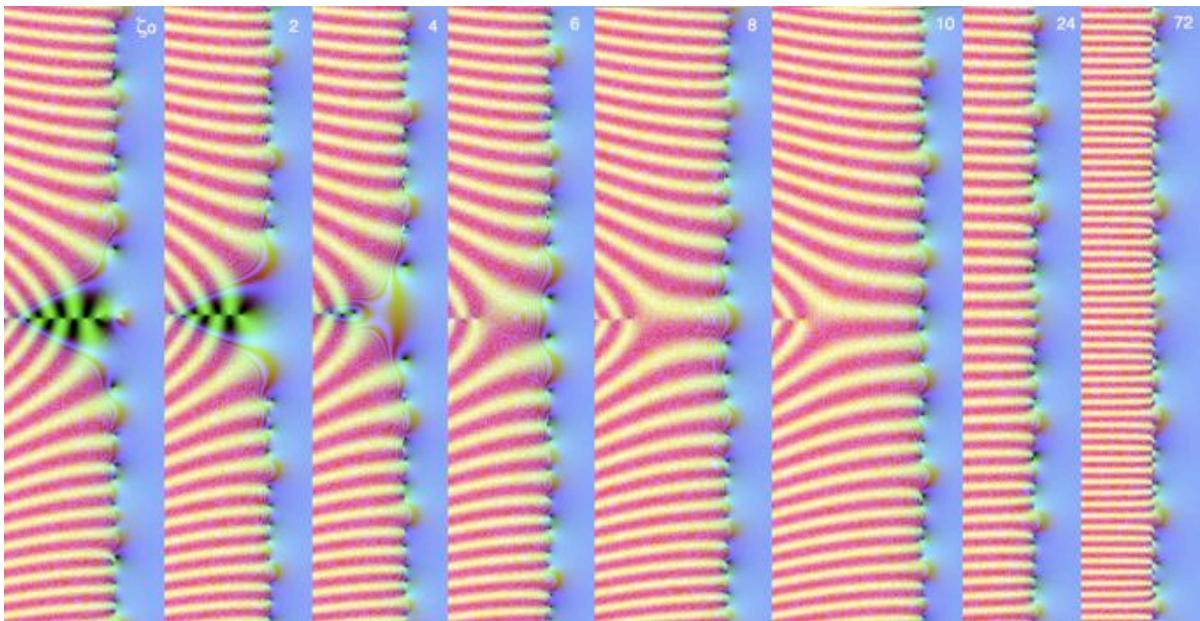

Fig 8: Profiles of the Dedekind zeta and Hecke *L*-functions for ***Z***[*i*], the extension to the Gaussian integers. The portraits require both series representation, and the functional equation and Mellin transform theta integrals.

We will look at those of the Gaussian integers ***Z***[*i*], defined by appending *i* to the integers, resulting in the lattice of complex numbers with integer real and imaginary parts. Here we have $N\alpha = \alpha\bar{\alpha} = |\alpha|^2$, so $\zeta_o = \sum_{0 \neq \alpha \in o \bmod o^\times} \frac{1}{(N\alpha)^z} = \frac{1}{4} \sum_{m,n \text{ not both } 0} \frac{1}{(m^2+n^2)^z} = \prod_{\omega \text{ prime}} \frac{1}{1-(N\omega)^{-z}}$, where $N\alpha$ is the norm of the ideal $Z[i]/\alpha Z[i]$, which is uniquely expressible as an Euler product of prime ideals. This has a functional equation $\pi^{-z}\Gamma(z)\zeta_o(z) = \pi^{-(1-z)}\Gamma(1-z)\zeta_o(1-z)$, although, lacking an eta analogue, convergence isn't assured in the critical strip $0 < x < 1$, so Mellin transforms are commonly used to define the function more accurately there.

Correspondingly we have Hecke *L*-functions defined as follows. Consider the multiplicative group $\chi : Z[i] \to S^1$, $\chi(\alpha) \to (\alpha/\bar{\alpha})^l$, $l \in Z$. To give the same value on every generator this requires *l* to be trivial on units, hence $1 = \chi(i) = \left(\frac{i}{-i}\right)^l = (-1)^l$, so $l \in 2Z$. We then have for each such *l* a Hecke *L*-function: $L(z,\chi) = \sum_{0 \neq \alpha \in o \bmod o^\times} \frac{\chi(\alpha)}{(N\alpha)^z} = \frac{1}{4} \sum_{m,n \text{ not both } 0} \frac{(\alpha/\bar{\alpha})^l}{(m^2+n^2)^z} = \prod_{\omega \text{ prime}} \frac{1}{1-\chi(\omega)(N\omega)^{-z}}$ where the primes are now those of Gaussian integers, units ±1 or ±$i$ times one of 3 types: $1+i$ or a real prime which isn't a sum of squares ($p$ mod $4 = 3$), or has sum of real part squared and imaginary part squared a prime ($p$ mod $4 = 1$). Again we have a functional equation:
$$\pi^{-(z+|l|)}\Gamma(z+|l|)L(z,\chi) = (-1)^l \pi^{-(1-z+|l|)}\Gamma(1-z+|l|)L(1-z,\chi).$$

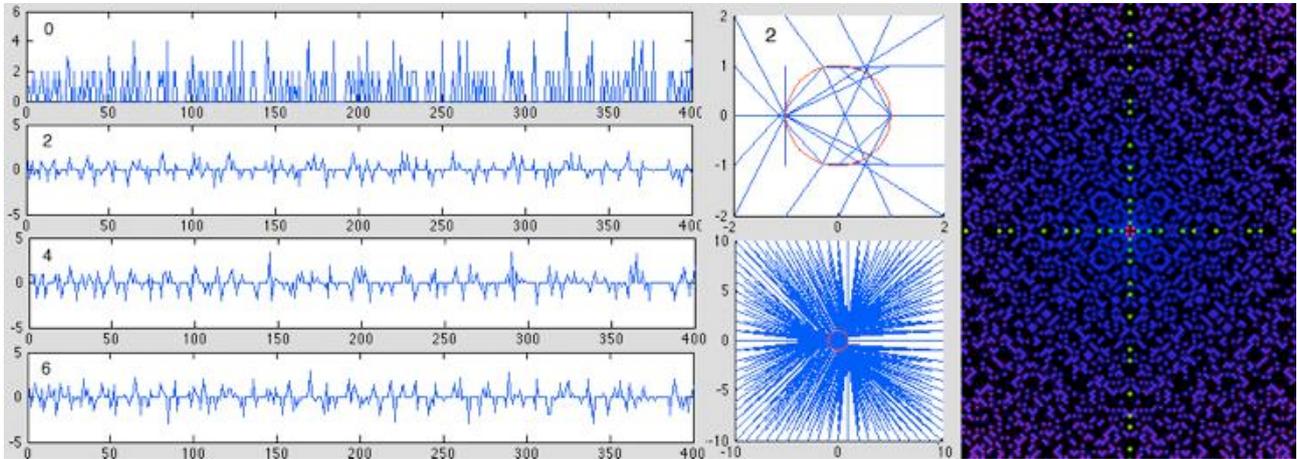

Fig 9: (Left) Profiles of the coefficients. Dedekind zeta (0) consists of the number of ways an integer can be represented as the sum of two integers divided by 4. The Hecke $L$-functions multiply these by the map $\chi: Z[i] \to S^1$, $\chi(\alpha) \to (\alpha/\bar{\alpha})^l$, $l \in 2Z$ to the unit circle illustrated (centre) for the case 2. Effectively this simply multiplies the angle of ($m+in$) by $2l$ and sets the modulus to 1 since $\chi(re^{i\theta}) = (e^{i\theta}/e^{-i\theta})^l = e^{2li\theta}$, $l \in 2Z$. It therefore plays a role similar to the Dirichlet characters in evenly distributing the coefficients. All the coefficients are real and all but zeta fluctuate in sign. (Right) distribution of the Gaussian primes [r (±1±$i$), g (0,±4$n$+3) (±4$n$+3,0), b ($m^2+n^2$)= 4$n$+1: 4$n$+$k$ prime]

The profiles of these functions with their analytic continuations are shown in fig 8, requiring, in addition to the functional equations, use of Mellin transform integral formulae in the critical strip:
$$\zeta_o(z) = \pi^z \Gamma(z) \int_1^\infty (y^z + y^{1-z}) \frac{\theta(iy)-1}{4y} dy, \ \theta(iy) = \sum_{m,n \in Z} e^{-\pi(m^2+n^2)y} = \left(\sum_{n \in Z} e^{-\pi n^2 y}\right)^2$$

$$L(z,\chi) = \pi^{z+|l|}\Gamma(z+|l|)\int_1^\infty (y^z + (-1)^l y^{1-z}) \frac{\theta_\chi(iy)}{4y} dy, \ \theta_\chi(iy) = \sum_{m,n \in Z} (m\pm in)^{2|l|} y^l e^{-\pi(m^2+n^2)y}$$

Counting the coefficients of the Dirichlet sum over the sums of squares, we find:
$$\zeta_0(z) = 1 + 2^{-z} + 0 + 4^{-z} + 2 \cdot 5^{-z} + 0 + 0 + 8^{-z} + \cdots$$

In terms of our original primes in $Z$, we can say they fall into three cases, which will carry over to Hasse-Weil zeta functions: (i) **split** ($p$ mod $4 = 1$) two square roots of -1 in the finite (Galois) field $F_{p^m}$ $m>1$ (see below); (ii) **inert** ($p$ mod $4 = 3$) no square root of -1 in $F_{p^m}$, $m$ odd but 2 if $m$ even; (iii) **ramified** ($p = 2$) one square root of -1. Confirmation for 2, 3, $3^2$, 5 and 7 is in appendix 2.

When we go back to Dedekind zeta's Euler product, we see that the product over Gaussian primes coincides exactly with an Euler product over integer primes incorporating the above cases and both generate the sum coefficients from unique prime power factorisations:

$$\prod_{\omega \text{ prime}} \frac{1}{1-(N\omega)^{-z}} = \left(\frac{1}{1-(1^2+1^2)^z}\right)\left(\frac{1}{1-(3^2)^z}\right)\left(\frac{1}{1-(1^2+2^2)^z}\right)\left(\frac{1}{1-(2^2+1^2)^z}\right)\cdots$$

$$= \frac{1}{1-(2)^z} \prod_{p \bmod 4=1} \frac{1}{(1-(p)^{-z})^2} \prod_{p \bmod 4=3} \frac{1}{1-(p)^{-2z}} = \left(\frac{1}{1-(2)^z}\right)\left(\frac{1}{1-(3)^{2z}}\right)\left(\frac{1}{(1-(5)^z)^2}\right)\cdots$$

$$= 1 + 2^{-z} + 0 + 4^{-z} + 2 \cdot 5^{-z} + 0 + 0 + 8^{-z} + \cdots$$

Alternatively, we can count the series terms directly in terms of a category mapping (functor) from commutative rings to sets, which preserves products and takes finite fields to finite sets (Baez). Effectively we are going to find how many ways make finite sets into semi-simple commutative rings, which are themselves always finite products of finite fields, which in turn have one field of $q$ elements when $q=p^n$, $p$ prime, and none otherwise, bringing in the powers of primes at a root level.

We can then make a general abstract Hasse-Weil zeta function $\zeta_s(z) = \sum_{n \geq 1} \frac{|Z_s(n)|}{n!} n^{-z}$ where $Z_s(n)$ are the species of different ways. To find the number of ways to make rings, we have to factor by the automorphisms of the finite fields that would make equivalent rings. The number of these turn out to be the number of automorphisms in each factor field times the number of permutations of equivalent factors. So we have for $n = 0$, none; $n = 1$, 1 (trivial ring an empty product of finite fields), $n = 2$, 1 ($F_2$); $n = 3$, 1 ($F_3$); $n = 4$, 2 ($F_4$ and $F_2 \times F_2$); $n = 5$, 1; $n = 6$, 1 ($F_2 \times F_3$); $n = 7$, 1, $n = 8$, 3 ($F_2 \times F_2 \times F_2$, $F_2 \times F_4$, $F_8$). Hence for all the cases up to 8 except 4 and 8 we have $n!/1$ ways, but for $n = 4$, we have $4!/2 + 4!/2 = 4!$ ways, the first from $F_4$ and the second from permutations of the $F_2$ factors. We find 8 similarly gives $8!/3 + 8!/2 + 8!/6 = 8!$ ways, so we find for the Riemann zeta function $\zeta(z) = \sum_{n \geq 1} \frac{n!}{n!} n^{-z} = \sum_{n \geq 1} n^{-z}$.

In the case of Dedekind zeta, each coefficient contains a number of ways combining the information from the number of roots of unity in each prime case with the above classification of the natural numbers, i.e. $n = 0$, 0; $n = 1$, 1x1!; $n=2$, 1x2!; $n = 3$, 0x3!; $n = 4$, 1x4!/2+1x4!/2=1x4!; $n = 5$, 2x5!; $n = 6$, (0x1)x6!; $n = 7$, 0x7!; $n = 8$, 1x8!/3 + 1x8!/2 + 1x8!/6 = 1x8! ways, leading again to: $\zeta_0(z) = 1 + 2^{-z} + 0 + 4^{-z} + 2 \cdot 5^{-z} + 0 + 0 + 8^{-z} + \cdots$.

This discussion leads on naturally to the next example of cubic curves where we see essentially the same picture of prime inertness, splitting or ramification, incorporated into an Euler product containing quadratic prime factors.

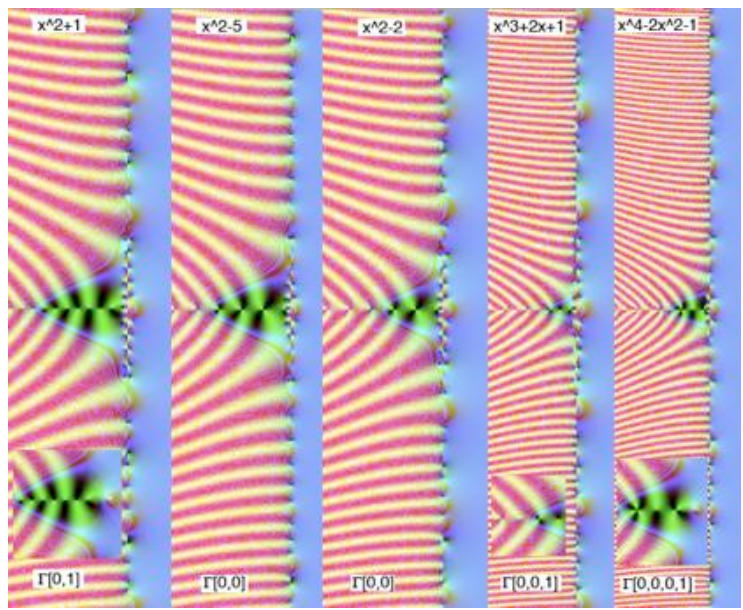

Fig 10: Dedekind zeta functions of a series of extension fields of polynomials portrayed with Dirichlet series and functional equation, but without the use of a Mellin transform in the critical strip, highlighting convergence failure of the Dirichlet series in the critical strip. Note the degenerate zeros in the left half plane caused by repeated gamma factors in the functional equation. Lower-right (inset) Computel Mellin transform portrait of the central valley, correcting the errors in the functional equation portrait.

**L-functions of Elliptic Curves**

The theory of elliptic curves and modular forms also generate *L*-functions (Booker 2008), which involve Euler products with quadratic factors in the denominator. In figs 11, 13 are illustrated a variety of abstract *L*-functions from the genus-1 *L*-function of the elliptic curve $y^2 + y = x^3 - 7x + 6$, through genus-2, 3 and 4 cases with repeated gamma factors causing multiple higher order zeros, to the *L*-function of a modular form based on the Ramanujan's Tau function, and many other cusp forms associated with elliptic curves. Simple scripts to list and generate L-functions of elliptic curves and diverse modular forms via Sage and PARI-GP using Tim Dokchitser's example files to generate the L-function coefficients and gamma factors for loading into RZViewer are included with the RZViewer package. Some simple Sage commands for elliptic curves and modular forms are illustrated in appendix 5.

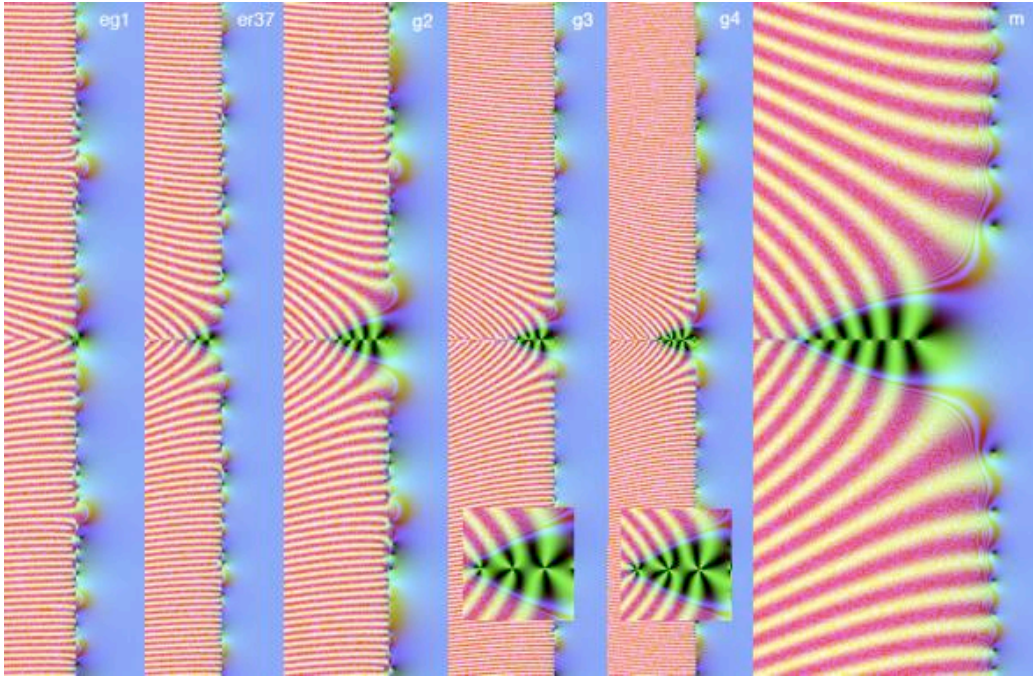

Fig 11: From left to right, $L$-functions of the genus-1 elliptic curve $y^2 + y = x^3 - 7x + 6$, the elliptic curve $y^2 + y = x^3 + 2x^2 + (19 + 8\omega)x + (28 + 11\omega)$, $\omega = (1 + \sqrt{37})/2$ over $K = Q(\sqrt{37})$, the genus-2 curve $y^2 + (x^3 + x + 1)y = x^5 + x^4$, the genus-3 curve $y^2 + (x^3 + x^2 + x + 1)y = x^7 + 2x^6 + 2x^5 + x^4$, the genus-4 curve $y^2 + (x^5 + x + 1)y = x^7 - x^6 + x^4$, and the modular cusp form $\Delta(z) = \sum_{n \geq 1} \tau(n)e^{2\pi i n z}$, of weight 12, the modular discriminant, using Ramanujan's Tau function

$$\tau(n) = (5\sigma(n,3) + 7\sigma(n,5))\frac{n}{12} - 35\sum_{k=1}^{n-1}(6k - 4(n-k))\sigma(k,3)\sigma(n-k,5), \text{ where } \sigma(n,k) = \sum_{d \backslash n} d^k.$$ This is identical to the unique cusp form of weight ≤ 12 over $SL(2,Z)=\Gamma_1(1)$, mg1p1w12 in the notation of fig 15, and so occupies a place among modular cusp forms akin to that of the Riemann zeta function among Dirichlet L-functions.

Fig 12: (Left) Examples of elliptic curves, (right) Group operation.

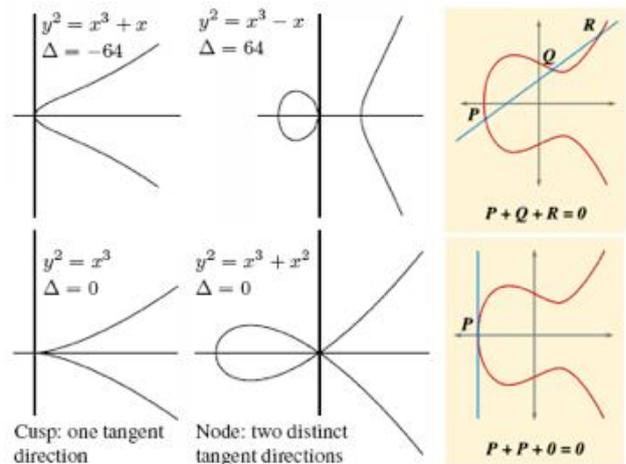

Hasse-Weil $L$-functions of elliptic curves $E$ are generated by taking the function $E(Q)$ over $Q$, or a field extension $F$, and estimating the number of rational points (Silverman 1986). Factoring mod $p$, for primes $p$, to get a set of $A_p$ points on the curve $E(F_p)$ in the finite prime field $F_p$, given up to a maximum of $p+1$ points in $F_p$ (including the point at infinity). We then let $a_p = p+1-A_p$ the number of missing points. For example, for the elliptic curve $y^2 + y = x^3 - 7x + 6$, (0,2), (1,0), (1,4), (2,0), (2,4), (3,1), (3,3), (4,1), (4,3), (∞,∞) are solutions mod 5, giving $a_5 = 5+1-10 = -4$.

Hence we can define:

$$L(E,z) = \sum_{n=1}^{\infty} a_n n^{-z} = \prod L_p(E,z), \quad L_p(E,z) = \begin{cases} \left(1 - a_p p^{-z} + p^{1-2z}\right)^{-1} & \text{good reduction} \\ \left(1 - a_p p^{-z}\right)^{-1} & \text{bad reduction} \end{cases}$$

where bad reduction i.e. a singularity of $E(F_p)$ results from repeated roots in $F_p$, when $a_p = \pm 1$, depending on the splitting or inertness of $p$ (rational or quadratic tangents of the node) for multiplicative reduction ($p|N$ but not $p^2$) of $E$, or is 0 if $p^2|N$ (additive reduction of the cusp), where $N$ is the conductor, the 'effective' product of bad primes. Setting $L^*(E,z) = N^{z/2}(2\pi)^{-z}\Gamma(z)L(E,z)$, we have the functional equation $L^*(E,z) = \varepsilon L^*(E, 2-z)$, where $\varepsilon = \pm 1$. The $a_n$ are generated from the Euler product, convergent for $x>3/2$.

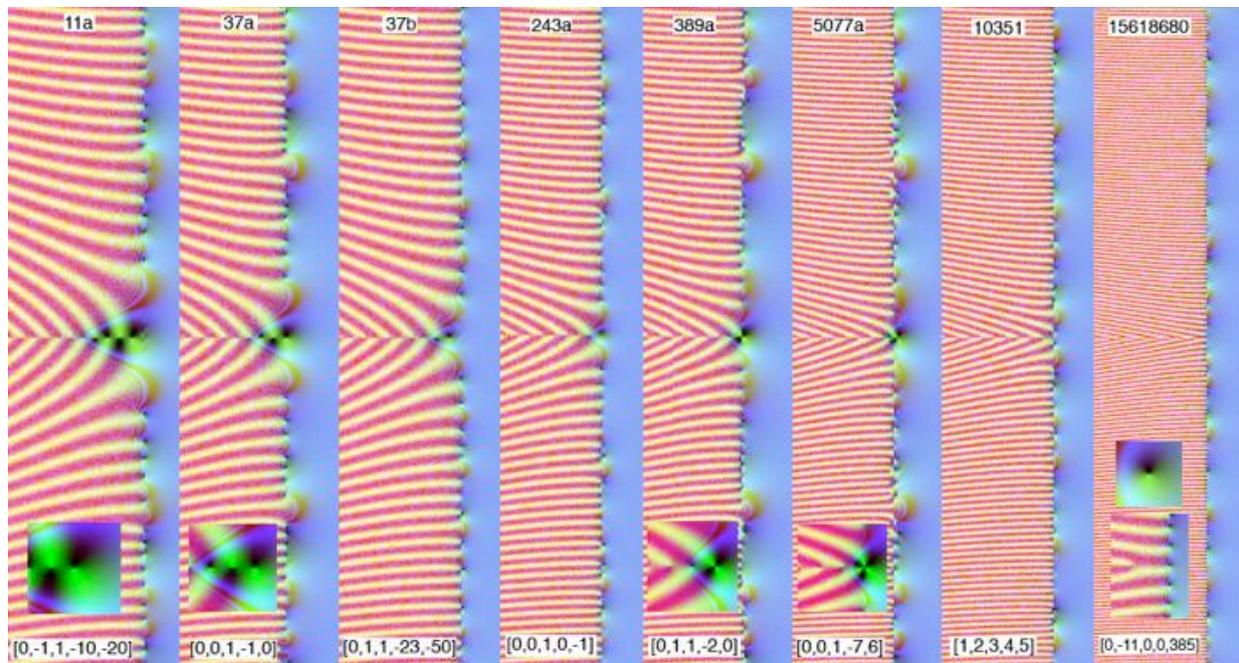

Fig 13: A menagerie of $L$-functions of elliptic curves over $Q$ classified by their conductors (above) and their defining equations (below) where $[a,b,c,d,e]$ corresponds to $y^2 + axy + cy = x^3 + bx^2 + dx + e$. Inset Sage renditions of rank 0 (non-zero at 1), 1, 2 and 3 cases on $[-2,2]^2$ illustrating the Birch and Swinnerton-Dyer conjecture in the multiple-ray angular variation round the point. Unlike the higher-genus cases of fig 11, where repeated gamma factors cause multiple higher order zeros here it applies only to $z = 1$. See appendix 4 for computational method comparisons.

The conductor, as the 'effective' product, differs from the discriminant - a product of all bad prime factors. It consists of factors 1 for good reduction, $p$ for multiplicative reduction, and $p^2$ for additive, except in the cases 2, 3 where the exponent may have an additional 'wild' component, increasing it up to 5 for 2 and 3 for 3, depending on the number of irreducible components (without multiplicity) of the 'special Neron fibre' (Tate, Silverman 1994).

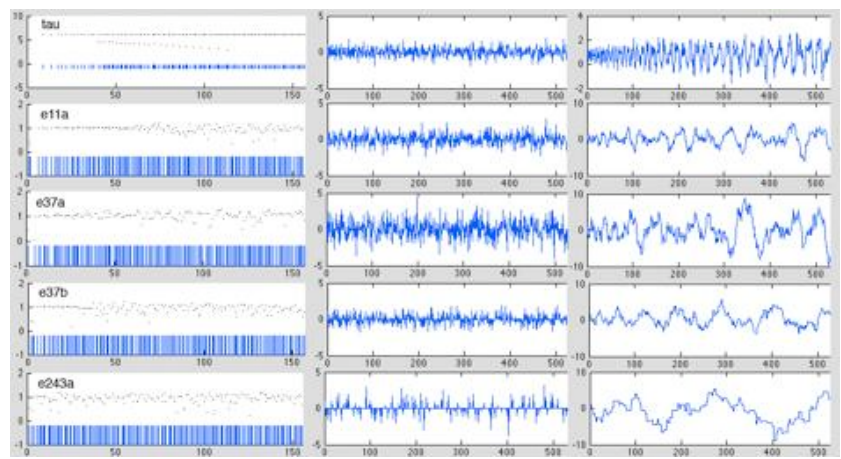

Fig 14: (Left): Newton's method on the Dirichlet series representation at 1000 terms for the modular form Delta and four elliptic curves of lowest conductor, in fig 13, show convergence to the critical line, for smaller imaginary values, similar to that of fig 6 for the Dirichlet $L$-function $L(61,2)$, with convergence diminishing, as the conductor becomes larger, due to longer fluctuations in the sign of the coefficients, illustrated in the rescaled coefficients (centre) and additive trends (right). Tim Dokchitser's Computel, now in Sage, can give a more accurate numerical calculation for individual zeros using Mellin transforms, however these work only for limited imaginary values ($<\pm30$ for 37a) so there is no obvious way to accurately test the generalized RH for these $L$-functions. Moreover the Mellin transform method depends on established functional equations and we are interested in Dirichlet series because they are possessed by both $L$- and non-$L$-functions, which may not have a functional equation.

A good example is the elliptic curve $y^2 = x^3 - 11x^2 + 385$ (Lozano-Robledo), with additive reduction on 2, 11, split multiplicative on 5 and inert multiplicative on 7 and 461:

$$L(z) = \left(1 - 5^{-z}\right)^{-1}\left(1 + 7^{-z}\right)^{-1}\left(1 + 461^{-z}\right)^{-1} \prod_{p \neq 2,5,7,11,461}\left(1 - a_p p^{-z} + p^{1-2z}\right)^{-1} = 1 - \frac{2}{3^z} + \frac{1}{5^z} - \frac{1}{7^z} + \cdots,$$

with conductor $N = 2^3 \cdot 11^2 \cdot 5 \cdot 7 \cdot 461 = 15618680$ and root number -1 (see fig 13).

Elliptic curves have a group multiplication connecting any two points on the curve to the third point of intersection of the line through them, as illustrated in fig 12. The Birch and Swinnerton-Dyer conjecture asserts that the rank of the abelian group $E(F)$ of points of $E$ is the order of the zero of $L(E, z)$ at $z = 1$. Even rank gives $\varepsilon=1$ and odd $\varepsilon=-1$. The group may also have finite torsion elements.

Although the function depends on a rather arcane definition, through an elliptic curve, and then a quadratic Euler product, the resulting Dirichlet series is a standard sequence of coefficients, which possesses a standard functional equation and can thus be portrayed as a meromorphic function in $C$ (analytic except for a finite number of simple infinities). For the elliptic curve $y^2 + y = x^3 - 7x + 6$, the first coefficients are: {1,-2,-3,2,-4,6,-4,0,6,8,-6,-6,-4,8,12,-4,-4,-12,-7,-8, …}.

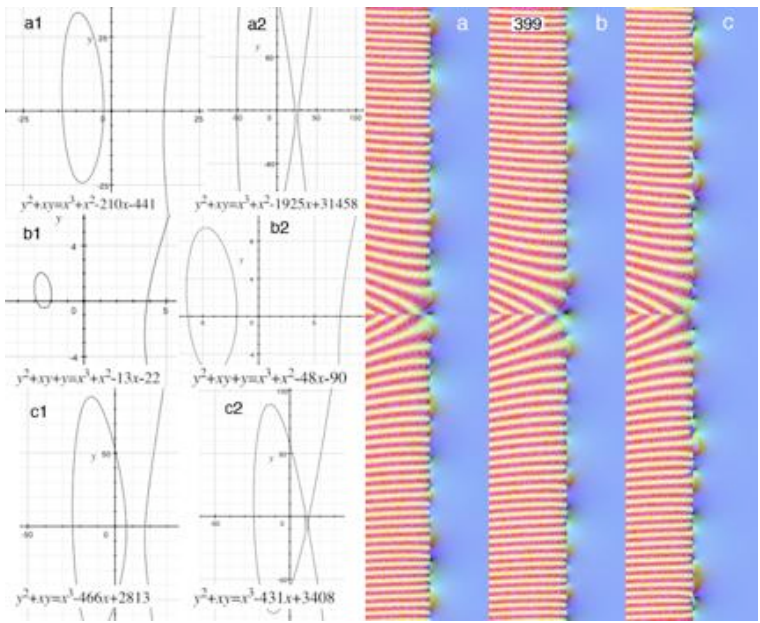

Fig 14b: *L*-functions of the elliptic curve conductor 399=3x7x19 come in three forms each of which has two elliptic curves associated with it. The space of modular cusp forms of weight 2 on gamma0 with conductor 399 has a dimension 53 and all three elliptic curve *L*-functions are linear combinations of several modular form basis functions (see fig 15).

If one takes the defining equation of an elliptic curve, one can generate an algebraic function, which is single-valued on a surface, enabling the elliptic curve to also be represented as a mapping of this surface. This parametrization, via the Weierstrass function and its derivative, defines a "fundamental parallelogram" in the complex plane, representing the two periodicities in the torus. The doubly periodic nature of the function and a one and three-holed torus (see modular forms) are illustrated below left, with the two periodicities illustrated on the one torus.

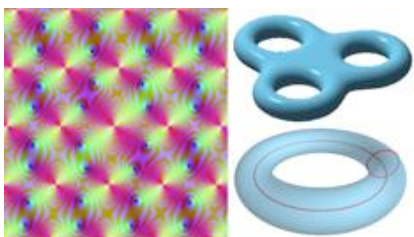

$$\wp(z;\omega_1,\omega_2) = \frac{1}{z^2} + \sum_{(m,n)\neq(0,0)}\left(\frac{1}{(z + m\omega_1 + n\omega_2)^2} - \frac{1}{(m\omega_1 + n\omega_2)^2}\right)$$

$$[\wp'(z)]^2 = 4[\wp(z)]^3 - g_2\wp(z) - g_3,$$

$$g_2 = 60 \sum_{(m,n)\neq(0,0)} (m\omega_1 + n\omega_2)^{-4}, \; g_3 = 140 \sum_{(m,n)\neq(0,0)} (m\omega_1 + n\omega_2)^{-6}$$

Elliptic functions over $C$ are thus genus-1 curves, topologically equivalent to embeddings of a torus in $PC \times PC$ where $PC$ is the complex projective plane or Riemann sphere derived by adding a single point at $\infty$ to $C$. Higher degree curves generate higher genus examples, as illustrated in fig 11.

**Modular and Automorphic Forms**
Complementing the *L*-functions of elliptic curves are those of modular forms. The toroidal nature of the elliptic function, causes it to be periodic on a parallelogram in *C*, resulting in a deep relationship with another kind of *L*-function. A modular function is a meromorphic function (analytic with poles) in the upper half-plane *H*, which is conserved by the modular group $SL(2,Z)$ of integer 2x2 matrices of determinant 1 i.e. $f(az+b)/(cz+d)=f(z)$. More generally we have modular of weight *w* (necessarily even) if $f(az+b)/(cz+d)= (cz+d)^w f(z)$. If it is holomorphic (fully analytic) in the upper half-plane (and at ∞) we say it is a *modular form*. If it is zero at ∞ we say it is a cusp form.

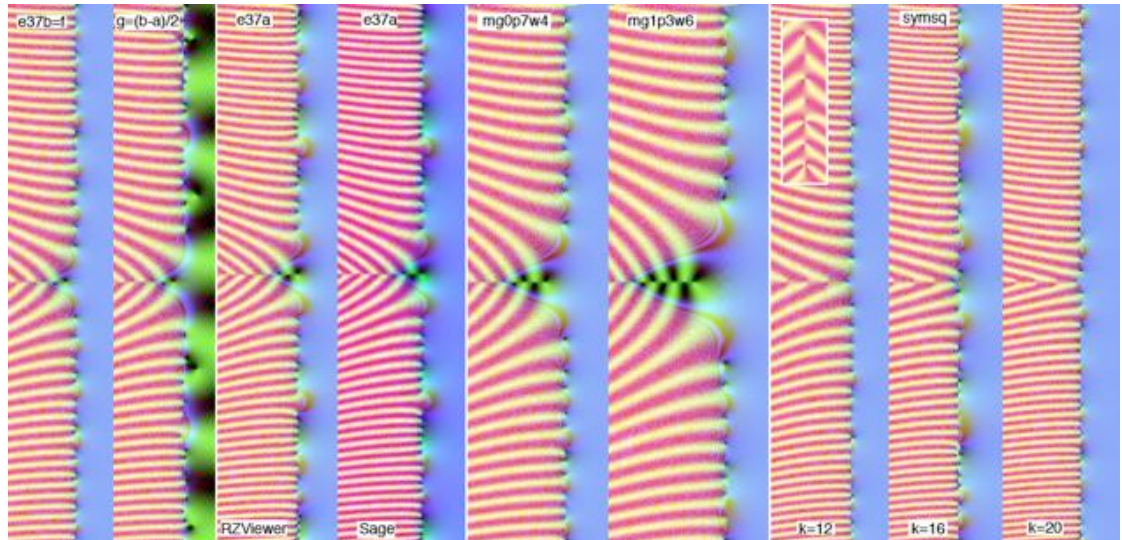

Fig 15: (Left): Functions *f*, *g* representing the 2 dimensions of modular cusp forms in $S_2(\Gamma_0(37))$. Here 37b has rank 0 with $\varepsilon = 1$, but 37a rank 1 with $\varepsilon = -1$. Consequently *g* has a complicated functional equation represented by $(b-a)/2$ for $x \geq 0$ and by $(b+a)/2$ for $x<0$. It also appears to have zeros manifestly off-critical. (Centre-left): Correspondence (see appendix 4) between portrait of e37a using RZViewer and an equivalent portrait using Computel Mellin transform algorithm via Sage. The left hand image takes 8 seconds and the right 4 hours on a Mac intel dual core at 2.1 GHz. (Centre-right) Two modular cusp forms mg0p7w4 and mg1p3w8 of gamma 0, 1, modulo 7, 3 and of weight 4, 6 respectively. In the latter case, and that of mg0p1w12 in fig 11, the gamma0 and gamma1 cusp forms are identical, but in general the gamma1 space has a higher dimension. For example there is only one elliptic *L*-function of conductor 5077 (see fig 13), but the gamm0 space has dimension 423 and the gamma 1 space dimension 1076535. If the weight is increased from 2 to 12 the dimension is even higher 11816505! (Right) Symmetric Square *L*-Functions of modular forms over $SL(2,Z)$ of weight $k = 12, 16, 20$ (Dummigan), having weight $4k-3$, with five Langlands gamma parameters $[0, 1, 1-k, 2-k, 2-2k]$. (Inset) $k=12$ negative real zeros showing varying degrees of degeneracy (rotated).

Since $f(z+1)=f(z)$, *f* is periodic, we can express it we can express it as a Fourier series in *z* or a Laurent series in *q* $f(z)=\sum_{n=-\infty}^{\infty}a_n e^{2\pi inz}=\sum_{n=-\infty}^{\infty}a_n q^n$. If *f* is meromorphic and has only simple poles we have only a finite number of negative powers of *q* and if *f* is holomorphic, we have a Taylor expansion $f(z)=\sum_{n=0}^{\infty}a_n e^{2\pi inz}=\sum_{n=0}^{\infty}a_n q^n$, $q=e^{2\pi inz}$. Using the Mellin transform $M(f,z)=\int_0^{\infty}f(t)t^{z-1}dt$, we can derive the *L*-function $L(f,z)=(2\pi)^2 M(f,z)/\Gamma(z)=\sum_{n=1}^{\infty}a_n n^{-z}$, which again has a functional equation. If $L^*(f,z)=N^{z/2}(2\pi)^{-z}\Gamma(z)L(f,z)$, then $L^*(f,z)=(-1)^{w/2}L^*(f,w-z)$, and $L^*$ is meromorphic on *C*.

In the case of weight $w = 2$ there is thus a correspondence between the functional equations of elliptic curves and modular forms. The Taniyama-Shimura modularity theorem asserts that every elliptic curve over *Q* has a modular form parametrization based on the conductor, essentially through the periodicities induced by its toroidal embedding, a relationship pivotal in the proof of Fermat's last theorem (Daney), where Andrew Wiles (1995) showed that any semi-stable elliptic curve (having only multiplicative bad reductions) is modular. But if we can find $x^n + y^n = z^n$ then the elliptic curve $Y^2 = X(X - x^n)(X + y^n)$ is semi-stable but not modular. Hence the proof!

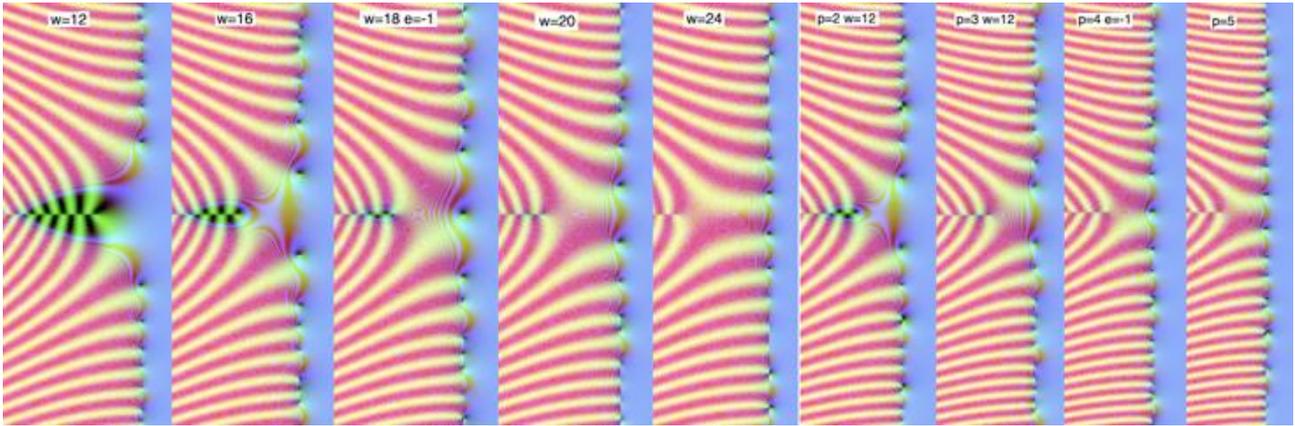

Fig 15a: (Left) Evolution of the principal Modular forms over *SL*(2,*Z*) with increasing weight. (Right) Forms over $\Gamma_0(N)$ for increasing levels *N* with weight 12 have a similar evolution, with increasing dimensions of old and new forms.

We can find the modular form corresponding to a given elliptic curve as follows (Lozano-Robledo). Consider the modular group and congruence subgroups:

$$SL(2,\mathbb{Z}) = \left\{ \begin{pmatrix} a & b \\ c & d \end{pmatrix} : a,b,c,d \in \mathbb{Z},\ ad-bc=1 \right\},$$

$$\Gamma_0(N) = \left\{ \begin{pmatrix} a & b \\ c & d \end{pmatrix} \in SL(2,\mathbb{Z}) : c \equiv 0 \bmod N \right\},$$

$$\Gamma_1(N) = \left\{ \begin{pmatrix} a & b \\ c & d \end{pmatrix} \in \Gamma_0(N) : a \equiv d \equiv 1 \bmod N \right\},$$

We now consider modular forms over congruence subgroups of [SL(2,Z)](#) as above. Note that $\Gamma(N) \subset \Gamma_1(N) \subset \Gamma_0(N) \subset SL(2,Z)$, (where $\Gamma(N)$ also has $b \equiv 0$), so that a form in $\Gamma_0(N)$ is also in $\Gamma_1(N)$. For any congruence subgroup $\Gamma_j$ there exists *N* so that $\Gamma_j(N) \subset \Gamma_j$ putting the form in $\Gamma_j(N)$. Since $M \setminus N \Rightarrow \Gamma_j(N) \subset \Gamma_j(M)$ a form in $\Gamma_j(M)$ is also in $\Gamma_j(N)$.

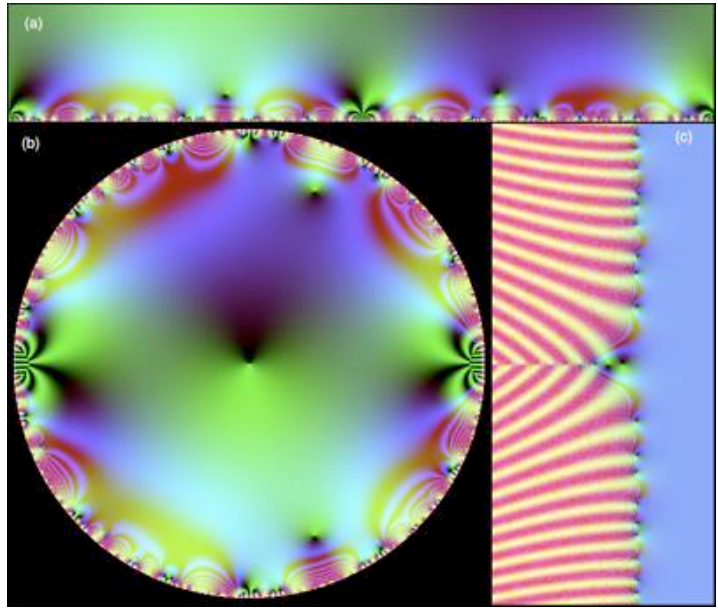

Fig 15b: (a) Modular cusp form over $S_2(\Gamma_0(26))$ represented as a Fourier series in the upper half-plane (a), as a function of *q* in the unit disc (b) and as a Dirichlet series *L*-function with a functional equation in the complex plane (c).

Setting $M_k(\Gamma_j(N))$ for the vector space of weight *k* modular forms and $S_k(\Gamma_j(N))$ for the subspace of cusp forms, we find that $M_2(\Gamma_0(11))$ is two dimensional and $S_2(\Gamma_0(11))$ is one-dimensional, generated by the function *f* with Taylor series in *q* having coefficients $a_n$={1, -2, 1, 2, 1, 2, 2, …} coinciding with those of *e*11*a*. The corresponding situation for $S_2(\Gamma_0(37))$ is a little more complicated, with *M* being three-dimensional generated by: $f(q) = q+q^3-2q^4-q^7-2q^9+\ldots$, $g(q) = q^2+2q^3-2q^4+q^5-3q^6+\ldots$ and $h(q)=1+2q/3+2q^2+8q^3/3+\ldots$, and *S* being two-dimensional, generated by *f*, *g* with corresponding attached *L*-functions as shown in fig 15. The dimension corresponds to the genus of a multi-hole torus embedding (Stein 2008). Turning to *e*37*a*, and *e*37*b* with coefficients *a*={1, -2, -3, 2, -2, 6, -1, 0, 6, …} and *b*={1, 0, 1, -2, 0, 0, -1, 0, -1, …}, we find that $b = f$ and $a = f - 2g$, confirmed by the Taniyama-Shimura theorem, noting that linear combinations of modular forms are modular. Notice that 37*b* has rank 0 with $\varepsilon = 1$ while 37*a* has rank 1, with $\varepsilon = -1$. $S_2(\Gamma_0(N))$ is the direct sum of two subspaces $S^+$ and $S^-$ because the linear operator $w_N$ on $S_k(\Gamma_0(N))$ $w_N(f)(z) = i^k N^{-k/2} z^{-k} f(-(Nz)^{-1})$ is self-dual and thus has eigenvalues ±1.

The $w_n$ eigenfunctions possess the same type of functional equation as elliptic curve $L$-functions. In this case we have an eigenform basis $a$ and $b$, however in the $q^n$ echelon basis generated by Sage, $g$ lies in neither subspace and has a composite functional equation (fig 15).

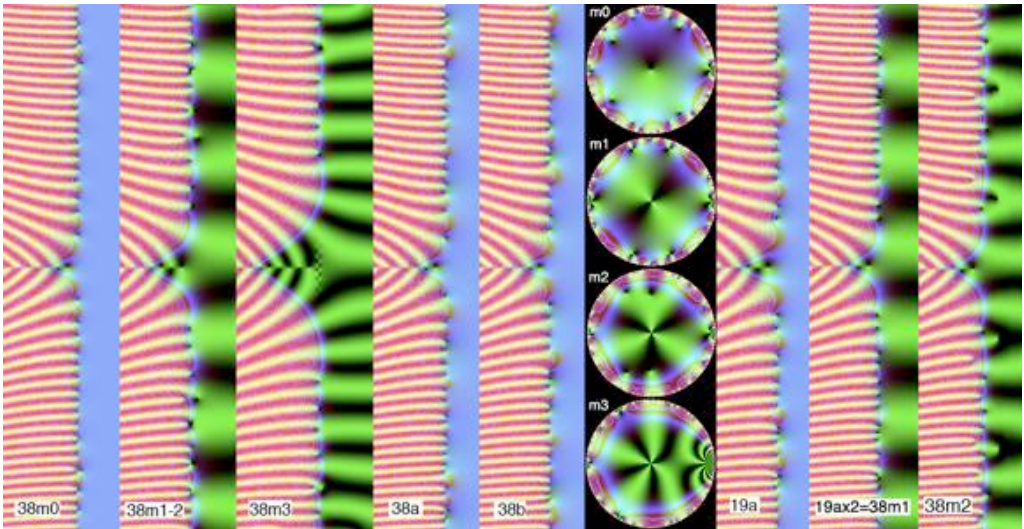

Fig 15c: (Left) $S_2(\Gamma_0(38))$ has basis vectors $m_0$-$m_3$, in echelon-form in $q^n$, e.g. $m_0=q-q^5-2q^6-q^7+\ldots$, $m_1=q^2-2q^6-2q^8+\ldots$, $m_2=q^3+q^5-2q^6+\ldots$, $m_3=q^4-3q^5+q^6+\ldots$, with the elliptic curves being combinations $a=m_0-(m_1-m_2)+m_3$ and $b=m_0+(m_1-m_2)+m_3$. In this case both elliptic curves have $\varepsilon=+1$. (Centre) the power series of $m_0$-$m_3$ in terms of $q$ in the unit disc. (Right) There are also two old eigenforms comprising $e19$ and $2^{-z}L(e19)=m_1$. The linear combination $m19=e19=m_0-2m_2-2m_3$, illustrates the fact that the space of modular forms of level $N$ includes those of $M$: $M|N$. We also have $e19x2=m_1$ giving 4 eigenforms, so can perform a functional equation reconstruction of the four $m_i$ by inverting the matrix defining the elliptic curve eigenforms in terms of the $q^n$ echelon basis. Although each of the $m$ functions shown left also have functional equations with $\varepsilon=+1$, $m_2$ lies in neither $S^+$ nor $S^-$, for $N=19$ or $38$, as with $g$ of $S_2(\Gamma_0(37))$. Again the $m_i$, being superpositions of eigenforms appear to have off-critical zeros, while the eigenforms $a$, $b$, $e19$ and $e19x2$ do not.

When we have a non-prime level $N$, there are both new forms and a spectrum of old forms arising from each of the factors of $N$. For each factor $M$, and $d|(N/M)$ we have the modular function $f(q^d)$ as well as $f(q)$ of level $M$. For example in the case of $N=38$, there are two old forms, which in terms of the Hecke operators, (see below) are eigenforms. These are the elliptic curve $L$-function $e19$ and the related $L$-function $2^{-z}L(e19)$, viz: $L_{e19}(z) = 1^{-z} + 0 \cdot 2^{-z} - 2 \cdot 3^{-z} + \cdots$, with $g(q) = q - 2q^3 + \cdots$. We also have $g(q^2) = q^2 - 2q^6 + \cdots$, with $L$-function $L_{e19x2}(z) = 2^{-z} + 0 \cdot 4^{-z} - 2 \cdot 6^{-z} + \cdots = 2^{-z}L_{e19}(z)$.

At another extreme, the space of modular forms over $SL(2,Z)$ of weight 12 has 2 dimensions, with basis vectors represented by the form $\Delta$ of the $\tau$ function illustrated in fig 11 and the normalized [Eisenstein series](#) $E_{12}$ where $E_{2k} = 1 + \dfrac{2}{\zeta(1-2k)}\sum_{n=1}^{\infty}\sigma(n,2k-1)q^n$, each of which is a modular form of weight $2k$ over $SL(2,Z)$. Eisenstein series are defined by $G_{2k}(z) = \sum_{m,n \in Z^2 \setminus (0,0)} (m+nz)^{-1}$ with $q = e^{2\pi i n z}$ as above, with further generalizations to $m, n \equiv 0 \mod N$ for $\Gamma_j(N)$. This doesn't have an $L$-function because the 1 makes it not a cusp form, but the coefficients generated by the divisor function do coincide with the sigma function $\zeta(z)\zeta(z-(2k-1))$ of fig 21, thus giving an illustration of the two eigenforms. In many ways the modular forms of weight 2 are atypical as Eisenstein series only begin with weight 4.

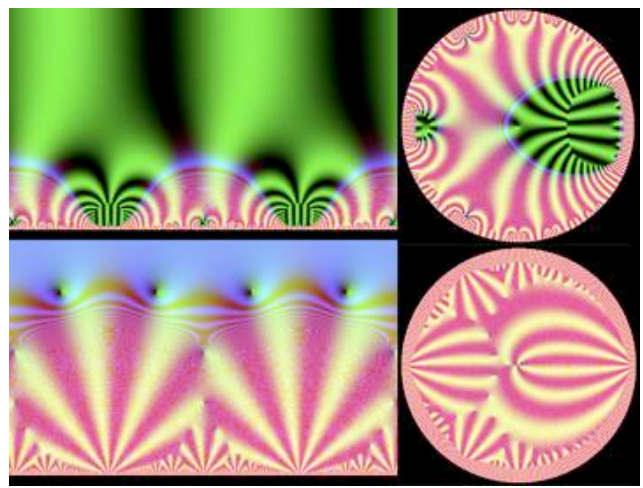

Fig 15d: Fourier and Taylor representations of the Tau and Eisenstein functions of weight 12 in addition to the $L$-function portraits of figs 11 and 21.

The symmetric square lift (see fig 15) is defined as follows. Given a form $f$ over $SL(2,Z)$ with Euler product, $L(z,f) = \prod_{p \text{ prime}} \left(1 - a_p p^{-z} + p^{-2z}\right)^{-1} = \prod_{p \text{ prime}} \left(1 - \alpha_f(p) p^{-z}\right)^{-1} \left(1 - \alpha_f(p)^{-1} p^{-z}\right)^{-1}$ where $a_p = \alpha_f(p) + \alpha_f(p)^{-1}$, the symmetric square lift is a $GL(3)$ form $\phi$ with Euler product

$$L(z,\phi) = \prod_{p \text{ prime}} \left(1 - A(1,p) p^{-z} + \overline{A(p,1)} p^{-2z} - p^{-3z}\right)^{-1} = \prod_{p \text{ prime}} \left(1 - \alpha_\phi(p) p^{-z}\right)^{-1} \left(1 - \beta_\phi(p) p^{-z}\right)^{-1} \left(1 - \gamma_\phi(p) p^{-z}\right)^{-1}$$

where $\begin{pmatrix} \alpha_\phi(p) & & \\ & \beta_\phi(p) & \\ & & \gamma_\phi(p) \end{pmatrix} = \begin{pmatrix} \alpha_f(p)^2 & & \\ & 1 & \\ & & \alpha_f(p)^{-2} \end{pmatrix}$ (Dummigan, Bian).

Each of the types of $L$-function discussed admit a functional equation determined by the Dirichlet series, a finite number of gamma (Hodge, or Langlands) parameters, determined by the underlying topology generating the Euler product, the conductor, and a sign factor (Dokchitser, Harron):

$$L^*(f,z) = N^{z/2} (2\pi)^{-z} \Gamma\left(\frac{z+\lambda_1}{2}\right) \cdots \Gamma\left(\frac{z+\lambda_d}{2}\right) L(f,z), \text{ then } L^*(f,z) = \varepsilon L^*(f, w-z)$$

where $|\varepsilon| = 1$, $\varepsilon = e^{2\pi i k/n}$ for Dirichlet $L$-functions $\varepsilon = \pm 1$ otherwise. The gamma factors can be used to define a generalized Mellin transform technique for describing $L$-functions in the critical strip for moderate $y$ values (Dokchitser). These types can be generalized in motivic $L$-functions (Deligne, Dokchitser). The Langlands program (1980) of automorphic forms includes a comparable explanation of Euler products involving polynomials of higher degree.

Modular forms that are eigenfunctions of all Hecke operators $T_n f = \lambda_n f = a_n f$, where

$$T_n f(z) = n^{k-1} \sum_{M \in \Gamma \backslash M_n} (cz+d)^{-k} f\left(\frac{az+b}{cz+d}\right), M_n = \left\{ A = \begin{pmatrix} a & b \\ c & d \end{pmatrix} : |A| = n \right\},$$

have the equivalent Euler product $L(f,z) = \prod_{p | N} (1 - a_p p^{-z})^{-1} \prod_{p \nmid N} (1 - a_p p^{-z} + p^{2k-1-2z})^{-1}$ (Cogdell, Lozano-Robledo).

Conveniently each eigenfunction satisfies all the Hecke operators $T_n$ simultaneously.

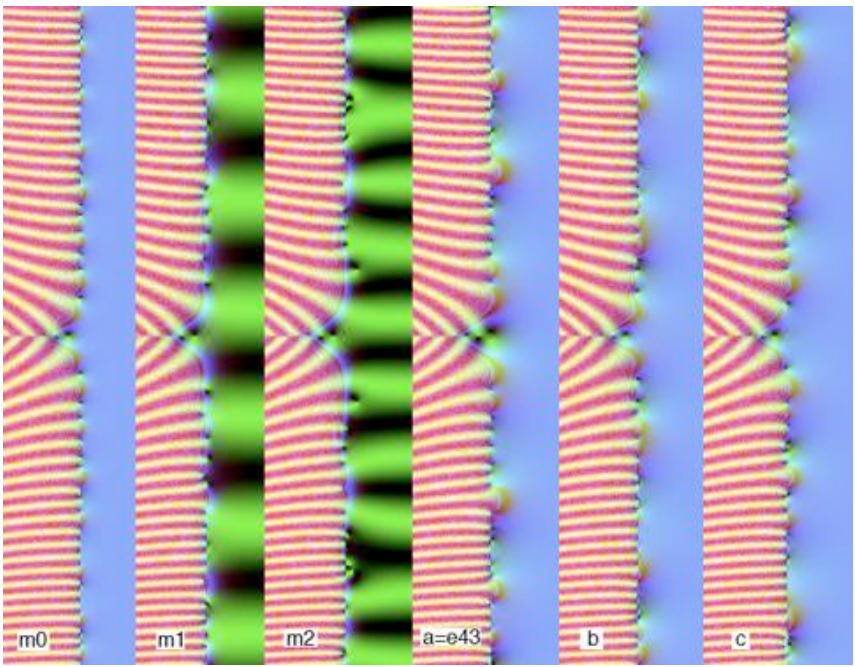

Fig 15e: The situation for level $N=43$ presents new features. The cusp space has three dimensions and there are also three eigenfunction newforms $a$, $b$, $c$, the first of which is identical to elliptic e43, which has root -1. The two additional root 1 eigenforms have Hecke eigenvalues $\pm\sqrt{2}$, and do not correspond to elliptic curves. Hence there are sufficient eignfunctions to represent the $q$-echelon basis forms in terms of the combined eigenfunction equations (left).

One can calculate Hecke operators and matrices in terms of the power series $f(q) = \sum_m a_m q^m$ as follows:

$T_p(f(q)) = \sum_{m \in Z} (b_m q^m = \sum_{m \in Z} (a_{mp} + p^{k-1} a_{m/p}) q^m$, $a_{m/p} = 0$, $m/p \notin Z$, $p$ prime, which is generally sufficient for determining eigenvalues and eigenforms.

For example for $N=43$, we have echelon basis: $f = q + 2*q^5 - 2*q^6 - 2*q^7 - q^9 + O(q^{10})$, $g = q^2 - 1/2*q^4 + q^5 - 3/2*q^6 - q^8 - 1/2*q^9 + O(q^{10})$, $h = q^3 - 1/2*q^4 + 2*q^5 - 3/2*q^6 - q^7 + q^8 - 1/2*q^9 + O(q^{10})$
Applying the above formula for $T_2$ to $f$, we have $b_1 = a_2 + 2.0 = 0$, $b_2 = a_4 + 2.0 = 2$,

$b_3 = a_6 + 2.0 = -2$, giving the first row of the Hecke matrix $T_2 = \begin{bmatrix} 0 & 2 & -2 \\ 1 & -1/2 & -3/2 \\ 0 & -1/2 & -3/2 \end{bmatrix}$, which has eigenvalues $a_0 = -2$ and $a_1 = \pm\sqrt{2}$. These eigenvalues give normalized eigenvectors $v = q - 2*q^2 - 2*q^3 + 2*q^4 - 4*q^5 + 4*q^6 + q^9 + O(q^{10})$, $w_{a1} = q + a_1*q^2 - a_1*q^3 + (2-a_1)*q^5 - 2*q^6 + (a_1 - 2)*q^7 - 2*a_1*q^8 - q^9 + O(q^{10})$, which are *newforms*, normalized eigenforms of level $N$ not arising from an $M < N$, defined as a linear combination of the above echelon basis functions, the first of which is the elliptic curve $e43a$. One can confirm they are eigenvectors using the same formula, e.g. $T_2(v) = -2v$.

By inverting the resulting basis transformation matrix $A = \begin{bmatrix} 1 & -2 & -2 \\ 1 & \sqrt{2} & -\sqrt{2} \\ 1 & -\sqrt{2} & \sqrt{2} \end{bmatrix}$, we can in turn express the echelon basis in terms of the eigenforms. Each elliptic function with conductor $N$ is thus associated with a *newform* - (Stein).

Computing the original echelon basis is more complicated (Stein). For forms over $SL(2,Z)$ such as those of weight 24 in fig 15a, we can derive a basis for the three dimensional modular space based on Eisenstein series namely: $E_4^6$, $E_4^3 E_6^2$, $E_6^4$ and then perform row operations to gain a $q^n$ echelon basis called the Miller basis, which has two dimensions of cusp forms. Computing the bases of cusp forms over $S_2(\Gamma_0(N))$ is complicated and most conveniently done using modular symbols, which are representations of homology classes of paths on the embedded multi-hole torus whose genus determines the dimension of the modular space, represented on the upper half plane between rational points on the real line (including the point at $\infty$). Modular symbols have relations such as $\{\alpha,\beta\} + \{\beta,\gamma\} + \{\gamma,\alpha\} = 0$, $\{\alpha,\alpha\} = 0$, $\{\alpha,\beta\} = -\{\beta,\alpha\}$ and are acted on by rational matrices $g\{\alpha,\beta\} = \{g(\alpha), g(\beta)\}$, can be readily computed using Sage by a technique derived by Manin using continued fractions, and are compatible with Hecke operators (Stein).

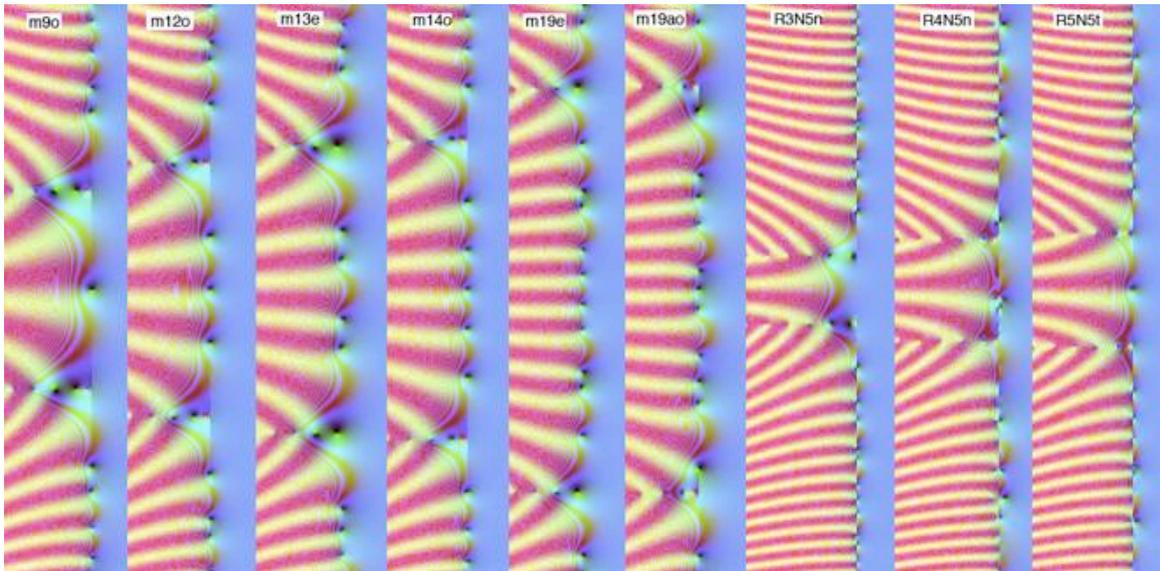

Fig 15f: *L*-functions of Maass forms over $PSL(2,Z)$ with eigenvalues 9.5336, 12.1730, 14.3585, 19.4847 (odd) and 13.7797, 19.4234 (even) and three forms over $\Gamma_0(5)$ with eigenvalues 3.2642, 4.8937, 5.4361.

Maass forms are modular differential functions satisfying the hyperbolic Laplace wave function $\Delta = -y^2\left(\frac{\partial^2}{\partial x^2} + \frac{\partial^2}{\partial y^2}\right)$, which commutes with SL(2,Z), and generates a vast spectrum of eigenforms, having complex gamma factors $\lambda_i = e \pm ir$, $i = 1,2$, where $e=0,1$ and $\varepsilon = 1,-1$ for even and odd functions respectively where the eigenvalue is $\frac{1}{4}+r^2$, and a slightly more complicated Fourier series $f(z) = \sqrt{y}\sum_{n=1}^{\infty} a_n K_{ir}(2\pi |n| y) e^{2\pi i x}$, with $K_{ir}$ the modified Bessel function. For N=11 there are around 1000 such forms over $\Gamma_0$ (Booker et. al. 2006, Farmer and Lemurell), which can be located by searching for eigenvalue hot spots. Several Maass form L-functions are illustrated in fig 15f.

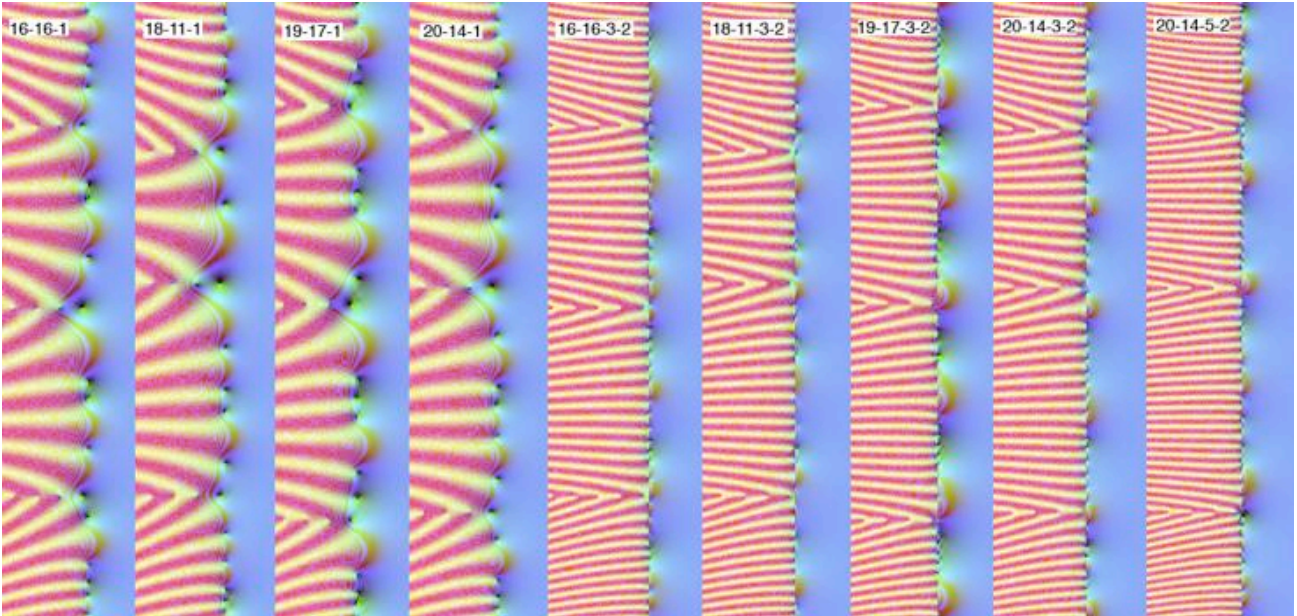

Fig 15g: Third degree transcendental Maass form L-functions. First views of the L-function profiles for the above parameters for three selections of Dirichlet character.

A new class of L-function (Bian 2010, Booker 2008) has been discovered, based on automorphic GL(3) Maass forms, which are written in terms of a three dimensional generalized upper half-plane w=XY, where $X = \begin{pmatrix} 1 & x_2 & x_3 \\ 0 & 1 & x_1 \\ 0 & 0 & 1 \end{pmatrix}$, $Y = \begin{pmatrix} y_1 y_2 & 0 & 0 \\ 0 & y_1 & 0 \\ 0 & 0 & 1 \end{pmatrix}$, $x_i, y_i \in R$, $y_i > 0$. The form $\varphi(w)$ is an eigenfunction of the Laplacian, which is preserved under SL(3,Z). This has an extended Fourier series, which can be used to define a complex L-function with a degree 3 Euler product

$$L(z, \varphi \times \chi) = \prod_{p \text{ prime}} \left(1 - A(1,p)\chi(p)p^{-z} + A(p,1)\chi^2(p)p^{-2z} - \chi^3(p)p^{-3z}\right)^{-1}$$

where $\chi(p)$ is a Dirichlet character 'twisting' the L-function and A(p,q) are Fourier coefficients of the Maass cusp form with eigenvalues $(\lambda_1, \lambda_2)$, $\text{real}(\lambda_i) = 1/3$.

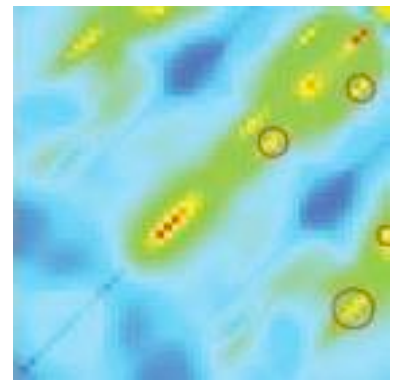

Fig 15h: Hot regions in (u,v) =[10, 20]2 (red) with 4 non-trivial degree 3 examples circled (18.902, 11.761), (16.741, 16.232), (20.021, 14.070), (19.179, 17.702) and quasi-trivial examples on the diagonal u = v at 13.779, 17.738, 19.423 (Bian)

In particular $A(1,p) = \overline{A(p,1)}$. Bian has found locations in parameter space (u,v), $\lambda_1 = (1+ui)/3, \lambda_2 = (1+vi)/3$ defined in terms of the gamma function imaginary parameters

$\alpha = -(u+2v)i/3$, $\beta = (2u+v)i/3$, $\gamma = (-u+v)i/3$ where non-trivial transcendental degree 3 $L$-functions, not being simply a product of degree 1 or degree 1 and 2 Euler products, nor symmetric square lifts (fig 15) of a quadratic Euler function, occur. These also obey a functional equation of the form $\Lambda(z, \varphi \times \chi) = \varepsilon_\chi^3 \Lambda(1-z, \tilde{\varphi} \times \bar{\chi})$, where ~ takes coefficients and gamma factors to their conjugates, and $\Lambda(z, \varphi \times \chi) = N^{z/2} \pi^{-(z-\alpha)/2} \Gamma\left(\frac{z-\alpha}{2}\right) \pi^{-(z-\beta)/2} \Gamma\left(\frac{z-\beta}{2}\right) \pi^{-(z-\gamma)/2} \Gamma\left(\frac{z-\gamma}{2}\right) L(z, \varphi \times \chi)$.

The upshot of this study of a reasonable spread of $L$-functions and non-$L$ counterparts, is that the non-trivial zeros lie on the weighted critical line only if they are generated by an underlying non-mode-locked primal distribution, despite the fact that the Dirichlet sum is over all integers and the relationship with the Euler product over primes holds only outside the critical strip in the right half-plane. This suggests widening the approach to consider more general classes of Euler products.

**Functions with Functional Equations but no Euler Product**
Davenport and Heilbronn (1936) devised an example of a zeta function possessing a functional equation (Bombieri and Ghosh, Titchmarsh) but no Euler product, with a set of non-trivial zeros wide of the critical line. Given $\xi = -\phi + \sqrt{1+\phi^2} = 0.2841$, $\phi = (1+\sqrt{5})/2$, the period 5 Dirichlet pseudo-character $\chi(5,1) = \{0, 1, \xi, -\xi, -1\}$ gives rise to a Dirichlet series having functional equation $L_{5,1}(z) = 5^{1/2-z} 2(2\pi)^{z-1} \Gamma(1-z) \cos(\pi z/2) L_{5,1}(1-z)$ (see appendix 3 for derivation).

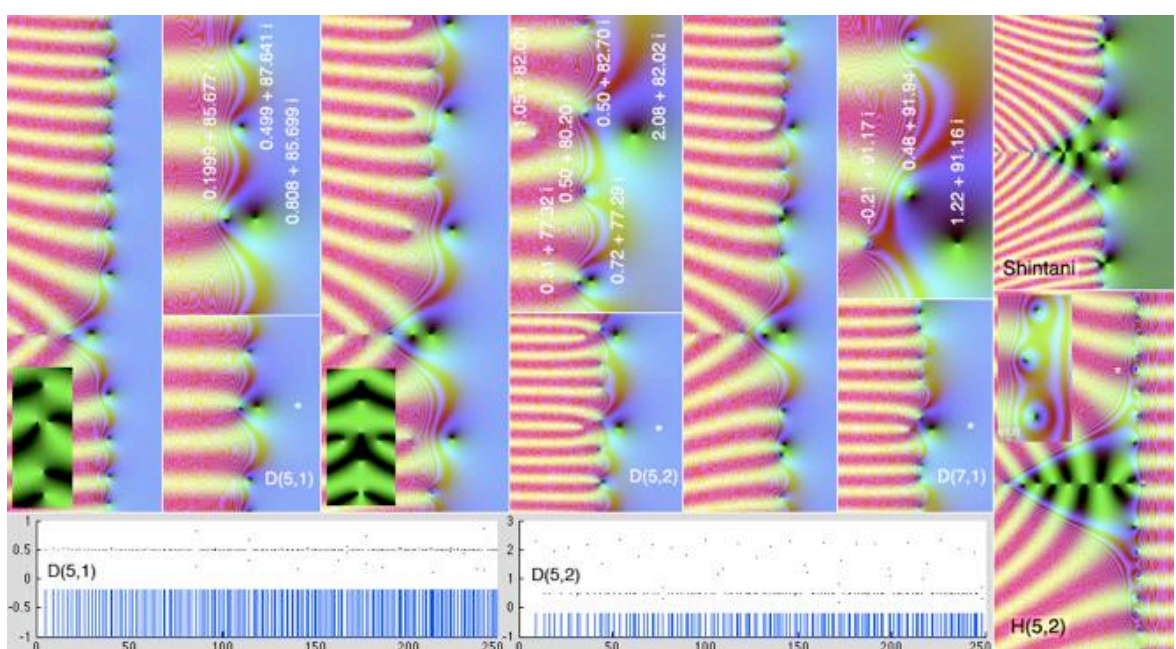

Fig 16: Davenport-Heilbronn functions (5,1), (5,2), and (7,1) possess functional equations demonstrating they are meromorphic on the complex plane, but lack an Euler product and have an array of non-trivial zeros off the critical line showing a functional equation is insufficient for RH. Significant is the large number of apparently critical zeros, with a sparse spread of off-critical ones in pairs, which shows why many critical zeros don't necessarily imply all. The symmetrical placing of the zeros about the critical line in (5,1) and (5,2) is confirmed in the local Xi function portraits inset. Intriguingly the majority of non-twin zeros appear to be on the critical line for the parameter value $\xi$, although all values have a quasi-functional equation expressed in the sum of DL-functions (appendix 3), so the symmetric form of the functional equation is a key. Hurwitz zeta (5/2) and a Shintani zeta function (right) likewise show sums of L-functions can have off-critical zeros.

This example was extended by Balanzario and Sanchez-Ortiz (2007) to a small class of functions, which likewise serve to demonstrate the existence of a functional equation is not a sufficient condition for non-trivial zeros to be critical. Out of these, two further examples are not simply derived Riemann zeta and Dirichlet $L$-functions and lack an Euler product, the pseudo-character

χ(5,2) ={0, 1, -1/ξ, 1/ξ, -1} and χ(7,1) ={0,1, -(1+α), -α, α, 1+α, -1}, where α~ 0.80194. These also give rich examples of off-critical zeros in a context where a valid functional equation allows us to accept all displayed zeros, including those in the half-plane $x<0$ are genuine. Intriguingly, there is a mix of a large number of critical and a sparse number of off-critical zeros defining a test case.

The Hurwitz zeta functions $\zeta(z,a) = \sum_{n=0}^{\infty}(n+a)^{-z}$ and their generalization in Shintani zeta functions also have off-line zeros with real values approaching every $0<x<1$ for some zero, despite having a functional equation, for rational $p/q$ in $(0,1)$

$$\zeta(1-z, p/q) = 2\Gamma(z)(2\pi q)^{-z}\sum_{k=1}^{q}\cos\left(\frac{\pi z}{2} - \frac{2\pi kp}{q}\right)\cdot\zeta(z, k/q),$$ when they are a sum over the

Dirichlet $L$-functions of period $q$, $\zeta(z, p/q) = q^z/m \sum_{\chi(q,k)\neq 0}\bar{\chi}(q,k)L(z,q,k)$, $m = |\{\chi(q,k) \neq 0\}|$, illustrating that sums of $DL$-functions are not necessary $L$. The original Davenport-Heilbronn example is derived as a sum of $DL$-functions through being a sum of Hurwitz zeta functions: $f(z) = 5^{-z}\left(\zeta(z,1/5) + \xi\cdot\zeta(z,2/5) - \xi\cdot\zeta(z,3/5) - \zeta(z,4/5)\right)$, which brings it very close to the period 5 $DL$-function with character $\{0, 1, i, -i, -1\}$ since the Dirichlet L-function

$L(z,q,k) = q^{-z}\sum_{k=1}^{q}\chi(k)\zeta(z, k/q)$. In fact the Hurwitz zetas have simple arithmetic coefficient sequences e.g. $\zeta(z,1/3) = (0+1/3)^{-z} + (1+1/3)^{-z} + (2+1/3)^{-z} + \cdots = 3^z(1 + 4^{-z} + 7^{-z}\cdots)$.

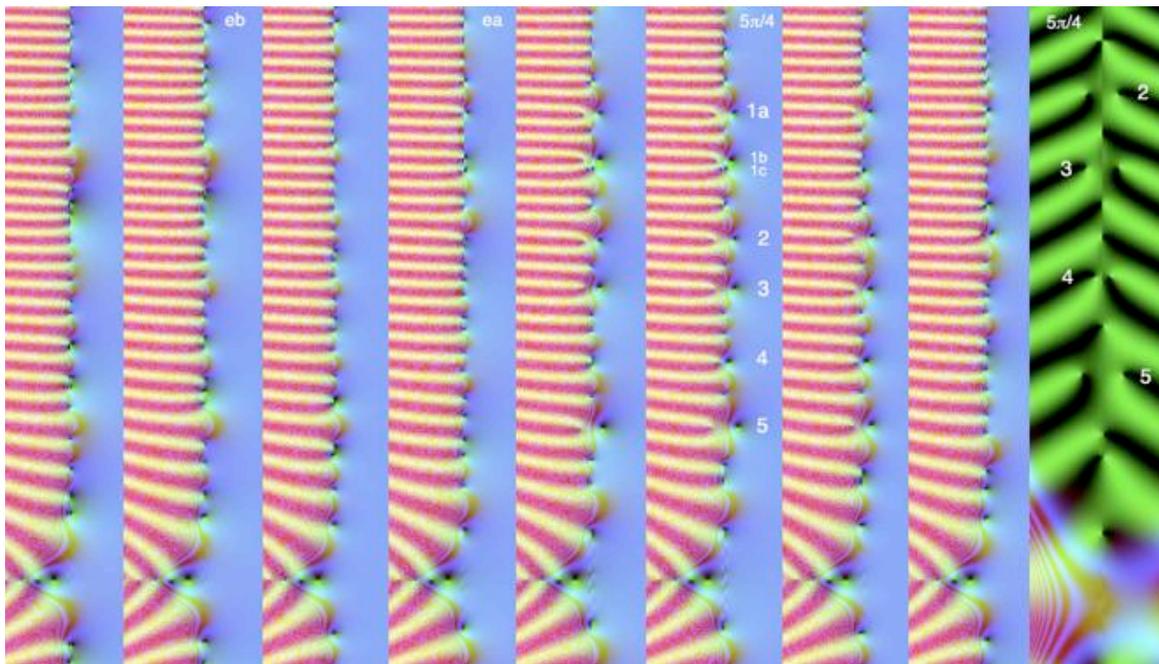

Fig 16a: Collisions and splitting into symmetrical off-critical zero pairs in a circular journey in the space of modular forms $S_2(\Gamma_0(38))$ illustrated in fig 15c. Far right Xi function confirming functional equation symmetry of the zero pairs.

**A Central Showcase: Modular Forms Meeting Elliptic Functions**
Let us return to the modular forms of fig 15c with $N=26$, an example, which shows the dynamical basis explicitly and makes the most sensitive test concerning how variations in the coefficients effect whether the zeros remain on the critical line. If we consider the four basis functions illustrated in fig 15c, these are echelon in powers of $q$, but are not necessarily Hecke eigenforms. In fact the two elliptic curve $L$-functions $e38a$ and $e38b$ as newforms of level $N=38$ and the elliptic curve $L$-function $e19$ of level $N=19$ stand as three of four basis eigenfunctions in the four-dimensional modular vector space. So what happens if we make a journey in the vector space of modular forms and examine the resulting '*L*-functions'?

In fig 16a are shown the series consisting of
$$m_0 + \sqrt{2}\cos(2\pi\theta)(m_1 - m_2) + \sqrt{2}\sin(2\pi\theta)(m_3), \ \theta = k/8, \ k = 0,\ldots,7.$$
For $k = 1, 3$ these coincide with the two elliptic curve *L*-functions, however as we circulate, we find that the critical zeros move down the critical line and successively collide to form symmetrical pairs of off-critical zeros, as confirmed by the Xi portrait. This is consistent with these non-eigenform *L*-functions again having a functional equation but not necessarily an Euler product, since even if the original basis functions did have an Euler product (and from fig 15c we can see that neither $m_1$-$m_2$ nor $m_3$ appear to) as we saw with the Davenport-Heilbronn examples, there is no reason why an arbitrary linear combination of such function should have a product representation. The uppermost of the collisions occurs very close to *e*38*a* This would imply that an arbitrarily small deviation from an eigenform with an Euler product would lead to a collision and violation of RH somewhere on the critical line. Notice also that even though *e*38*a* and *e*38*b* correspond to modular eigenforms of *N*=38, as noted and have the same functional equation, linear combinations of them do not necessarily satisfy RH (e.g. $b-a=2(m_1-m_2)$, which from fig 15c clearly doesn't).

This behavior is relatively tame in that the critical strip zeros appear to move along the line until they collide and form symmetrical off-line pairs, consistent with the functional equation. As noted above, despite the tameness, we would expect a violation of RH somewhere on the critical line for any arbitrarily small deviation from an eigenform.

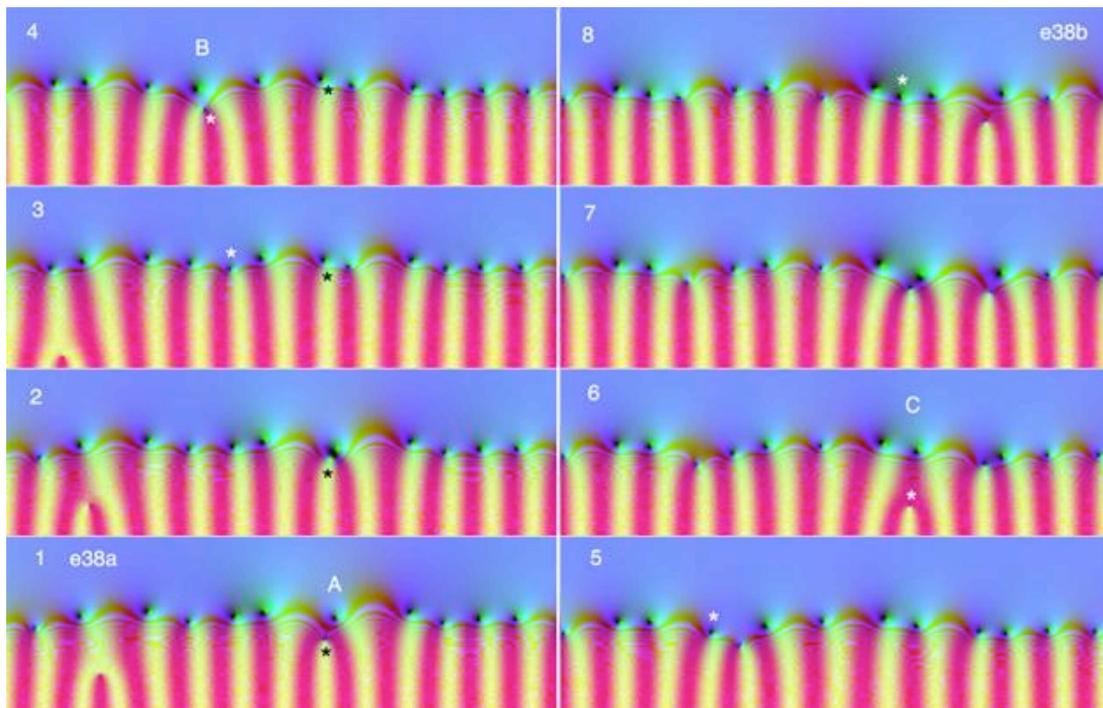

Fig 16b: Threaded dynamics in continuous variation of Euler products between elliptic curve *L*-functions *e*38*a* and *e*38*b*. (A) A Gollum zero has a close encounter with a critical zero before taking a position as a critical zero beside it. (B) A critical zero descends below its left-hand neighbour and rejoins the critical line on the other side. (C) The Gollum zero which began at the far left moves to the centre and enters the critical line.

Supposing now we reverse the situation and demand that the varying functions do have an Euler product. Rather than forming a linear combination of *m*0, *m*1-*m*2 and *m*3, we simply take a linear combination of the prime coefficients $a_p$ forming the Euler products of the elliptic curve *L*-functions *e*38*a* and *e*38*b*, which are themselves linear combinations $c_p = ta_p + (1-t)b_p$ of the modular forms and eigenfunctions. We can then regenerate the appropriate Dirichlet series coefficients using the recurrence realtion $c_{p^e} = c_p c_{p^{e-1}} - (p \mid N) p c_{p^{e-2}}$, $c_n = \prod c_{p_i^{e_i}}$, $n = p_1^{e_1} \cdots p_k^{e_k}$, $k > 1$. Each of these will have an Euler product, but not necessarily an exact functional equation, since they are no longer a linear combination of *e*38a and *e*38b. Significantly the two *L*-functions and their intermediates have good convergence in the critical strip by comparison with several examples that follow.

When we do this, we find that some very bizarre dynamics have set in. We will use a 'naked' portrait, simply analytically continuing the finite approximation to the negative real half-plane because of the lack of an exact functional equation. We now find that the zeros have a dynamic variation, leaving the critical line and exchanging positions both with one another and with shadow or 'gollum' zeros in the negative real half-plane, which would normally be eliminated by the functional equation representation.

This has two implications. It implies that the *L*-functions with critical zeros have special 'quantum' properties equivalent to being eigenfunctions that result in their having both and Euler product and a functional equation and it is only the eigenstate that conforms to RH. Secondly there is a topological issue in the ordering of the zeros, which can result in topological braiding of the zeros under continuous transformation.

**Seeking Examples with Product Formulae**
We now explore functions that do have a product representation to seek further examples outside the class of *L*-functions. Let us first consider the product of integers:

$$f(x) = \prod_{n=2}^{\infty} (1 - n^{-z})^{-1} = \sum_{n=1}^{\infty} \phi(n) n^{-z}, \phi(n) = \text{ unordered factorizations of } n.$$

The number of unordered factorizations of *n* with largest part at most *m* can be calculated from the recursion relation $\phi(m,n) = \sum_{\substack{d \backslash n \\ d \leq m}} \phi(d, n/d)$ (Hughes & Shallit 1983, Knopfmacher & Mays 2006).

The 'unordered' factorizations (in which different orders are not distinct) consist of all possible distinct $n^{-z}$ terms in the product, which become coefficients of a given *n*. For example for 12 we have 4: 12, 6.2, 4.3, 3.2.2, written in descending order of the factors involved. With zeta and the *L*-functions, because prime factorization is unique, there is only one such term, so the coefficients of zeta are all 1. In this case, the coefficients vary widely, according to the distribution shown in fig 17, and we will expect to see significant phase-locking in the imaginary waves.

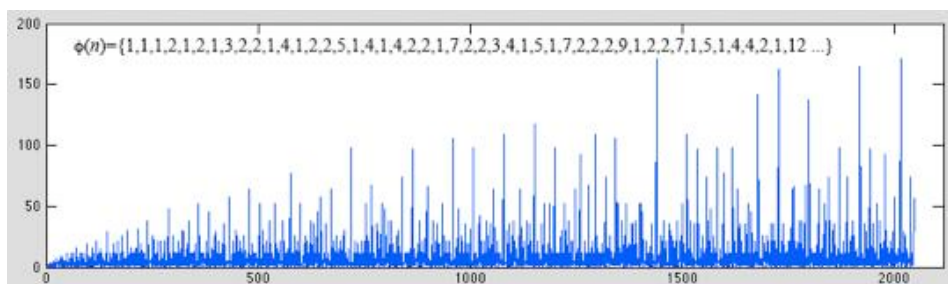
Fig 17: Unordered factorizations as a function of *n*.

To give an exploratory profile based on alternating series, we examine the related functions:
$\sum_{n=1}^{\infty} \phi(n)(-1)^{n-1} n^{-z}$, the alternating variant, the function $\sum_{n=1}^{\infty} (\phi(n) - 2(1 - n \bmod 2)\phi(n/2)) * n^{-z}$ which is derived from that of *f*(*x*) in the same way as eta is derived from zeta by subtracting $2^{1-z}$ times the series from itself, and $(1-2^{1-z})^{-1}\sum_{n=1}^{\infty}(\phi(n)-2(1-n\bmod 2)\phi(n/2))*n^{-z}$, the zeta series re-derived from the previous alternating series. We also for a comparison investigated the series

$$f_{-}(x)=\prod_{n=2}^{\infty}(1+n^{-z})^{-1}=\sum_{n=1}^{\infty}\varphi(n)n^{-z}, \varphi(n)= \text{ unordered factorizations of } n \text{ with alternating powers},$$

viz {1,-1,-1,0,-1,0,-1,-1,0,0,-1,0,-1,0,0,1,-1,0,-1,0,0,0,-1,1,0,0,-1,0,-1,1,-1,-1,0,0,0,1 …} calculated by sorting factorizations into bins, removing duplicates and checking against the above method.

As can be seen from fig 18, both the Newton's method Matlab portrait and the function plot using the software application developed by the author (see below) show quasi-regular variations in the position of the zeros, consistent with substantial phase-locking caused by the fluctuating (repeated) coefficients. The two representations of the factorization zeta function top and second bottom show a degree of consistency, which can be compared with the naked and analytic versions of zeta itself in fig 19. By contrast, the last pair, which end up having quasi-random coefficients close to 0, 1 and -1, the portrait is similar to the random coefficient Dirichlet sequence of fig 5.

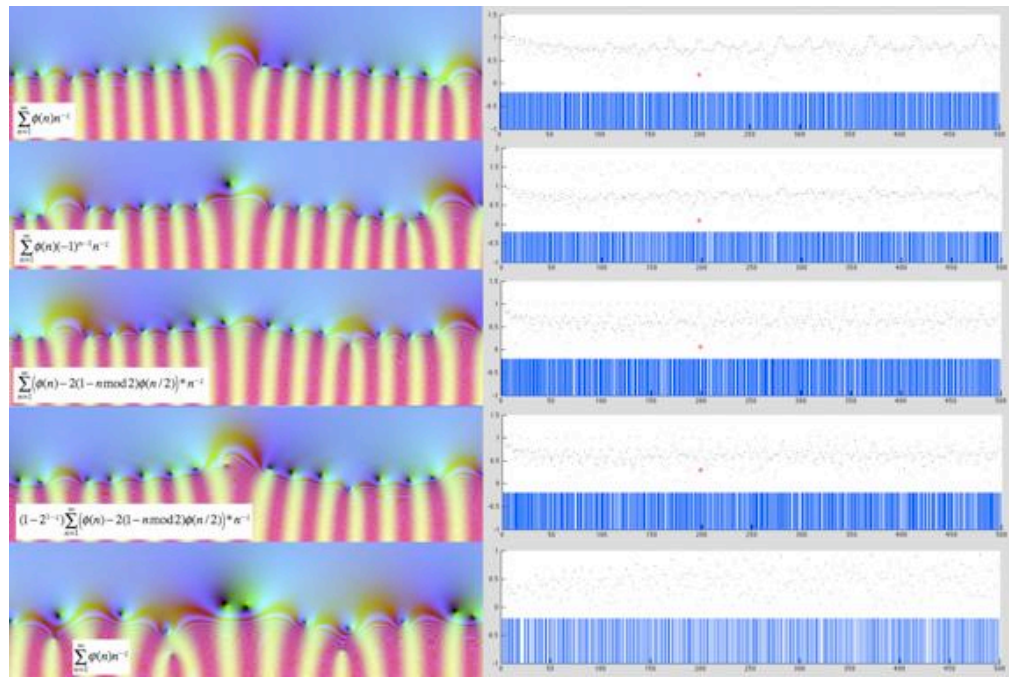

Fig 18: Dirichlet functions derived from an integer rather than a prime product show evidence of phase locking arising from their erratically increasing coefficients. Left naked plots in the region of 200, starred at the right where Newton's method is used to seek for zeros. All are at 1024 series terms except the bottom pair at 250 terms.

We now examine more closely how Euler products of primes with varying coefficients might fare when encoded back into Dirichlet series. We can't take the products directly because these are unstable in the critical strip, and each involves an infinite number of terms in the sum, however, if we define a set of coefficients $\chi(p)$ for each prime then we have for each term

$$\frac{1}{1-\chi(p)p^{-z}}=\sum_{n=0}^{\infty}(\chi(p)p^{-z})^n=1+\chi(p)p^{-z}+(\chi(p))^2(p^2)^{-z}+\cdots$$

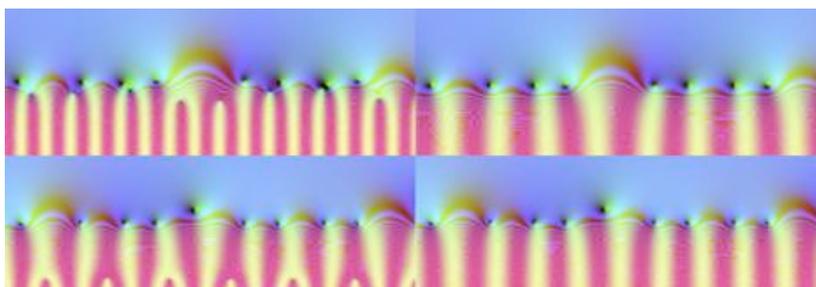

Fig 19: Naked and analytic portraits of zeta (above) and eta (below) show that eta's alternating series naked representation is true to its equivalent representation analytically in the critical strip, while zeta's shows distortions of the zeros caused by non-convergence of the absolute series.

Since each base $n$ in the Dirichlet series has a unique prime factorization $n=p_1^{m_1}\ldots p_k^{m_{1k}}$ it will have a uniquely-defined sum coefficient $\chi(n)=(\chi(p_1))^{m_1}\ldots(\chi(p_k))^{m_k}$, so we can define the series. Any such coefficient set is completely multiplicative, as *L*-function characters are. More generally we could define $\chi(p,n)$ separately for each $p^n$, which would be multiplicative, but not completely.

In fig 20 are shown a set of examples where similar variations of the completely multiplicative coefficients have been chosen to those of the additive coefficients in fig 5. In all cases apart from (b) the zeros are distributed off the critical line, implying an Euler product may not be sufficient to cause the zeros to be on the critical line.

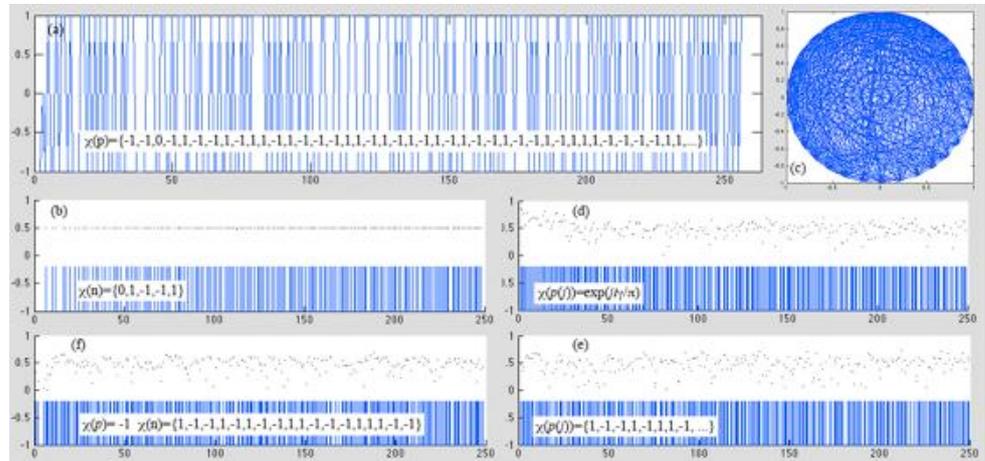

Fig 20: (a) The multiplicative coefficients of *L*(5,3) (b) form an irregular distribution on the primes (d) a golden mean angle variation on the prime product coefficients (c) and a Morse-Thule fractal recursive distribution of {1,-1} have zeros off the critical line, with significant indications of phase-locking as does the lambda function (f) with all prime multiplicative coefficients -1, but this function does have zeros on the critical line x=1/4 through its analytic expression in terms of zeta indicating lack of convergence of the naked approximation.

This brings us to a major caveat about representing approximations as naked functions in the 'forbidden' zone $x < 0$, or even in the critical strip without a guarantee of convergence.

The lambda function $\sum_{n=1}^{\infty}\frac{\lambda(n)}{n^s} = \prod_{p\text{ prime}}\left(1+p^{-s}\right)^{-1}$, where $\lambda(n)=(-1)^{\Omega(n)}$ $\Omega(n)$= number of prime factors of *n* with multiplicity, as shown in fig 20(f), has multiplicative coefficients $\chi(p)=-1$, and can immediately be expressed in terms of the zeta function, since $\dfrac{\prod\limits_{p\text{ prime}}\left(1-p^{-s}\right)}{\prod\limits_{p\text{ prime}}\left(1-(p^2)^{-s}\right)} = \dfrac{\zeta(2s)}{\zeta(s)}$.

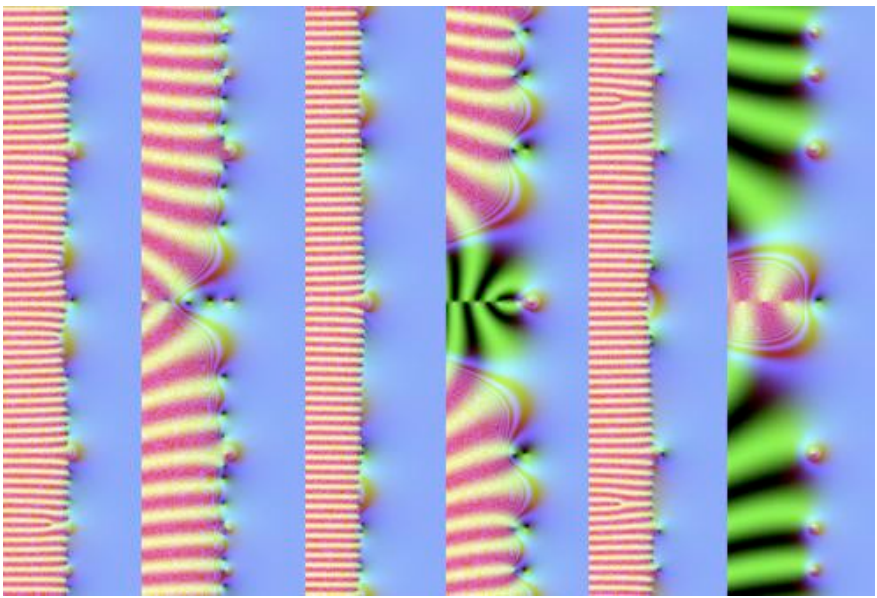

Fig 21: Failure of adequate convergence to the analytic continuation in lambda (left 2) sigma (centre 2) and mu (right 2) at 1024 function terms.

However the naked function is only marginally similar to its analytic continuation as shown in fig 21 (left) which has zeros on the line $x = ¼$ and singularities at the position of the zeta zeros, which are barely expressed at 1024 naked function terms. Even more pertinently, the two other derived functions sigma and mu have naked approximations with a very low degree of convergence to their analytic continuations.

Sigma $\sum_{n=1}^{\infty} \frac{\sigma(n)}{n^s} = \zeta(s)\zeta(s-1)$, $\sigma(n) = \sum_{d|n} d$ is defined in terms of the divisor function and is equivalent to a product of zetas. As shown centre in fig 21, its features and double zeros are barely apparent at 1024 terms.

Finally (right) we have mu, $\sum_{n=1}^{\infty} \frac{\mu(n)}{n^s}$, $\mu(n) = \begin{cases} (-1)^k, & n \text{ has } k \text{ distinct prime factors of multiplicity 1} \\ 0 & \text{otherwise} \end{cases}$

Which, from convolutions $\sum_{n=1}^{\infty} \frac{f(n)}{n^s} \sum_{n=1}^{\infty} \frac{g(n)}{n^s} = \sum_{n=1}^{\infty} \frac{(f*g)(n)}{n^s}$, where $(f*g)(n) = \sum_{d|n} f(d)g(n/d)$ resolves to $\sum_{n=1}^{\infty} \frac{\mu(n)}{n^s} = \frac{1}{\zeta(s)}$, since $\frac{1}{\zeta(s)} \zeta(s) = 1$, and $\mu*1 = \varepsilon$, $\varepsilon(n) = \begin{cases} 1, & n=1 \\ 0, & n>1 \end{cases}$. Consequently this has Euler product

$$\sum_{n=1}^{\infty} \frac{\mu(n)}{n^s} = \prod_{p \text{ prime}} (1-p^{-s}).$$

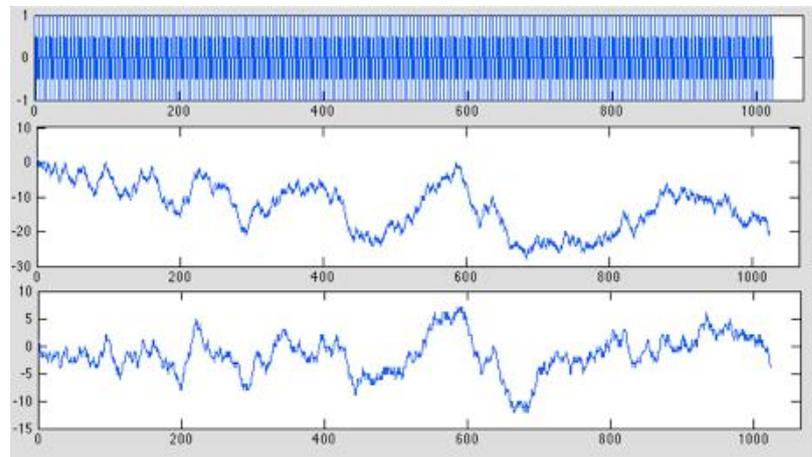

Fig 22: Trends in the added coefficients are emphasized when they are summed. L(5,3) has stable periodic sums, while lambda and mu have erratic long-term trends, which disrupt their convergence in the critical strip when calculated raw.

Fig 22 shows why convergence is bad for sigma and mu, which, unlike eta and the *L*-functions of non-trivial Dirichlet characters, illustrated by *L*(5,2) which have regularly alternating sum coefficients, have coefficients with large amounts of drift to the positive and negative, compromising their convergence. Many series generated from simple or cyclic multiplicative coefficients share such irregularities in the sum coefficients because of their varied impact on each integer through its prime factorization.

Fig 23: Trends in the multiplicative coefficients, emphasized by summing terms, are complex, even for Dirichlet *L*-functions, where the sum coefficients are strictly periodic (inset). The coefficient chain proceeds from red to blue. This is because the defining relationships depend on prime distributions modulo *n*. For example $\chi$ (4,2)={0,1,0,-1} defines $\chi$ (*p*)=1 if *p* mod 4 =1 and -1 if *p* mod 4 = 3, which also explains why the Gaussian primes of fig 9, which consist of separate *p* mod 4 = 1 and *p* mod 4 = 3 prime distributions also give rise to valid *L*-functions.

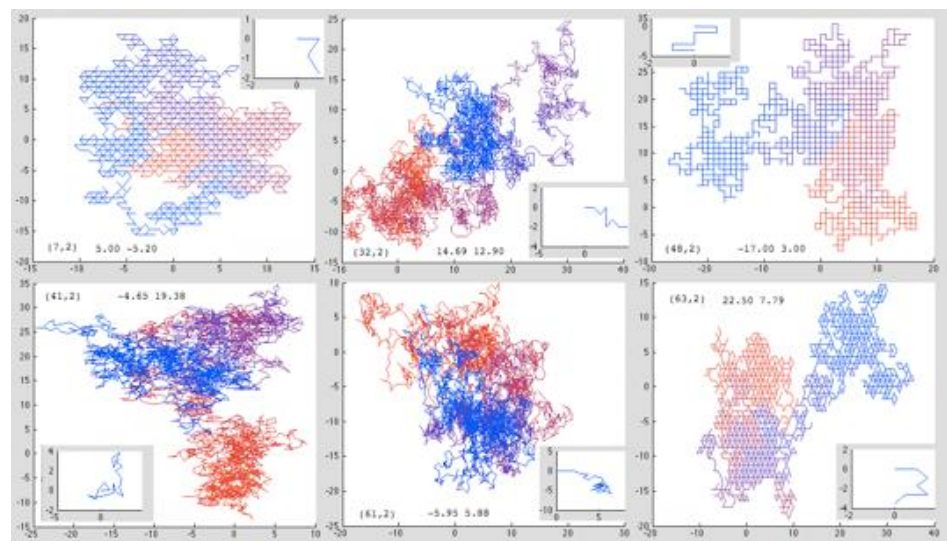

Ideally one would like an eta analogue to guarantee convergence of these functions in the critical strip, however solutions such as applying additional product terms to zeta*lambda for each +1 and -1 coefficient for real multiplicative coefficients produces spurious functions with zeros on x = 0 and x = 1/4, displaying the failure of convergence that occurs with the Euler product itself.

Alternatively, as we have seen with the Dedekind zeta and Hecke *L*-functions, we can try to represent the function in the critical strip, or even more widely, using a Mellin integral transform representation. However finding a suitable theta function can prove problematic (Garrett). Tim Dokchitser's Computel PARI script, now incorporated into Sage uses just such a sophisticated Mellin technique to explore a variety of abstract *L*-functions, but there is no implication the technique would extend to any non-*L* Dirichlet series derived from an Euler product.

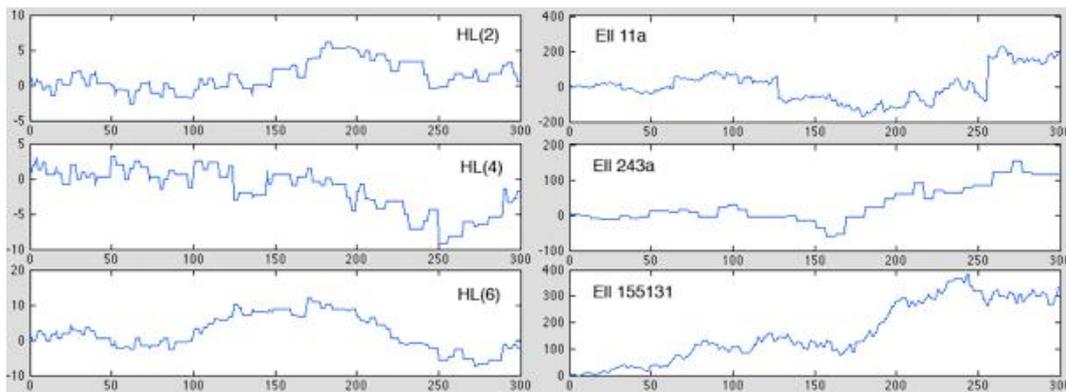

Fig 24: Trends in the product coefficients for three Hecke *L*-functions and three elliptic curve *L*-functions likewise show erratic trends through prime encoding. The Hecke coefficients are real, but have non-integer values not o absolute value 1. The elliptic curve *L*-function coefficients have varying real integer definitions depending on good or bad reduction and consequently have Dirichlet series coefficients that are not completely multiplicative.

Conversely, turning to the multiplicative coefficients, we see that it is no easy task to find criteria here which distinguish *L*-functions apparently satisfying RH from completely multiplicative functions which violate it, because the multiplicative coefficients for simple periodic Dirichlet series are encrypted through the primes into complex irregular sequences, with Brownian-like prime walks in their summed terms, as illustrated for a spread of *L*-functions in figs 23 and 24.

We are beginning to see why there is a 'conspiracy' among *L*-functions, as pointed out by Brian Conrey (2003). Essentially the *L*-functions show us arcane forms of Dirichlet series coefficients, which provide additional encrypted keys to the master lock of the Riemann zeta function's zeros on the critical line, but apart from confirming that equi-distributed primal encoding provides additional keys, emphasizing the primal basis for Riemann zeta's non-trivial zeros, they don't shed closer light on RH itself.

The fundamental difficulty here is that outside the known *L*-functions, with both Euler product representations and functional equations, thus also defining Mellin transforms in the critical strip, we have no guarantee of a convergence for the Dirichlet series defined by arbitrary Euler products, so cannot be sure we have a representation in the critical strip, let alone zeros on the critical line, which some obvious examples, such as lambda manifestly do not have.

At the core of the problem of both the Dirichlet series and Euler product coefficients is that, for each there is a countable infinity of them determining the locations of a countable infinity of zeros. In fact, we have a simple countable set of linear equations $\sum_{n=0}^{\infty} a_n n^{-z_k} = 0$ to solve for each zero $z_k$. In the countably infinite dimensional space $\{a_1, a_2, a_3, \ldots\}$ there will be a parametric transformation of the coefficients which moves only one zero, so that a continuous path can be defined in the space in which the zeros remain critical, resulting in a continuous transformation of *L*-functions. Given that two Dirichlet *L*-functions have a closely identical distribution of zeros, see fig 25, one may be able to find a path in coefficient space between them if no two zeros have to collide. There are additional topological constraints determining a discrete and varying number of gamma factors in the functional equations of individual *L*-functions, which cannot be continuously transformed. We thus now investigate how some elementary continuous transformations of Dirichlet *L*-functions affect the topological relationships of the zeros.

**Dynamically Manipulating the Non-trivial Zeros in and out of the 'Forbidden Zone'**
To get a closer view of the dynamics of the zeros we now investigate breaking out of the boundaries imposed both by the *L*-functions and the taboos created by the assumption that Dirichlet functions can be depicted for negative real parts only if they have a formal functional equation. We will thus continue to 'unashamedly' use a 'naked' depiction of Dirichlet series in the 'forbidden zone', particularly those closely approximating Dirichlet *L*, or *DL*-functions containing 'alternating' or 'rotating' coefficients, without exerting any functional equation, or expecting complete multiplicativity, so that we can see the dynamics of how zeros of such functions change under a continuous transformation between *L*-functions. This is justified by a theorem (Balanzario and Sanchez-Ortiz), which states that for a small enough change in the coefficients, there is a correspondingly small change in the locations of the zeros.

Since we are actually using finite approximations (see also Borwein et. al.), the functions we visualize will all be tame and we will only consider 'convergence' in the sense of how the assumed infinite limit might behave. For relatively small imaginary values in the tens to hundreds, no more than a few hundred terms of the Dirichlet series are needed to get a good approximation, which clearly show zeros on the line for the naked equivalents of these *L*-functions because of their convergence in the critical strip.

We will view how both changes in the character cast of *L*-functions and changes in which the usual integer values of the base exponents may be continuously shifted to adjacent real and complex values, to enable a continuous transformation between differing integer values. Rather than attempting a one-process solution to the dynamics, our aim is to explore emergent features of the dynamics of zeros under such continuous transformations of Dirichlet series, as a clue to the hidden complexity, which cannot be seen when the zeros are ostensibly fixed on the critical line and attempts are made to find abstract criteria which define those having only critical zeros aiming at an abstract proof of RH. We thus explore four examples, using richly different types of tame and wild continuous transformation. These are much better viewed as movies from the link below, but here a series of stills with path diagrams will have to suffice.

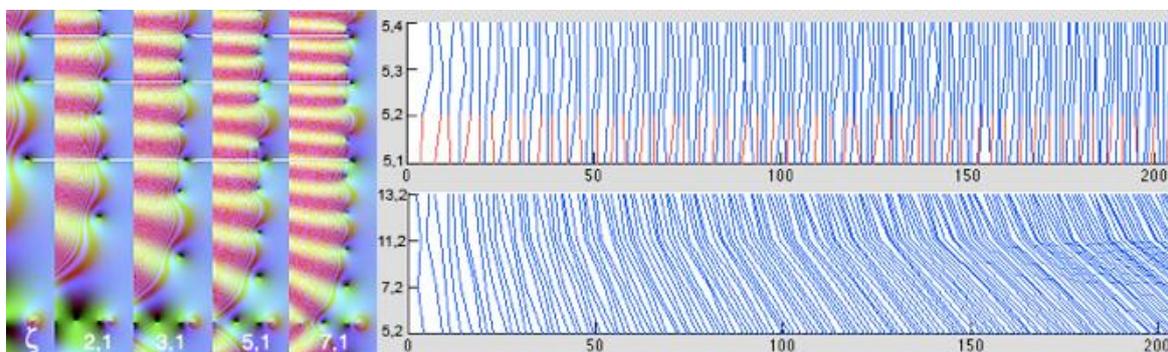

Fig 25: (Left) Zeta and a series of principal *DL*-functions showing retention of the zeta zeros with additional regular zeros on *x*=0. (Lower right) Trends in the non-trivial zeros for a series of prime non-principal characters, showing consistent trends. (Upper right) The four characters of period 5 showing the periodic zeros on *x*=0 for 5,1 occupy corresponding positions to the non-trivial zeros of 5,2, 5,3 and 5,4 indicating the *L*-functions do not distinguish between periodic and non-trivial in forming their Fourier transforms of the prime distribution modulo 5.

To gain a view of the justification of this process, fig 25 shows a series of *DL*-functions along with the distributions of their unreal zeros. Note that the principal *DL*s have the same non-trivial zeros as zeta, being derived from it by additional prime product factors, e.g. $L(5,1) = (1 - 5^{-z})\zeta(z)$, with additional periodic zeros on *x*=0. By contrast, the non-principal prime *DL*s have only non-trivial zeros on *x*=1/2. As the prime numbers increase, the non-trivial zeros become more densely distributed, but the distribution for non-principal *DL*s for a given prime coincides with the total distribution periodic and zeta of the principal *DL*s indicating the functions are treating all the zeros alike as part of an effective Fourier transform of the respective modular prime distributions.

This means that many *DL*s do not have all their effective unreal zeros on the critical line, but have a neatly phase-locked distribution with one periodic set sequestered from the irregular distribution on *x*=1/2. This kind of distribution is also shared by all *DL*s of non-prime period with multiple factors. We will examine continuous transformations between *DL*s where the intermediate stages are reasonably convergent series to elucidate the transitional dynamics.

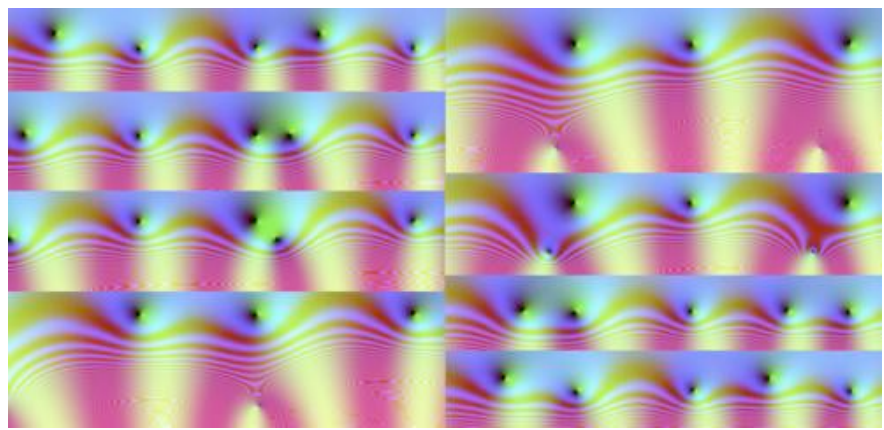

Fig 26: A continuous rotation from eta (top left) through zeta (between bottom left and top right) and back to eta (bottom right) shows the periodic zeros crossing the critical line and plunging asymptotically into the negative real half-plane as we cross zeta, subsequently being picked up by rising zeros spaced between the originals.

The first and tamest example is making a simple rotation between zeta and eta by using the multiplicative function connecting them: $f(z,\theta) = (1 - (1 + e^{i\theta})2^{-z})$, $\theta = 0, \cdots, 2\pi$. One can immediately see this will have periodic zeros at $z = \ln(1 + e^{i\theta})/\ln 2$ since
$1 - (1 + e^{i\theta})2^{-z} = 0 \Leftrightarrow (1 + e^{i\theta})2^{-z} = 1 \Leftrightarrow (1 + e^{i\theta}) = 2^z = e^{z\ln 2} \Leftrightarrow z\ln 2 = \ln(1 + e^{i\theta})$, and that $x \to -\infty$, as $\theta \to \pi$, so that the periodic zeros on the line $x = 1$ in eta will plunge into the negative real half-plane as we cross zeta. This is confirmed in fig 26, where we are able to use the functional equation throughout.

We now need to examine continuous transformations of functions that cannot be represented for negative real values using a functional equation, so we need to understand the consequences of using naked Dirichlet series functions in the forbidden zone. With one transitional exception in the last case, the examples are broadly confined to alternating series which are well-defined and convergent in the critical strip $0 < x < 1$ however for $x < 0$ these can become singular as the number of terms in the series increases in the limit to infinity.

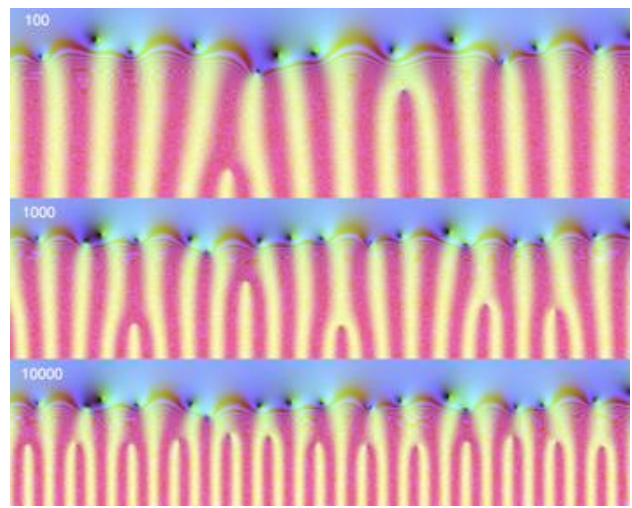

Fig 27: The non-*L*-function with $\chi$ = {0,1,0, -1,0,0,0,1,0,-1} in the neighbourhood of ½+230*i* represented naked with 100, 1000 and 10000 series terms, shows increasingly extreme fluctuations in the forbidden zone with increasing numbers of 'gollum' zeros falling closer and closer to $x = 0$. However the approximation to the zeros in the critical strip $0 < x < 1$ is sufficiently close by 1000 that little subsequent change is observable by including 10000 terms.

Fig 27 gives a portrait of a section of the function in fig 1 when the iterations are increased through two orders of magnitude, and as can be seen, there is increasing instability in the forbidden zone with increasing numbers of 'gollum' zeros closer and closer to $x = 0$.

However, while the approximation is inaccurate for 100 terms at this imaginary range, by 1000 terms, increasing the terms to 10000 has little effect on the zeros in the critical strip, showing that a finite approximation suffices, as a numerical analytic tool, if the number of terms is over a suitable bound, which varies with the imaginary value of the neighbourhood being investigated.

For the next example, we explore a continuous transformation between $L(6,2)$ with character cast $\chi = \{0,1,0,0,0,-1\}$ and $\chi = \{0,1,0,0,-1,0\}$, which corresponds to the alternating Dirichlet series having the arithmetic progression $\sum_{n=0}^{\infty}(-1)^n(1+3(n-1))^{-z}$. The corresponding series for 2 is $L(4,2)$ with character series $\chi = \{0,1,0,-1\}$, all of whose non-trivial zeros are on the critical line, but this is not the case for the above series, as shown in fig 5. In fact, this is the eta version of the Hurwitz zeta function $\zeta(z,1/3)$, which we have seen is a sum of the period 3 Dirichlet series, and does not have zeros on the critical line.

To make a continuous transformation requires moving off the natural numbers as bases. We will move continuously around a semicircular loop running firstly down the real axis from $6n+5$ to $6n+4$ and then anti-clockwise around $6n + \frac{1}{2}(9 + e^{i\theta})$, $\theta = \pi, \cdots, 2\pi$.

The first part of the trajectory gives us a good idea of why the Riemann hypothesis might be true due to mode-locking, as the zeroes each follow local orbits approximating rotations, under the continuous transformation, so the zeros which are lined non-periodically on the critical line and periodically on $x = 0$ lose their phase relationships when we move from $6n+5$ to $6n+4$, thus throwing the critical zeros 'offline'.

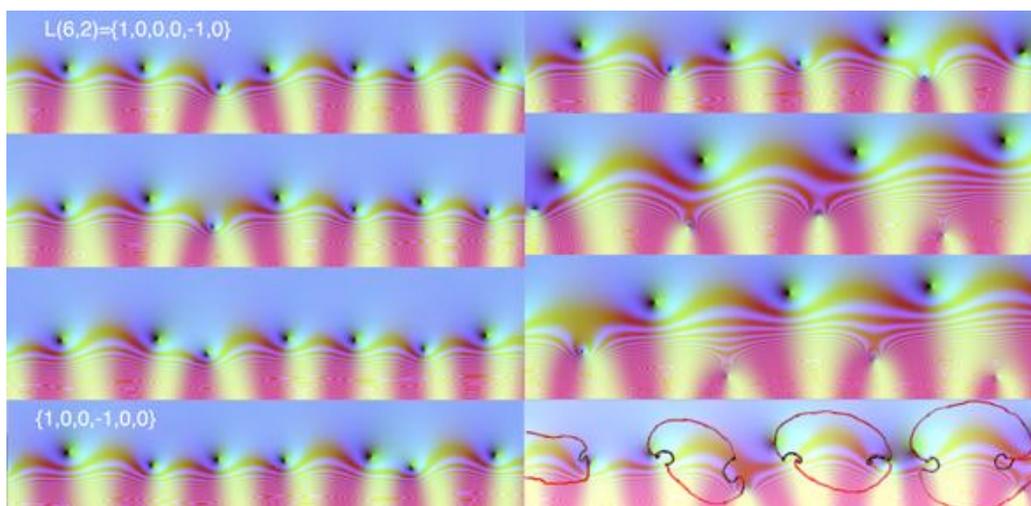

Fig 28: Continuous transformation between $L(6,2)$ and the alternating arithmetic Dirichlet series with base $(1+3(n-1))^{-z}$ left running down the real axis and right anti-clockwise around a semi-circle.
The bottom right image shows the orbits along the real line in black and round the semi-circle in red.

However this neat picture is confounded by the dynamics we perceive on the semi-circular track, where pairs of zeros exchange places, plunging both deep into the 'forbidden' zone, and well up into the positive half-plane. Significantly the neat distinction between the critical non-periodic zeros and those on $x = 0$, which satisfies the generalized Riemann hypothesis, that if a zero is in $0 < x < 1$ then it is on $x = \frac{1}{2}$, is also lost, because critical and periodic zeros are interchanged. The dynamics of the zeros on the far left and right are not fully elucidated and may be periodic, or otherwise, as we shall see in following examples.

Note that the movement of the zeros is path-dependent and that continuous transformations can exchange the roles of critical, periodic forbidden 'gollum' zeros.

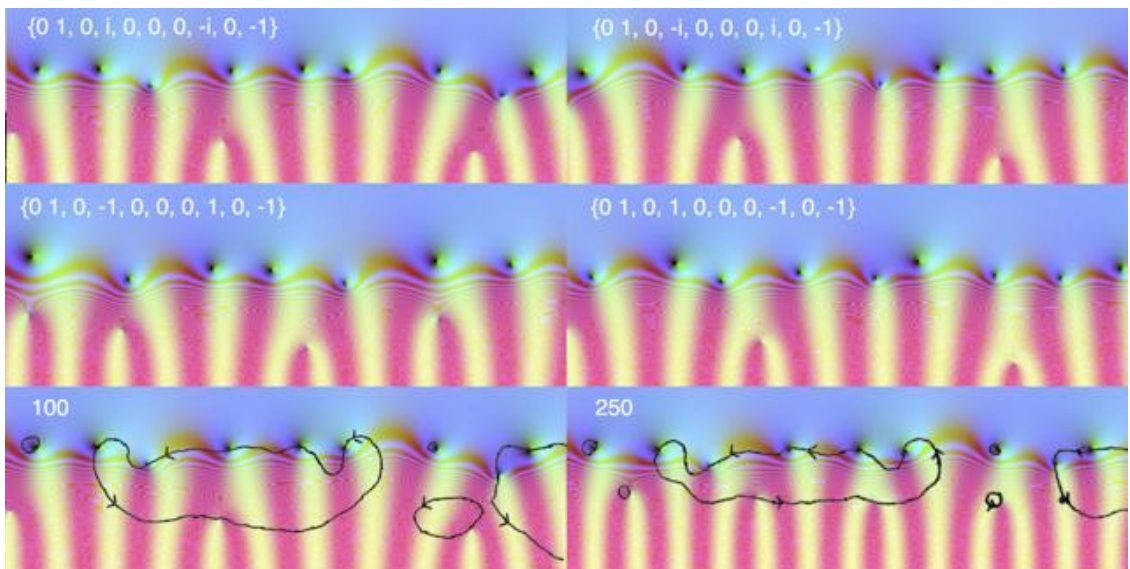

Fig 29: Cyclic rotation of the characters from $L(10,4)$ to $L(10,2)$ and back demonstrates 'transmigration' of the zeros, one step to the left each character cycle, in which critical, periodic and gollum zeros exchange positions. To assess the validity of approximation, using naked functions in the forbidden zone, the bottom two images compare the orbits for 100 and 250 iterations in a neighbourhood of $\frac{1}{2} + 24i$. The central transmigration orbit is preserved and an additional gollum carrier zero has entered the loop.

The next example takes this further into the wilderness, by examining changes in the character terms rather than the positions of the integer bases. We start with the character cast for $L(10,4)$ $\chi = \{0,1,0,i,0,0,0,-i,0,-1\}$ and apply $\chi(\theta) = \{0,1,0,e^{i\theta},0,0,0,-e^{i\theta},0,-1\}$, $\theta = \pi/2, \cdots, 5\pi/2$, a cyclic rotation, passing through $L(10,2)$ $\chi = \{0,1,0,-i,0,0,0,i,0,-1\}$ at $\theta = 3\pi/2$ and the non $L$-function $\chi = \{0,1,0,-1,0,0,0,1,0,-1\}$, at $\theta = \pi$.

As we move around the cycle, two of the critical zeros remain in small local closed orbits, but the rest, including both critical and periodic zeros, pass in a chain, from one to another, stepping once to the left for each complete cycle of rotation, exchanging places on the way with one of the 'gollum' zeros which should not exist in the negative real half-plane. In all it takes 10 cycles of the characters for a zero to move across the field of view.

This shows us that the dynamics of the critical zeros under continuous transformation of the function cannot be understood without taking into account the 'gollum' zeros in the naked representation that are eliminated in the functional equation representations of zeta and the $L$-functions, which have legitimate zeros only on $y = 0$, $x = 0$ and $x = \frac{1}{2}$.

Fig 30: Braiding in one cycle through $L(10,4)$.

The braiding of the zeros on cycling the characters shows us that there is a topological reason why the zeros have to move off the line, as we pass along this particular parametric loop in function space from one $DL$-function to another. There is no way a set of crtical zeros can make a braided transformation topologically on the critical line without either passing off it, or colliding to form transition states with degenerate multiplicities, which will lead to a pair of off-line zeros emerging from the collision.

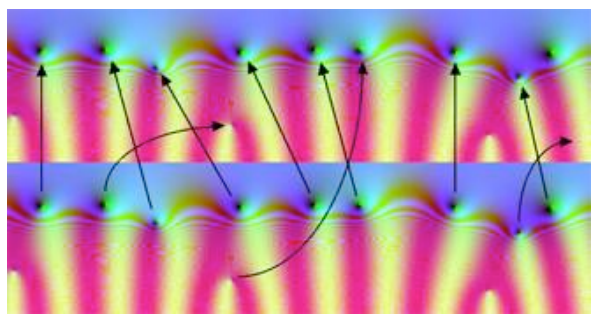

To address the problem of the naked finite functions in the forbidden zone not necessarily converging to the limit in fig 27 are orbits for two different numbers of terms, 100 and 250, in a neighbourhood of ½+24$i$. Intriguingly, although the increased number of terms has given rise to a greater number of gollums falling closer to the imaginary axis, the phase portrait of the orbits retains homology. The major transmigration cycle of the critical strip zeros is preserved by utilizing additional gollums as carriers. Other gollum and critical zeros have local closed orbits in both cases.

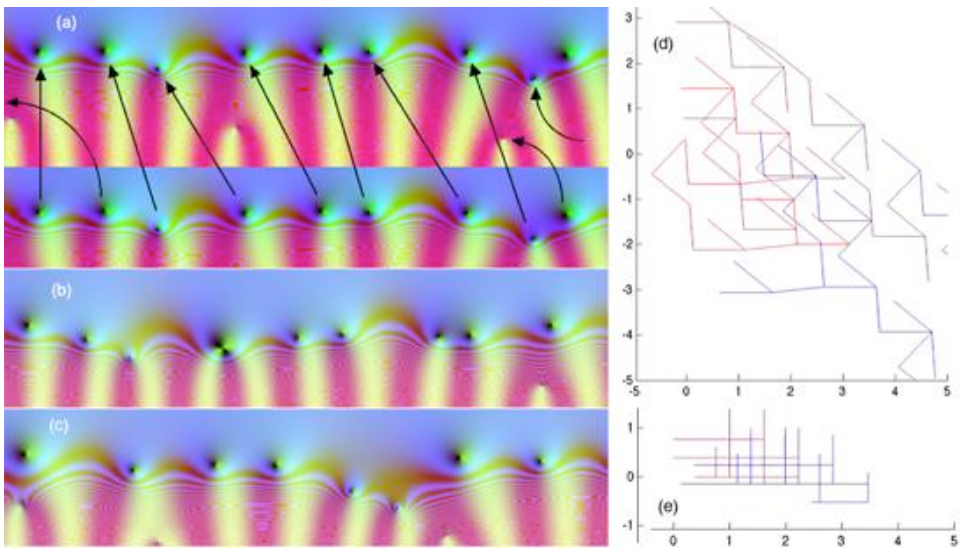

Fig 31: The Euler product may not be sufficient. (a) Braiding in the rotational parametrization also has off-critical intermediates, with near collisions (b) and involvement of gollum zeros. Although this could be due to convergence problems from the erratic coefficient sum (d), even the tame linear translation (c) has severe off critical states although its additive coefficient sum (e) supports convergence in the critical strip.

However it is clear also that the intermediate states have coefficients which are not completely multiplicative, so one could seek a set of intermediates which still possess an Euler product as a test of whether the existence of an Euler product is sufficient for the non-trivial zeros to be critical. We can do this by rotating the characters of each prime and then defining series coefficients in terms of their prime factorizations: $f(z,\theta) = \sum_{n=1}^{\infty} a_n n^{-z}$, $a_n = \prod_{p_i} \chi(p_i,\theta)$, $\chi(\theta) = \{0,1,0,\varphi(\theta),0,0,0,-\varphi(\theta),0,-1\}$.

This will have another offshoot problem of producing non-periodic Dirichlet series, which may have convergence problems because of wandering in their signs and angles.

However, when we do this, we find that the Euler product is not sufficient, and indeed cannot be for the same topological reasons outlined above, namely that for this loop also there is a topological braiding of the zeros. For the rotation $\varphi(\theta) = e^{i\theta}$, $\theta \in [\pi/2, 5\pi/2]$ as in fig 31, this also results in off-critical intermediate states. Even the tame linear translation $\varphi(\theta) = \theta$, $\theta \in [-1,1]$, we find severe off-critical intermediate states.

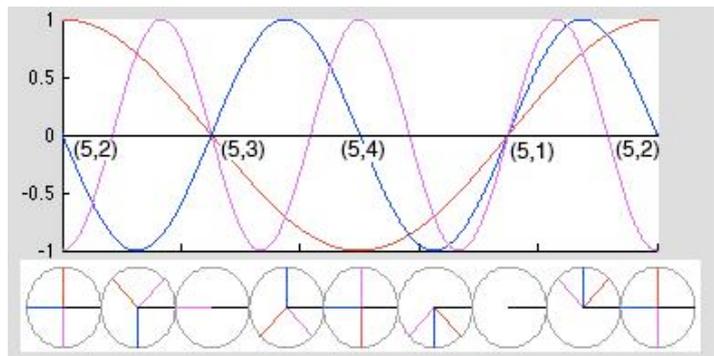

Fig 32: Character cast rotation of $L(5)$ passes through all four $DL$-functions albeit with naked convergence issues for $L(5,1)$ and its neighbouring functions. Above is the trends in the imaginary part and below the distribution of the individual characters.

To take these examples to a fireworks finale, we have an example of another rotation of a character cast, this time of $L(5)$ rotating each character by the factor implied by their position on the unit circle $L(5,\theta) = \{0,1,e^{i\theta},e^{3i\theta},e^{2i\theta}\}, \theta = \pi/2,\cdots,5\pi/2$.

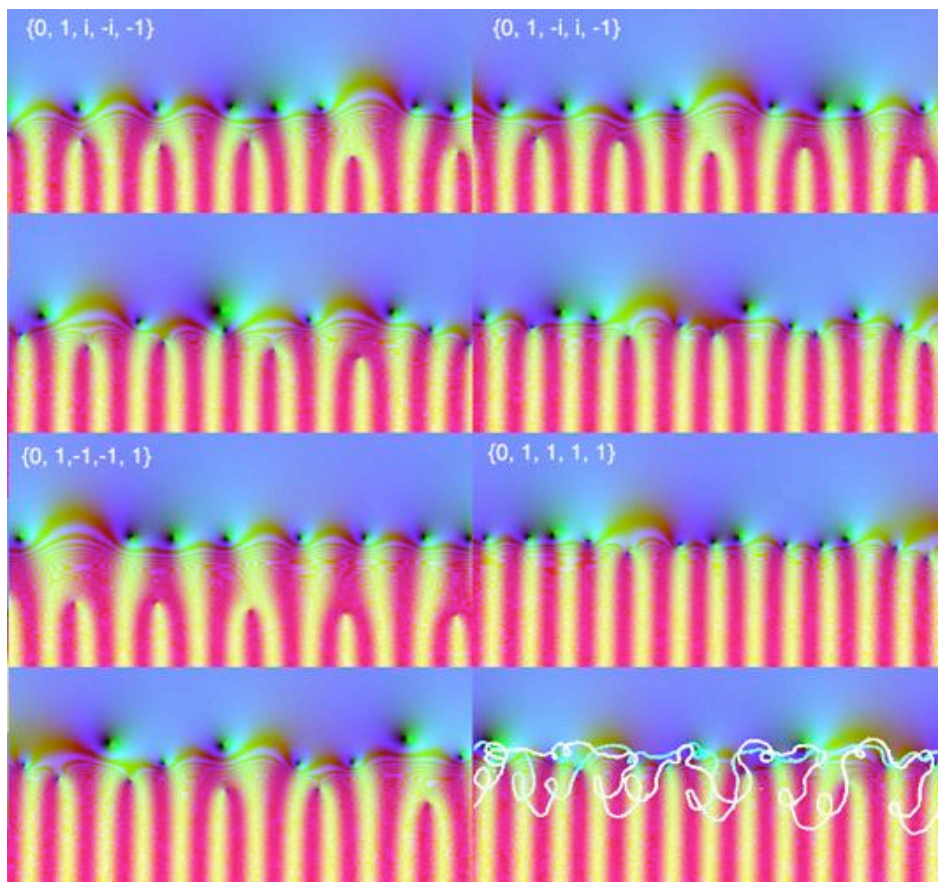

Fig 33: The character cast rotation of $L(5)$ displaying two complex entwined orbits again linking critical and 'gollum' zeros. The entire zero set moves to the left each rotation cycle.

This is pushing the boundaries, as it involves potentially 'non-covergent' intermediates, but it passes through all four $L$-functions $L(5,2) = \{0,1,-i,i,-1\}$ at $\pi/2$, $L(5,3) = \{0,1,-1,-1,1\}$ at $\pi$, $L(5,4) = \{0,1,i,-i,-1\}$ at $3\pi/2$ and $L(5,1)=\{0,1,1,1,1\}$ at $2\pi=0$, and gives an excellent example of complex orbits in motion. A diagram of the rotations is shown in fig 33 with the $L$-functions located. The orbits in this case have become very complex with many enclosed loops and 'entwined in the sense that successive zeros pass through almost identical paths before diverging again. As in the previous example the orbits involve critical zeros being carried into the position of gollum zeros and vice versa as they pass along the same orbit.

Fig 34: Fractal patterns in the Dirichlet coefficients of the Euler product encoded version of the above character rotation.

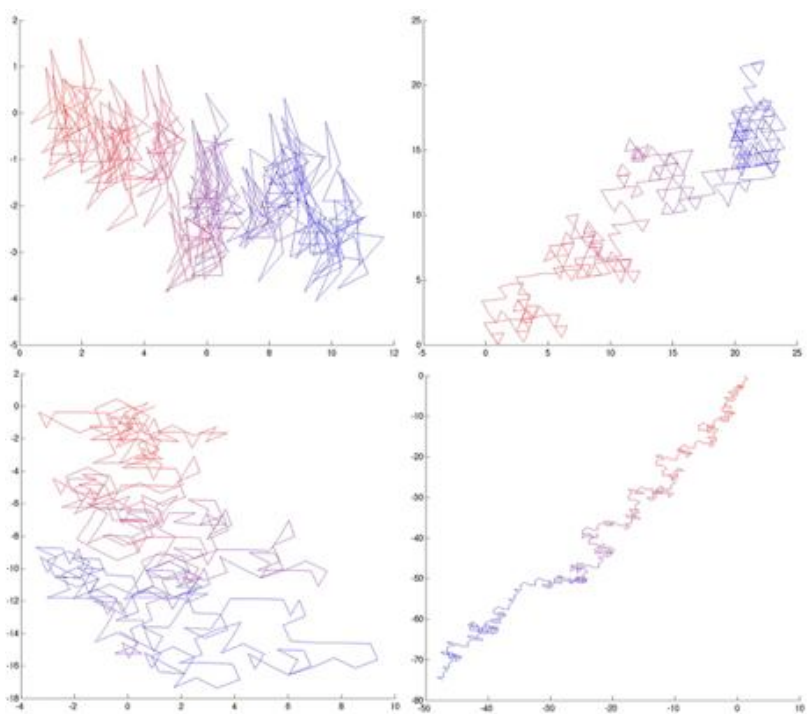

**Finding Coefficient Paths with On-Critical Zeros**

The essential problem, as we have already noted is that the paths we are taking in function space are simple translations and rotations of only a few variables of what is a countably infinite function space $\{a_1, a_2, \ldots a_n, \ldots\}$. To find out what patterns of coefficients would make a transition between $DL$-functions which do keep the non-trivial zeros on the critical line, even if some Gollum or other zeros appear as well, we reverse the problem and seek a specific solution to a finite set of sum coefficients for a series of parameter values making a transition between $L(5,2)$ and $L(5,4)$.

We seek $\sum_{n=1}^{N} a_n^\theta n^{-z_j} = 0$, $z_j^\theta = \frac{1}{2} + i\rho_j^\theta$, $j = 1, \cdots, k$, where $\theta$ is a parameter defining a transition between the $j$th zeros of $L(5,2)$ and $L(5,4)$. Solving these linear conditions for the first $n-1$ zeros with the accessory scaling constraint $\sum_{n=1}^{N} a_n^\theta = 1$, we have:

$$\begin{pmatrix} 1^{z_1^\theta} & \cdots & n^{z_1^\theta} \\ \vdots & \ddots & \vdots \\ 1^{z_{n-1}^\theta} & \cdots & n^{z_{n-1}^\theta} \\ 1 & \cdots & 1 \end{pmatrix} \begin{pmatrix} a_1^\theta \\ \vdots \\ a_n^\theta \end{pmatrix} = \begin{pmatrix} 0 \\ \vdots \\ 0 \\ 1 \end{pmatrix}, \text{ with solution } \begin{pmatrix} a_1^\theta \\ \vdots \\ a_n^\theta \end{pmatrix} = \begin{pmatrix} 1^{z_1^\theta} & \cdots & n^{z_1^\theta} \\ \vdots & \ddots & \vdots \\ 1^{z_{n-1}^\theta} & \cdots & n^{z_{n-1}^\theta} \\ 1 & \cdots & 1 \end{pmatrix}^{-1} \begin{pmatrix} 0 \\ \vdots \\ 0 \\ 1 \end{pmatrix}$$

Rotating the solutions so that $a_1^\theta = 1$ gives a unique conjugate pair of evolving series coefficients, with critical zeros varying linearly between those of $L(5,2)$ and $L(5,4)$.

Running a solution in Matlab gave reasonable matrix inversions up to the order of 30-40 zeros. The series of coefficients begin with a signature almost identical to $\chi(5,2) = \{0,1,-i,i,-1,0\}$:
{1, 0.07-1.00i, 0.0+1.0i, -0.99-0.07i, 0.01+0.01i, 0.98+0.07i, -0.06-1.00i, -0.19+0.95i, -1.15- 0.07i, -0.24 -0.14i, ...}
and its conjugate to $\chi(5,4)$, however they rise to an exponential hump towards the end, indicating why larger matrix solutions cause overflow. The intermediate sates display an intriguing feature suggesting evidence of the primes even at this level of approximation. The initial and end states have zeros principally on the critical line coinciding with those of $L(5,2)$ and $L(5,4)$ but intermediate states wih zeros which are critical but do not conform to a prime distribution display a dominant pattern of additional zeros emerging from the forbidden zone into the positive half-plane.

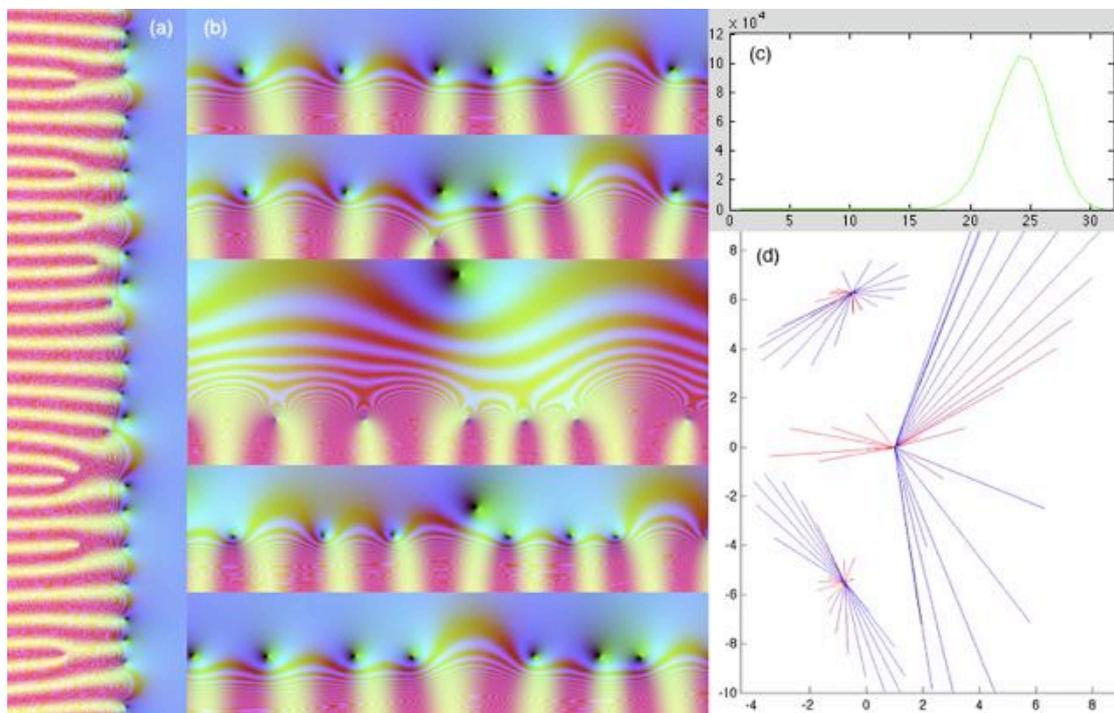

Fig 35: (a) Finite Dirichlet series with the same first 30 critical zeros as $DL(5,2)$. (b) Transition to $DL(5,4)$ preserves criticality of the moving zeros but a new series of zeros make a circuit into the positive half plane. (c) Absolute coefficient trend, (d) Complex coefficient maps for three states beginning, intermediate and end.

**Conclusion**

The Riemann hypothesis cannot be fully understood without decoding how the interference of the imaginary logarithmic wave functions results in the distribution of the zeros. The confinement of the zeta and *L*-function zeros to the critical line is clearly a result of the primes being asymptotically as close to evenly distributed in relation to the logarithmic integral as they can possibly be, given that they cannot be evenly distributed and be prime. This is already established as a necessary and sufficient condition for RH and the main reason the proof of RH has been sought as a means to establish the asymptotic prime distribution is because of the 'fatal attraction' of the symmetry of the Riemann zeta and xi functions. However the functional equation which is the basis of this symmpoetry also permits symmetrical off-critical zero pairs and we have seen these occurring in various counterexamples with functional equations.

The portraits of convergence of the Dirichlet sums imply that the Riemann hypothesis is a consequence of minimal phase-locking in the imaginary wave functions at $x=1/2$, caused by the prime distribution, making convergence to zero possible in the asymptotic limit for cis($y$ln$n$) at the one 'index value' of $n^{-1/2}$, determining the power law trend in the absolute values of the terms.

The diverse types of abstract *L*-function demonstrate the involvement of primes in equidistributed ways through their Euler products in other coefficient series also having critical zeros but with differing imaginary values and real weights determining the critical line. These however do not necessarily aid the proof of RH because they represent encoded forms of the same asymptotic prime distribution, forming a regress into more and more rarified encodings of the same prime distribution. In this sense there is a conspiracy among them, in the form of Chinese whispers echoing RH upon itself, preventing them shedding new light on the hypothesis.

While this conspiracy remains discrete when dealing with Dirichlet *L*-functions and those of elliptic curves, when we come to examine modular forms, we find we have continuous spaces of functions, providing a way to distinguish forms that might satisfy RH from those which do not.

The modular eigenforms show why only some functions among a continuous distribution of functions admit both a functional equation and have an Euler product and also appear to satisfy RH, while neither of these on their own appear to be sufficient. This casts the *L*-function conspiracy problem into a quantum form in which the *L*-functions like the solutions of the Schrödinger wave equation have perfect re-entrant properties and thus only these would satisfy RH.

RH may thus be true, but unprovable, as a type of Turing halting problem, because the, despite the apparent symmetry of the zeta zeros, RH can be confirmed only by infinite computation as in the much simpler Collatz conjecture. If this is so, the zeta zeros, asymptotic prime distribution and that of the Farey fractions are logically equivalent definitions, but no one can be proved to establish the truth of the others, thus giving all three a similar status to the axiom of choice, as an additional number postulate about asymptotic infinities. The root fact is thus the prime distribution, from which the distribution of the zeta zeros, and Farey fractions follows, rather than the proof of the prime distribution following from RH.

**Appendix 1: Mediants and Mode-Locking**

In dynamical mode-locking, any irrational rotation close enough to a rational fraction of a revolution undergoing dynamical feedback becomes locked to the periodicity of that rotation, forming a series of intervals of mode-locked states, with a residual set of points in between retaining their unperturbed irrational motion. These mode locked intervals form a continuous fractal monotone increasing function, constant on intervals surrounding each rational number, called the Devil's staircase illustrated in fig 5. Mode-locking is manifest in many processes where dynamic periodicities interact, including the non-mode-locked orbits (to Jupiter) of the remaining asteroids, because the mode-locked ones were thrown into chaotic orbits and collided with planets, the ordered mode-locking of the Moon's day to the month, and Mercury's day to 2/3 of its year.

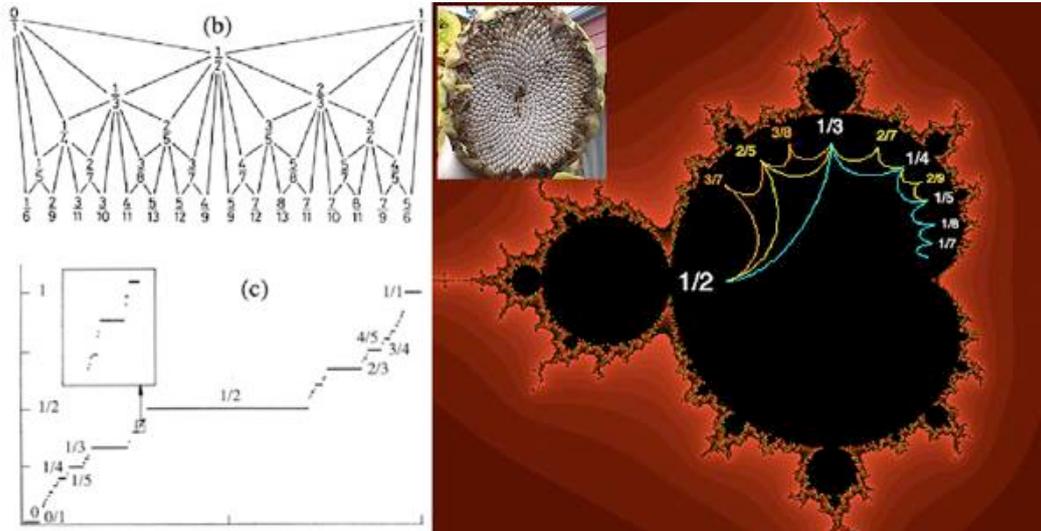

Fig 36: Left Farey Tree and Devils Staircase. Right: Mediant-based mode-locking in the Mandelbrot set bulbs, defining their fractional rotation, the periodicity in each bulb and the number of its dendrites - e.g. 1/2 and 1/3 span 2/5. Spirals in the sunflower follow Fibonacci numbers, minimizing mode-locking by approximating the golden angle $\gamma\pi$.

As shown in fig 36, the bulbs on the Mandelbrot set follow the fractions on the Farey tree, adding fractions as mediants $\frac{p}{q} + \frac{r}{s} = \frac{p+r}{q+s}$. This can be seen by counting the number of their dendrites, which also corresponds to the periodicity of the attractor in each bulb. Mediants correctly order the fractional rotations between 0 and 1 into an ascending Farey sequence, providing a way of finding the fraction with smallest denominator between any two other fractions. A way of seeing why this is so is provided by using a discrete process to represent the periodicities or fractional rotations. For example if we have 2/3 = [110] and combine it with ½ = [10] by alternating, we get [11010], or 3/5. The Golden Mean $\gamma = \frac{-1 \pm \sqrt{5}}{2} = 0.618, -1.618$, by virtue of its defining relation $\frac{1}{\gamma} = 1 + \gamma$ is the limit of the ratios of successive Fibonacci numbers 1, 1, 2, 3, 5, 8, 13, 21 … for which $F_{n+1} = F_n + F_{n-1}, F_0 = F_1 = 1$, is non-mode-locked. The Farey tree leads to other Golden numbers, if we alternate left and right as we descend, following a series of Fibonacci fractions. Numbers g, such as these, avoid becoming mode-locked because their distance from any fraction of a given denominator q exceeds a certain bound: $\left| g - \frac{p}{q} \right| > \frac{\varepsilon}{q^2}, \varepsilon \sim \frac{1}{\sqrt{5}}$.

The golden numbers can most easily be described in terms of continued fractions, which, when truncated represent the closest approximation by rationals: $n = a_0 + \cfrac{1}{a_1 + \cfrac{1}{a_2 + \cfrac{1}{\ldots}}} = [a_0, a_1, a_2, \ldots]$.

Golden numbers end in a series of 1's thus having slower convergence to fractions of a given denominator than any other numbers. The Golden Mean itself is [1,1,1,…]. More generally, the Farey Tree has straightforward natural rules of parental and descendent inheritance, using continued fractions – e.g. 2/5 = [2,2]=[2,1,1] has descendents 3/7=[2,3] and 3/8 = [2,1,2] each gained by adding 1 to the last term in the two equivalent formulations. Any fraction or quadratic irrational has an eventually repeating continued fraction – e.g. $1/\sqrt{3} = [1,\overline{1.2}]$.

**Appendix 2: Finite Fields and Square Roots of -1**

Since $F_p$ is just the field of integers mod $p$, we can calculate the squares of each reside to determine how many roots of unity each has, for example: $F_2=\{0, 1\}$ with $1=1^2 = -1$
$F_3=\{0, 1, 2\}$ with $2 = -1$ and the squares $\{0, 1, 4\} = \{0, 1, 1\}$ do not contain any 2's.
$F_5 =\{0, 1, 2, 3, 4\}$ gives $\{0, 1, 4, 9, 16\} = \{0, 1, 4, 4, 1\} = \{0, 1, -1, -1, 1\}$ two square roots of 4.
$F_7 =\{0, 1, 2, 3, 4, 5, 6\}$ gives $\{0, 1, 4, 9, 16, 25, 36\} = \{0, 1, 4, 2, 2, 4, 1\}$ do not contain any 6's.

However, $F_{p^m}$ needs to be defined using an irreducible polynomial. Taking $F_9$ as a key example we want to check there are 2 square roots of -1. We need to find a degree 2 ($9=3^2$) polynomial $f(x)$ which is irreducible in $F_3$ and look at $F_3/f(x)$. Examining all $f(x) = x^2 + ax + b$, $a,b \in F_3$, we can confirm $x^2 +1$, $x^2 + x + 2$, $x^2 + 2x + 2$ have no $F_3$ zeros over $\{0, 1, 2\}$, so we can take $\alpha : \alpha^2 +1 = 0$ in $F_3$ attached to $F_3$ itself viz $\{0,1,2,\alpha,1+\alpha,2+\alpha,2\alpha,1+2\alpha,2+2\alpha\}$. Taking squares here, we have: $\{0,1,2,\alpha^2,1+2\alpha+\alpha^2,4+4\alpha+\alpha^2,4\alpha^2,1+4\alpha+4\alpha^2,4+8\alpha+4\alpha^2\}$
$= \{0,1,2,\alpha^2,2\alpha,3+4\alpha,4\alpha^2,1+4\alpha-4,8\alpha\} = \{0,1,2,-1,2\alpha,\alpha,-1,\alpha,2\alpha\}$
So we do indeed get 2 square roots of -1!

**Appendix 3: Derivation of Davenport Heilbronn**

Consider the period 5 quasi-character $\chi_\xi = \{0,1,\xi,-\xi,-1\}$ (Bombieri and Ghosh, Titchmarsh). For $\xi = i$, we can use this to generate each of the period 5 Dirichlet characters term by term by $\chi_k(n) = (\chi_\xi(n))^k$, $k = 0,\cdots,4$. For any completely multiplicative function $\psi$, we can then define $L_{k,\psi}(z) = \sum_{n=1}^{\infty} \chi_k(n)\psi(n)n^{-z}$ and set $f_{\xi,\psi}(z) = \sum_{n=1}^{\infty} \chi_\xi(n)\psi(n)n^{-z} = \frac{1}{2}\big((1-i\xi)L_{1,\psi}(z)+(1+i\xi)L_{3,\psi}(z)\big)$, where $\chi_1 = \{0,1,i,-i,-1\}$, $\chi_3 = \overline{\chi}_1$, since the other two characters are symmetric and cancel out.

For real $\xi$, we can write the above as $f_{\xi,\psi}(z) = \frac{1}{2}\sec(\theta)\big(L_{\chi,\psi}(z)e^{-i\theta} + L_{\overline{\chi},\psi}(z)e^{-i\theta}\big)$, $\chi = \chi_1$, where $\xi = \tan(\theta)$, giving quasi-character $\chi(\theta) = \{0,1,\tan(\theta),-\tan(\theta),-1\}$. Setting $\xi_\pm = -\phi \pm \sqrt{1+\phi^2}$, $\phi = (1+\sqrt{5})/2$, we then get two solutions with $\xi_- = 1/\xi_+$ satisfying the functional equation $\frac{\pi^{-z/2}}{5}\Gamma\left(\frac{1+z}{2}\right)f_\xi(z) = \frac{\pi^{-(1-z)/2}}{5}\Gamma\left(\frac{1+(1-z)}{2}\right)f_\xi(1-z)$.

**Appendix 4: A Comparison of Computational Methods**
An analysis of the various computer methods for depicting *L*-functions is revealing of the strengths and weaknesses of each. Tim Dokchitser's Computel algorithm incorporated into Sage, and now included in RZViewer as an adjunct package, using advanced generalized Mellin transforms, is highly accurate at depicting the zeros at the centre and in the critical strip up to ~1±35*i* but then undergoes a catastrophic transition, losing the zeros entirely, thus displaying lack of convergence for large imaginary values, as illustrated in fig 17.

The method is summarized as follows. Given a motivic *L*-function *L(f,z)* we consider

$L^*(f,z) = A^z \gamma(z) L(f,z)$, $\gamma(z) = \Gamma\left(\dfrac{z+\lambda_1}{2}\right)\cdots\Gamma\left(\dfrac{z+\lambda_d}{2}\right)$, $A = N^{1/2}/(2\pi)$ and seek $\phi(t)$ such that

$\gamma(z) = \int_0^\infty \phi(t) t^z \dfrac{dt}{t}$ i.e. $\phi(t)$ is the inverse Mellin transform of $\gamma(z)$. Then

$L^*(z) = A^z \sum\limits_{n=1}^\infty a_n \dfrac{\gamma(z)}{n^z} = \int_0^\infty \sum\limits_{n=1}^\infty a_n \phi\left(\dfrac{nt}{A}\right) t^z \dfrac{dt}{t}$. The method then uses the inverse Mellin transform to find $\phi(t)$ based on the residues of the individual gamma factors. Three separate methods, a Taylor formula for small $t$, an approximant for mid-range $t$, and an asymptotic formula for large $t$ are then used to calculate $\phi(t)$, and its generalization $G_z(t) = t^{-z} \int_t^\infty \phi(x) x^z \dfrac{dx}{x}$ used to calculate $L$ and its derivatives.

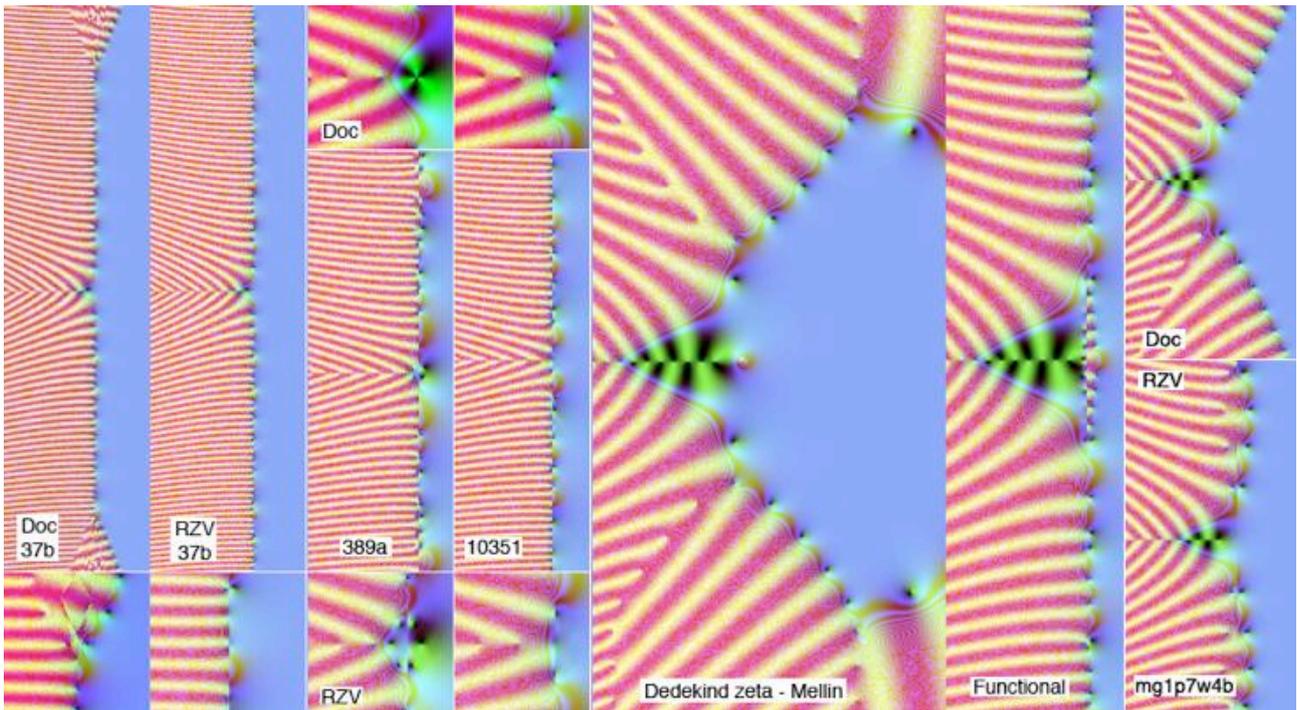

Fig 37: Overview of the computer methods. (Left pair) Dokchitser's Mellin transform method compared with applying the functional equation to the series in RZViewer. (Centre pair) Mellin and functional equation views of the central zero. (Right pair) Dedekind zeta on $Z[i]$ Author's Mellin transform method and the functional equation. (Far right) both computel and the functional equation break down dealing with a modular form with positive coefficients.

This can be compared in fig 37 with the functional equation method $L^*(f,z) = \varepsilon L^*(f, w-z)$, $|\varepsilon| = 1$, simply expressing the left half-plane in terms of the Dirichlet series $L(f,z) = \sum\limits_{n=1}^\infty a_n n^{-z}$. Here the lack of convergence is in the neighbourhood of the critical line, particularly in the neighbourhood of the origin, when the coefficients are not rapidly alternating in sign or regularly rotating in complex angle. On a dual core Intel Mac, the functional equation run in RZViewer is approximately 2,200 times faster than the Computel Mellin transform running in Sage.

The two methods are thus complementary, with the Mellin transform excellent for the neighbourhood of the origin and a bounded region of the critical strip and the functional equation and Dirichlet series good for a global profile and approximate investigation of the critical strip outside the bounds of Mellin transform convergence.

The catastrophic breakdown of the Computel method for large imaginary values appears to be characteristic of the convergence limits of Mellin transforms generally, as a very similar profile results, as show in fig 17 when a more elementary specific Mellin transform:

$$\zeta_o(z) = \pi^z \Gamma(z) \int_1^\infty (y^z + y^{1-z}) \frac{\theta(iy) - 1}{4y} dy, \quad \theta(iy) = \sum_{m,n \in Z} e^{-\pi(m^2+n^2)y} = \left( \sum_{n \in Z} e^{-\pi n^2 y} \right)^2$$ used by the author,

following Garrett (2011), for Dedekind zeta on $Z[i]$, as shown right if fig 17, compared with the functional equation: $\pi^{-z}\Gamma(z)\zeta_o(z) = \pi^{-(1-z)}\Gamma(1-z)\zeta_o(1-z)$ on $\zeta_o(z) = \frac{1}{4} \sum_{m,n \text{ not both } 0} \frac{1}{(m^2+n^2)^z}$.

## Appendix 5: Useful Sage and PARI-GP Commands

Both are accessible inside Sage 4.7 using Sage and GP terminal sessions.

### (a) Sage commands to find equations of all elliptic curves of a given conductor

```
c = CremonaDatabase()
c.allcurves(399)
{'a1': [[1, 1, 0, -210, -441], 1, 2], 'a2': [[1, 1, 0, -1925, 31458], 1, 2],
'b1': [[1, 1, 1, -13, -22], 1, 2], 'b2': [[1, 1, 1, -48, 90], 1, 2], 'c2': [[1,
0, 0, -466, 2813], 0, 2], 'c1': [[1, 0, 0, -431, 3408], 0, 2]}
```

### (b) GP script to generate elliptic curve L-function coefficients for importation to RZViewer

```
elleq    = [0, -1, 0, -651, 6228]
ellorig  = ellinit(elleq);
gred     = ellglobalred(ellorig);
ell      = ellchangecurve(ellorig,gred[2]);
conductor = gred[1];        \\ conductor for the exponential factor
gammaV   = [0,1];           \\ list of gamma-factors
weight   = 2;               \\ L(s)=sgn*L(weight-s)
sgn      = ellrootno(ell); \\ sign in the functional equation
a(k)     = ellak(ell,k);   \\ L-series coefficients a(k)
print("Elliptic curve : ", elleq);
print("Conductor      = ", conductor);
print("Root number    = ", sgn);
for(i=1,100, print(a(i),","));
```

### (c) Sage Modular form and Elliptic curve commands to generate basis and eigenfunctions

```
M=ModularForms(SL2Z,12, prec=6);
M.dimension()
2
M.basis()
[
q - 24*q^2 + 252*q^3 - 1472*q^4 + 4830*q^5 + O(q^6),
1 + 65520/691*q + 134250480/691*q^2 + 11606736960/691*q^3 +
274945048560/691*q^4 + 3199218815520/691*q^5 + O(q^6)
]
M=ModularForms(Gamma0(37),2, prec=10);
M.dimension()
3
M.basis()
[
q + q^3 - 2*q^4 - q^7 - 2*q^9 + O(q^10),
q^2 + 2*q^3 - 2*q^4 + q^5 - 3*q^6 - 4*q^9 + O(q^10),
1 + 2/3*q + 2*q^2 + 8/3*q^3 + 14/3*q^4 + 4*q^5 + 8*q^6 +
16/3*q^7 + 10*q^8 + 26/3*q^9 + O(q^10)
]
M=CuspForms(Gamma0(37),2, prec=10);
M.basis()
```

```
[
q + q^3 − 2*q^4 − q^7 − 2*q^9 + O(q^10),
q^2 + 2*q^3 − 2*q^4 + q^5 − 3*q^6 − 4*q^9 + O(q^10)
]

S = CuspForms(39)
T2 = S.hecke_matrix(2); T2
T5 = S.hecke_matrix(5); T5
CuspForms(39).newforms('a')
```

**(d) Sage script to generate modular form *L*-function coefficients for importation to RZViewer**

```
M=ModularForms(Gamma0(399),2, prec=100);
c=M.basis()
d=c[0]
for j in range (p):
f=float(d[j])
print(f)

Can also be used to print out elliptic curve coefficients

E=EllipticCurve('39a')
E = EllipticCurve([1,1,0,-4,-5]); E
d=E.q_eigenform(300)
```